%% file: nlconsensus_double_column_rev4.tex
\documentclass[10pt,final,twocolumn]{IEEEtran}
\usepackage[english]{babel}
\title{Distributed Control of Networked~Dynamical~Systems: Static~Feedback, Integral~Action and~Consensus}

\author{Martin Andreasson, Dimos V. Dimarogonas, Henrik Sandberg and Karl~H.~Johansson

\thanks{This paper was first submitted for review October 1 2012. The authors are supported in part by the European Commission by the HYCON2 project, the Swedish Research Council (VR) and the Knut and Alice Wallenberg Foundation. The second author is also affiliated with the Centre for Autonomous Systems at KTH and is supported by the VR 2009-3948 grant. 

The Authors are with the ACCESS Linnaeus Centre, KTH Royal Institute of Technology, 11428 Stockholm, Sweden.
Corresponding author: Martin Andreasson, mandreas@kth.se.  Preliminary versions of this paper were presented as \cite{Andreasson2012_acc} and \cite{Andreasson2012_cdc}. }  
}

\IEEEoverridecommandlockouts

\usepackage{mathrsfs}
\usepackage{amsfonts}
\usepackage{amssymb}
\usepackage{amsmath}
\usepackage{graphicx}
\usepackage{amsthm}
\usepackage{commath}
\usepackage{tikz}
\usepackage{tkz-graph}
\usetikzlibrary{arrows}
\usepackage{flushend}
\usepackage{psfrag}
\usepackage{pgfplots}
\usepackage{todonotes}
\usepackage{blockdiagram}
\usetikzlibrary{external}
\newtheorem{theorem}{Theorem}
\newtheorem{corollary}[theorem]{Corollary}
\newtheorem{lemma}[theorem]{Lemma}
\newtheorem{remark}{Remark}
\newtheorem{proposition}{Proposition}

\newtheorem{assumption}{Assumption}

\DeclareMathOperator*{\diag}{diag}

\DeclareMathOperator*{\sgn}{sgn}

\newcommand{\ud}{\,\mathrm{d}}

\newcommand{\beq}{\begin{equation}}
\newcommand{\eeq}{\end{equation}}
\newcommand{\bq}{\begin{eqnarray}}
\newcommand{\eq}{\end{eqnarray}}
\newcommand{\bqn}{\begin{eqnarray*}}
\newcommand{\eqn}{\end{eqnarray*}}
\newcommand{\bee}{\begin{enumerate}}
\newcommand{\eee}{\end{enumerate}}

\newlength\fheight 
\newlength\fwidth 

\begin{document}
\maketitle

\begin{abstract}
This paper analyzes distributed control protocols for first- and second-order networked dynamical systems. We propose a class of nonlinear consensus controllers where the input of each agent can be written as a product of a nonlinear gain, and a sum of nonlinear interaction functions. By using integral Lyapunov functions, we prove the stability of the proposed control protocols, and explicitly characterize the equilibrium set. 
We also propose a distributed proportional-integral (PI) controller for networked dynamical systems. The PI controllers successfully attenuate constant disturbances in the network. We prove that agents with single-integrator dynamics are stable for any integral gain, and give an explicit tight upper bound on the integral gain for when the system is stable for agents with double-integrator dynamics. 
Throughout the paper we highlight some possible applications of the proposed controllers by realistic simulations of autonomous satellites, power systems and building temperature control.
%
%
\end{abstract}
\noindent
\section{Introduction} 
 
\subsection{General motivation} 
Distributed or decentralized control is in many large-scale systems the only feasible control strategy, when sensing and actuation communication is limited. In this paper we distinguish between decentralized control strategies where agents only have access to local measurements, from distributed control strategies where agents have access to local measurements and the measurements from neighboring agents. 
To attenuate unknown disturbances as well as to stabilize systems in the presence of model imperfections, proportional-integral (PI) controllers are commonly employed in such plants. However, it is still an open problem under which conditions distributed PI-controllers can stabilize a plant in general \cite{MorZaf1989}. 
One important class of systems that require integral action to attenuate disturbances are electric power systems \cite{machowski2008power}. Due to sudden load and generation changes as well as model imperfections, a proportional frequency controller cannot reach the desired reference frequency in general. To attenuate such static errors, integrators are commonly employed.
Due to the inherent difficulties with distributed PI control, automatic frequency control of power systems is typically carried out at two levels, an inner and an outer level. In the inner control loop, the frequency is controlled with a proportional controller against a dynamic reference frequency. In the outer loop, the reference frequency is controlled with a centralized PI controller to eliminate static errors. While this control architecture works satisfactorily in most of today's power systems, future power system developments might render it unsuitable. For instance, large-scale penetration of renewable power generation increases generation fluctuations, creating a need for fast as well as local disturbance attenuation. Decentralized control of power systems might also provide efficient control under islanding and self-healing features, even when communication between subsystems is limited or even unavailable \cite{Senroy2006}. Distributed integral action for power system frequency control, as well as numerous other physical systems such as mobile robotic systems, can often be well modeled by consensus-like protocols for double-integrator dynamics. Motivated by this fact, we develop in this paper a more general framework for consensus protocols with integral action for agents with double-integrator dynamics and velocity damping.  

Static feedback controllers, in contrast to integrators, have been used in a variety of distributed control problems. For example they are commonly used in the consensus problem. The consensus  problem
has applications in formation control \cite{Fax2004}, flocking \cite{Liu03stabilityanalysis} and rendezvous \cite{DistCtrlRobotNetw}, amongst others. While static linear feedback controllers can solve the basic problems, nonlinear feedback controllers are of great importance when, e.g., dealing with connectivity constraints and collision avoidance \cite{Dimarogonas2008}. Furthermore, some distributed control problems have inherent nonlinear dynamics, e.g.,  heat diffusion where heat capacities are nonlinear. Much attention has been devoted to nonlinear consensus protocols in recent years. However, to the best of our knowledge, none of the previous work characterizes the convergence point of nonlinear consensus protocols explicitly. In this paper we introduce a nonlinear consensus protocol for single and double-integrator dynamics. We prove the asymptotic stability of the dynamics induced by the protocol, and explicitly characterize its equilibrium set.

\subsection{Related work}
Nonlinear interaction functions for consensus problems are  well-studied \cite{Saber_consensusprotocols, Moreau2005}, in applications to, e.g., consensus with preserving connectedness \cite{Dimarogonas2008} and collision avoidance \cite{Dimarogonas2008}. 
Sufficient conditions for the convergence of nonlinear protocols for first-order integrator dynamics are given in \cite{Ajorlou2011625}.
Consensus on a general nonlinear function value, referred to as $\chi$-consensus, was studied in \cite{cortes2006}, by the use of nonlinear gain functions. The literature on $\chi$-consensus has been focused on agents with single-integrator dynamics. However, as we show later, our results can be generalized to hold also for double-integrator dynamics. 
Consensus protocols where the input of an agent can be separated into a product of a positive function of the agents own state were studied in \cite{Bauso2006918} for single integrator dynamics. 
In \cite{Munz2011}, position consensus for agents with double-integrator dynamics under a class of nonlinear interaction functions and nonlinear velocity-damping is studied. In contrast to \cite{Munz2011}, we study undamped consensus protocols for single- and double-integrator dynamics using integral Lyapunov functions. 
In \cite{Xie2007}  double-integrator consensus problems with linear non-homogeneous damping coefficients are considered. 
We generalize these results to also  hold for a class of nonlinear damping coefficients.

Multi-agent systems, as all control processes, are in general sensitive to disturbances. When only relative measurements are available, disturbances may spread through the network. It has for example been shown by \cite{bamieh2012coherence} that vehicular string formations with only relative measurements cannot maintain coherency under disturbances, as the size of the formation increases. In \cite{young2010robustness} the robustness of consensus-protocols under disturbances is studied, but limited to the relative states of the agents. However, none of the aforementioned references consider disturbance rejection. In \cite{xiao2007distributed} however, the steady-state error for first-order consensus dynamics is minimized by convex optimization over the edge-weights.
A similar approach is taken in \cite{lin2012optimal}, where the application is vehicle-platooning. 
 In \cite{zelazo2011edge}, an optimal sensor placement problem for consensus problems is formulated, minimizing the $\mathcal{H}_2$ gain of the system. However, these approaches eliminate output errors if the disturbances are constant only in special cases, as no no integral control is employed. 

Consensus with integral action is studied in \cite{Freeman2006} for agents with single-integrator dynamics. It was shown that the proposed controller attenuates constant disturbances. In \cite{Yucelen2012}, the authors take a similar approach to attenuate unknown disturbances. In both references the analysis is limited to agents with single-integrator dynamics.
Our proposed PI controller is related to the consensus controllers in \cite{Cheng2008,Hong2007}. However, these references do not consider the influence of disturbances.

\subsection{Main contributions}
The main contributions of this paper are threefold. Firstly,
we propose  a class of nonlinear consensus protocols where the input of an agent can be decoupled into a product of a positive gain function of the agents own state, and a sum of interaction functions of its neighbors' relative states. 
Nonlinear consensus protocols with double-integrator dynamics and state-dependent damping are also considered.  
 The equilibria are characterized for both the first- and second-order nonlinear consensus protocol.
Secondly, we propose a distributed PI controller for multi-agent systems with single integrator dynamics and a corresponding controller for agents with damped double-integrator dynamics. Necessary and sufficient conditions for stability are derived, given that the controller gains are uniform. 
Lastly, we demonstrate some applications of our theoretical results to satellite control, mobile robots, green buildings and frequency control of power systems. 

\subsection{Outline}
The remainder of this paper is organized as follows. In Section \ref{sec:prob_for} we introduce four motivating examples, which illustrate the theoretical results developed later in the paper. In Section \ref{sec:model} we define a mathematical model for a class of distributed control systems, and the considered problem is stated. In Section \ref{sec:nonlinear} we analyze the problem of distributed control for a class of nonlinear control problems. We explicitly characterize the equilibria by integral equations, and provide sufficient stability criteria by integral Lyapunov functions. In section  \ref{sec:integral} we analyze the problem of distributed control by integral action, where we give necessary and sufficient stability criteria. In Section \ref{sec:applications} we show the feasibility of the results by applying them to the examples in Section \ref{sec:prob_for}. The paper ends by concluding remarks in Section \ref{sec:conclusions}.

\section{Motivating applications}
\label{sec:prob_for}

\subsection{Thermal energy storage in smart buildings}
\label{subsec:green_intro}
Thermal energy storage has emerged as a possible method for energy-efficient temperature regulation in buildings, as discussed in \cite{Zalba2003251}. By using substances which undergoes a phase transition near the desired maximum temperature in the building, the temperature may be kept below the maximal desired temperature. While the heat capacity of the air in a building is approximately constant, the total heat capacity of the room is highly nonlinear due to the thermal energy storage. The endothermic  and exothermic processes of the phase transitions may be modeled by nonlinear heat capacities, which take the form of a continuous approximation of a Dirac delta function at the temperature of the phase transition. The model fits well with a consensus protocol for agents with single-integrator dynamics with nonlinear gain and interaction functions. By a nonlinear extension of Fourier's law  \cite{incropera2011fundamentals}, the room temperatures are thus well-described by the following nonlinear differential equation
\begin{equation}
\dot{T}_i = - \gamma_i(T_i) \sum_{j\in \mathcal{N}_i} a_{ij}(T_i-T_j),
 \label{eq:ex_building_temperature}
\end{equation}
where $T_i$ is the temperature of room $i$, $a_{ij}(T_i-T_j)$ is the heat conductivity between room $i$ and $j$. Here $a_{ij}(\cdot)$ is a nonlinear function $\forall (i,j)\in\mathcal{E}$, and $1/\gamma_i(T_i)$ is the temperature-dependent heat capacity of room $i$, capturing the dynamics of the energy storage. It is of interest to determine the asymptotic temperature in the rooms given their initial temperatures.

\subsection{Autonomous space satellites}
\label{ex:autonomous_space}
Groups of autonomous space satellites may solve tasks in space that require coordination. E.g., for a solar power plant in space, this could involve formation control of mirrors, reflecting the sunlight to a solar panel. If the agents are far away from any fixed reference frame, it may be assumed that the satellites only have access to their distance and velocity relative to their neighboring satellites. It is however often important to analyze the dynamical behavior of the satellites from a common reference frame, e.g., from the earth. Even if the control laws are linear in the relative velocities in the satellites reference frame, they are in general nonlinear in another reference frame. More specifically, the dynamics of a group of $N$ satellites are assumed to be governed by Newton's second law of motion, resulting in a set of second-order dynamical systems. 
 The control signal is the power applied by each agent's engine, $P_i$. However, the acceleration in an inertial reference frame is $a_i={P_i}/{|v_i|}$ due to $P_i= \langle F_i,v_i \rangle $, and the force $F_i$ being parallel to $v_i$, where $v_i$ is agent $i$'s velocity. We assume that the agents only have access to relative measurements. 
The following dynamics describe the interactions of the satellites
\begin{align}
\begin{aligned}
\dot{x}_i&=v_i \\
\dot{v}_i 		&=  - \dfrac{1}{|v_i|+c} \sum_{j \in \mathcal{N}_i} \left[ a_{ij} \left( x_i-x_j \right) + b_{ij} \left( v_i-v_j \right) \right],
\end{aligned}
\label{nnlsecorder3}
\end{align}
where $a_{ij}(\cdot)$ and $b_{ij}(\cdot)$ are nonlinear interaction functions, $i=1, \dots n$ and $\mathcal{N}_i$ denotes the neighbor set of satellite $i$, and where $c$ is an arbitrarily small constant which ensures the boundedness of $\dot{v}_i$ when $v_i=0$. Here $x_i$ and $v_i$ denote the position and velocity of satellite $i$.
 This motivates the analysis of consensus protocols for agents with double-integrator dynamics and nonlinear gain and interaction functions. 

\subsection{Mobile robot coordination under disturbances}
\label{subsec:ex_mobile_robots}
As all control systems, mobile robot systems are susceptible to disturbances. In general, even constant disturbances cause a robot formation with only relative position and velocity measurements to drift. We will consider the particular control objective of reaching position-consensus.
To address the issues caused by disturbances to the robots, a distributed PI controller is proposed.  We consider robots with the second-order dynamics and the following controller:
\begin{align}
\begin{aligned}
\dot{x}_i &= v_i \\
\dot{v}_i &= u_i  -\gamma v_i + d_i \\
u_i &= -\sum_{j\in \mathcal{N}_i} \left(b(x_i-x_j) +  a \int_0^t (x_i(\tau)-x_j(\tau)) \ud \tau  \right), 
\end{aligned}
\label{eq:ex_mobile_robots}
\end{align}
where $x_i$ is the position, $v_i$ is the velocity. The constants $a, b>0$ are controller parameters, while $\gamma > 0$ is a damping coefficient, and $d_i$ is a constant disturbance. We will investigate when distributed PI controllers can attenuate static disturbances in mobile-robot networks. 
Furthermore, given the system-specific damping coefficient $\gamma$, we will characterize under which conditions on $a$ and $b$ asymptotic consensus is reached.

\subsection{Frequency control of power systems}
\label{subsec:intro_power}
Power systems are among the largest and most complex dynamical systems  created by mankind \cite{machowski2008power}. The interconnectivity of power systems poses many challenges when designing controllers. 
We model the power system as interconnected second-order systems, known as the swing equation. The swing equation has been used, e.g., in studying transient angle stability of power systems \cite{DoerflerBullo2011} and fault detection in power systems \cite{Shames20112757}, and is given by
\begin{eqnarray}
 \label{eq:swing_intro}
m_i\ddot{\delta}_i + d_i \dot{\delta}_i = -\sum_{j\in \mathcal{N}_i} k_{ij}(\delta_i - \delta_j) + p^m_i + u_i,
\end{eqnarray} 
where $\delta_i$ is the phase angle of bus $i$, $m_i$ and $d_i$ are the inertia and damping coefficient respectively, $p_i^m$ is the electrical power load at bus $i$ and $u_i$ is the mechanical input power. The coefficient $k_{ij} = |V_i||V_j|b_{ij}$, where $|V_i|$ is the absolute value of the voltage of bus $i$, and $b_{ij}$ is the susceptance of the line $(i,j)$. The frequency of the power system is denoted $\omega_i=\dot{\delta}_i$. 
Maintaining a steady frequency is one of the major control problems in power systems. If the frequency is not kept close to the nominal operational frequency, generation and utilization equipment may cease to function properly. The frequency is maintained primarily by automatic generation control (AGC), which is carried out at different levels. In the first level, which is carried out locally at each bus, the power generation  is controlled by the deviation from a dynamic reference frequency.  At the second level, which is carried out by a centralized controller, the reference frequency is controlled based on the deviation of the frequency from a reference frequency at a specific location. When communication is unavailable or limited, a decentralized controller is a possible alternative. 
 A simple decentralized frequency control with integral action would take the form:
\begin{align}
u_i &= a ( \omega^{\text{ref}} - \omega_i(t)) + b \int_0^t ( \omega^{\text{ref}} - \omega_i(t')) \ud t',
 \label{eq:decentralized_control_22}
\end{align}
where $ \omega^{\text{ref}} $ is the reference frequency.
The control objective is to ensure that the system frequency reaches the nominal operational frequency, i.e.,
$
\lim_{t\rightarrow \infty}\omega_i = \omega^{\text{ref}} \; \forall i \in \mathcal{V}.
$

\section{Problem formulation}
\label{sec:model}
In this section we  formalize the previously mentioned problems. We introduce a unified mathematical notation which includes all problems previously described.

\subsection{Notation}
Let $\mathcal{G}$ be a graph. Denote by $\mathcal{V}=\{ 1,\hdots, n \}$ the vertex set, and by $\mathcal{E}=\{ 1,\hdots, m \}$ the edge set of $\mathcal{G}$. Let $\mathcal{N}_i$ be the set of neighboring nodes of node $i$. 
Denote by $\mathcal{B}$ the vertex-edge adjacency matrix of $\mathcal{G}$, and let $\mathcal{L}$ be its Laplacian matrix. For undirected graphs it holds that $\mathcal{L}=\mathcal{B}\mathcal{B}^T$. Throughout this paper we will assume that $\mathcal{G}$ is connected, as motivated by the applications considered. 
We denote by $\mathbb{R}^-$/$\mathbb{R}^+$ the open negative/positive real axis, and by $\bar{\mathbb{R}}^-$/$\bar{\mathbb{R}}^+$ its closure. 
Let $\mathbb{C}^-$/$\mathbb{C}^+$ denote the open left/right half complex plane, and $\bar{\mathbb{C}}^-$/$\bar{\mathbb{C}}^+$ its closure. 
We will denote the scalar position of agent $i$ as $x_i$, and its velocity as $v_i$, and collect them into column vectors $x=(x_1, \hdots, x_n)^T$, $v=(v_1, \hdots, v_n)^T$. We denote by $c_{n\times m}$ a vector or matrix of dimension $n\times m$ whose elements are all equal to $c$. $I_{n}$ denotes the identity matrix of dimension $n$. 
A function $ f(\cdot)$ with domain $\mathcal{X}$ is said to be globally Lipschitz (continuous) if there exists $K\in \mathbb{R}^+:\; \forall x,y \in \mathcal{X}: \; \norm{ f(x)- f(y)} \le  K \norm{x-y}$.

\subsection{System model}
We consider agents endowed with either single-integrator dynamics
\begin{align}
\label{eq:dynamics_single}
\dot{x}_i&= d_i + u_i,
\end{align}
or double-integrator dynamics
\begin{align}
 \begin{aligned}
\dot{x}_i&=v_i  \\
\dot{v}_i&= d_i + u_i,
\end{aligned}
\label{eq:dynamics_double_2} 
\end{align}
where $d_i$ is a constant disturbance. 
\subsection{Objective}
The main objectives of this paper are twofold. Our first objective is to characterize the stability of nonlinear feedback protocols, and to determine under which nonlinear feedback protocols the consensus point of the agents may be determined a priori, in the absence of disturbances.  The second objective is the design of consensus protocols robust to disturbances. We will focus on constant but unknown disturbances.
In both cases, the overall objective will be for all agents to converge to a common state, i.e. $\lim_{t\rightarrow \infty} x_i(t)=x^* \; \forall i \in \mathcal{V}$ for single-integrator dynamics, and $\lim_{t\rightarrow \infty} v_i(t)=v^* \; \forall i \in \mathcal{V}$ for double-integrator dynamics. 


\section{Distributed control with static nonlinear feedback}
\label{sec:nonlinear}
Although the consensus problem has received tremendous amounts of attention in the last decade, the vast majority of the studied consensus protocols are linear. In this section we define a class of nonlinear consensus algorithms where the input of each agent can be decoupled into a product of a nonlinear gain function and a nonlinear interaction function. We show that several interesting properties of the system arise due to this coupling. We first study consensus for single-integrator dynamics by nonlinear protocols in Section~\ref{subsec:nonlinear_single_int}. In Section~\ref{sec:nnlsecorder} we consider a nonlinear consensus protocol for agents with double-integrator dynamics. Finally, in Section~\ref{sec:consensus_damping} a nonlinear consensus protocol with nonlinear, state-dependent damping for agents with double-integrator dynamics is considered. 

\subsection{Consensus for single-integrator dynamics}
\label{subsec:nonlinear_single_int}
Consider agents with dynamics \eqref{eq:dynamics_single}, where $d_i=0 \; \forall i \in \mathcal{V}$, and where $u_i$ is given by:
\begin{equation}
u_i = - \gamma_i(x_i) \sum_{j\in \mathcal{N}_i} a_{ij}(x_i-x_j).
\label{eq:ag}
\end{equation}
\begin{assumption}
\label{ass:gamma}
The gain $\gamma_i$ is continuous and 
$ 0< \underline{\gamma} \le  \gamma_i(x)\le \bar{\gamma} \quad \forall i \in\mathcal{V},\; \forall x \in \mathbb{R}.
\label{eq:ass5}$ \end{assumption}
\begin{assumption}
\label{ass:alpha}
The interaction function $a_{ij}(\cdot)$ is Lipschitz continuous $ \forall (i,j) \in \mathcal{E}$, and:
\begin{enumerate}
\item
$a_{ij}(-y) = -a_{ji}(y) \quad \forall  (i,j)\in \mathcal{E},\; \forall y \in \mathbb{R}$
\item
$a_{ij}(-y) = - a_{ij}(y)  \quad \forall  (i,j)\in \mathcal{E},\; \forall y \in \mathbb{R}$
\item
 $ a_{ij}(y)   > 0 \quad \forall  (i,j)\in \mathcal{E} ,\; \forall y>0,$
\end{enumerate}
 \end{assumption}
\begin{remark}
Assumption \ref{ass:alpha} guarantees that the agents move in the direction of their neighbors, as well as symmetry in the flow. A consequence of Assumption \ref{ass:alpha} is that $\alpha_{ij}(0) =0 \;  \forall  (i,j)\in \mathcal{E}$, ensuring that any consensus point where $x_i=x_j \; \forall i,j\in \mathcal{V}$ is an equilibrium. 
\end{remark}
We are now ready to state the main result of this section.
\begin{theorem} 
\label{prop:1}
Given $n$ agents with dynamics \eqref{eq:dynamics_single} with $d_i=0 \; \forall i \in \mathcal{V}$, and $u_i$ given by \eqref{eq:ag}, where  $\gamma_i(\cdot)$ and $a_{ij}(\cdot)$ satisfy Assumptions \ref{ass:gamma} and \ref{ass:alpha}, respectively, for all $i\in \mathcal{V}$ and for all $(i,j)\in \mathcal{E}$. 
 Then the agents converge asymptotically to an agreement point $\lim_{t\rightarrow \infty} x_i(t)=x^*  \;\forall\; i \in \mathcal{V}$, where $x^*$ is uniquely determined by the integral equation for any  $x_i(0)=x_i^0, i=1, \dots, n$,
\begin{equation}
\sum_{i\in \mathcal{V}}  \int_0^{x^0_i} \frac{1}{\gamma_i(y)} \ud y = \int_0^{x^*} \sum_{i\in \mathcal{V}}  \frac{1}{\gamma_i(y)} \ud y, 
 \label{eq:sol}
\end{equation}
\end{theorem}

\begin{proof}
The proof is Lyapunov based, and requires finding an appropriate Lyapunov function, which guarantees convergence to an equilibrium set. This equilibrium set is then characterized by a time-invariant function. Such a function is given by
$ E(x) = \sum_{i\in \mathcal{V}}  \int_0^{x_i} {1}/{\gamma_i(y)} \ud y.
$ Differentiating $E(x)$ with respect to time yields
\begin{align*}
{ \dot{E}(x(t))} &= - \left[ \frac{1}{\gamma_1(x_1)}, \hdots,  \frac{1}{\gamma_n(x_n)}\right] \Gamma(x) \mathcal{B} a({\mathcal{B}^T	x}) \\ 
&= -{1}_{1\times n}\mathcal{B}a(\mathcal{B}^T{x}) = 0,
\end{align*}

where  $\Gamma(x)=\diag([\gamma_1(x_1), \dots, \gamma_n(x_n)])$, and $a(\cdot)$ is taken component-wise.
Hence $E(x)$ is time-invariant and the agreement point $x^*$, if existing, is given by \eqref{eq:sol}. By Assumption \ref{ass:gamma}, $E(x^* 1_{n\times 1})$ is strictly increasing in $x^*$, and hence \eqref{eq:sol} admits a unique solution.
Now consider  the following candidate Lyapunov  function:
\begin{eqnarray}
V(x) = \sum_{i\in \mathcal{V}} \int_{x^*}^{x_i} \frac{y-x^*}{\gamma_i(y)} \ud y 
\label{eq:lyap_switch},
\end{eqnarray}
where $x^*$ is the agreement point given by \eqref{eq:sol}. It can easily be verified that $V(x^*1_{n\times 1}) = 0$. To show that $V(x)>0$ for $x\ne x^*1_{n\times 1}$, it suffices to show that $ \int_{x^*}^{x_i} {(y-x^*)}/{\gamma_i(y)} \ud y > 0 \; \forall i \in \mathcal{V}, \; \forall x\ne x^*1_{n\times 1}$. Consider first the case when $x_i>x^*$, where 
$
\int_{x^*}^{x_i} ({y-x^*})/{\gamma_i(y)} \ud y = \int_{0}^{x_i-x^*} {z}/{\gamma_i(z+x^*)} \ud z > 0,
$
by the change of variable $z=y-x^*$. 
The case when $x_i<x^*$ is treated analogously, where 
$
\int_{x^*}^{x_i} {(y-x^*)}/{\gamma_i(y)} \ud y = \int_{0}^{x^*-x_i} {z}/{\gamma_i(x^*-z)} \ud z > 0,
$
by the change of variable $z=x^*-y$. This also implies that $V(x)=0 \Rightarrow x=x^* 1_{n\times 1}$.
Now consider $\dot{V}(x)$: 
\begin{align}
&\dot{V}(x) = - \sum_{i \in \mathcal{V}}  \frac{x_i-x^*}{\gamma_i(x_i)}  \gamma_i(x_i) \sum_{j\in \mathcal{N}_i} a_{ij}(x_i-x_j) \nonumber  \\
&=  - \sum_{i \in \mathcal{V}} x_i  \sum_{j\in \mathcal{N}_i} a_{ij}(x_i-x_j) 
+ \sum_{i \in \mathcal{V}} x^*  \sum_{j\in \mathcal{N}_i} a_{ij}(x_i-x_j).
 \label{eq:proof_th_nl_lyap}
\end{align}
Due to the symmetry property in Assumption \ref{ass:alpha}, the first term of \eqref{eq:proof_th_nl_lyap} may be rewritten as 
$
 \sum_{i \in \mathcal{V}} x_i  \sum_{j\in \mathcal{N}_i} a_{ij}(x_i-x_j) = \sum_{i \in \mathcal{V}} \sum_{j\in \mathcal{N}_i} x_i  a_{ij}(x_i-x_j) =  \frac 12 \sum_{i \in \mathcal{V}}  \sum_{j\in \mathcal{N}_i} (x_i-x_j)  a_{ij}(x_i-x_j)  
$. 
Clearly the second term of \eqref{eq:proof_th_nl_lyap} satisfies $ \sum_{i \in \mathcal{V}} x^*  \sum_{j\in \mathcal{N}_i} a_{ij}(x_i-x_j)=0$ due to Assumption~\ref{ass:alpha}. Hence, $\dot{V}(x)$ may be rewritten as
$
\dot{V}(x) =  - \frac 12 \sum_{i \in \mathcal{V}}  \sum_{j\in \mathcal{N}_i} (x_i-x_j)  a_{ij}(x_i-x_j)   < 0,
$ 
unless $x_i=x_j \; \forall i,j \in \mathcal{V}$. 
Therefore the agents converge to  $x_i=x^* \; \forall i \in \mathcal{V}$. 
\end{proof}

\begin{remark}
The agreement protocol \eqref{eq:ag} has a physical interpretation. If we consider the smart building problem in Section \ref{subsec:green_intro}, and let $x_i$ be the temperature of the rooms, ${1}/{\gamma_i(x_i)}$ is the temperature-dependent heat capacity of the nodes. Analogously, $a_{ij}(\cdot) $ is the thermal conductivity of the walls between rooms $i$ and $j$, being dependent on the temperature gradient between the rooms. The invariant quantity $E(x)=\sum_{i\in \mathcal{V}}  \int_0^{x_i} {1}/{\gamma_i(y)} \ud y$ is the total stored thermal energy.
\end{remark}
\begin{remark}
The convergence of the dynamics \eqref{eq:dynamics_single}, and the control \eqref{eq:ag} were proven in \cite{Shi20091165}. However, as opposed to that reference, we characterize here explicitly the equilibrium set. This gives a priori knowledge about the point of convergence. 
 Furthermore, in \cite{Shi20091165}, the convergence was studied with a Lyapunov function consisting of the difference of the maximal and minimal state, as opposed to the Lyapunov function consisting of integral equations used in our proof. This new Lyapunov function facilitates the proof by not relying on non-smooth analysis.
\end{remark}

\subsection{Consensus for double-integrator dynamics}
\label{sec:nnlsecorder}
Consider agents with double-integrator dynamics \eqref{eq:dynamics_double_2}, where $d_i=0 \; \forall i \in \mathcal{V}$, with $u_i$ given by
\begin{align}
u_i 		&=  - \gamma_i(v_i) \sum_{j \in \mathcal{N}_i} \left[ a_{ij} \left( x_i-x_j \right) + b_{ij} \left( v_i-v_j \right) \right]. \label{nnlsecorder3}
\end{align}
We show that under mild conditions, the controller \eqref{nnlsecorder3} achieves asymptotic consensus. The following Theorem generalizes the literature on linear second-order consensus as in \cite{Ren2008Book}, by extending the analysis of linear gains and interaction functions to nonlinear ones. The nonlinear analysis covers a much richer class of control inputs than the corresponding linear analysis, and we are again able to explicitly characterize the equilibrium set. 
By using an integral Lyapunov function, we are thus able to prove that the agents reach consensus for the nonlinear consensus protocol also in the case of double-integrator dynamics.
\begin{theorem}
\label{prop:nlsecorder}
Consider agents with dynamics \eqref{eq:dynamics_double_2}, with $d_i=0 \; \forall i \in \mathcal{V}$, and $u_i$ given by \eqref{nnlsecorder3}, where $a_{ij}(\cdot)$ and $\gamma_i(\cdot)$ satisfy Assumptions \ref{ass:gamma} and \ref{ass:alpha}, and $ b_{ij}(\cdot)$ satisfies Assumption \ref{ass:alpha}, replacing $a_{ij}$ with $b_{ij}$, for all $i\in \mathcal{V}$ and for all $(i,j)\in \mathcal{E}$, respectively.
The agents achieve consensus with respect to $x$ and $v$, i.e., $|x_i-x_j|\rightarrow 0,|v_i-v_j|\rightarrow 0\;\forall i,j\in\mathcal{G} \text{ as } t\rightarrow \infty$  for any initial condition $(x(0), v(0))$. Furthermore, the velocities converge to $\lim_{t\rightarrow \infty}  v_i(t) = v^* \;\forall\; i \in \mathcal{V}$, uniquely determined by 
\begin{align}
\label{eq:consensus_velocity_double_integrator}
 \sum_{i\in \mathcal{V}} \int_0^{v^0_i} \frac{1}{\gamma_i(y)} \ud y = \int_0^{v^*} \sum_{i\in \mathcal{V}}  \frac{1}{\gamma_i(y)} \ud y.
\end{align} \end{theorem}

\begin{proof}
The convergence analysis relies on Lyapunov techniques. For characterizing the equilibrium set, a time-invariant function is introduced.
We write \eqref{nnlsecorder3} in vector form as
\begin{align*}
\dot{ x} &= v \\ 
\dot{ v} &= - \Gamma(v) \left[ \mathcal{B}a(\bar{x}) + \mathcal{B}b({\mathcal{B}^T v}) \right],
\end{align*} 
where 
$\bar{x} = \mathcal{B}^T x$, and  $a(\cdot)$ and $b(\cdot)$ are taken component-wise, and \\ 
 $\Gamma(x)=\diag([\gamma_1(x_1), \dots, \gamma_n(x_n)])$. 
Consider now the following candidate Lyapunov function, also used in \cite{Munz2011}, however in another context,
$$
 V(\bar{x},v) \allowbreak= \allowbreak \sum_{i\in \mathcal{V}} \left( \int_{v^*}^{v_i} \frac{(y{-}v^*)}{\gamma_i(y)} \ud y\right)\allowbreak + \allowbreak\sum_{(i,j)\in \mathcal{E}}   \int_0^{\bar{x}_{ij}} a_{ij}(y)\ud y.
$$ 
Here $v^*$ is the common velocity of the agents at the equilibrium, which we will show is given by \eqref{eq:consensus_velocity_double_integrator}. 
It is straightforward to verify that $V([0_{1\times m}, v^*1_{1\times n}]^T)=0$. By following the proof of the positive semi-definiteness of $V(x)$ in the proof of Theorem \ref{prop:1}, replacing $x_i$ and $x^*$ with $v_i$ and $v^*$, respectively, the positive semi-definiteness of $ \sum_{i\in \mathcal{V}} ( \int_{v^*}^{v_i} {(y{-}v^*)}/{\gamma_i(y)} \ud y) $ follows. For showing the positive semi-definiteness of the second term of $V(\bar{x},v)$, it suffices to show that $\int_0^{\bar{x}_{ij}} a_{ij}(y)\ud y > 0 \; \forall (i,j) \in \mathcal{E}$. For $\bar{x}_{ij}>0$, this inequality clearly holds. When $\bar{x}_{ij}<0$ we have
$
\int_0^{\bar{x}_{ij}} a_{ij}(y)\ud y = - \int^0_{\bar{x}_{ij}} a_{ij}(y)\ud y = \int^0_{\bar{x}_{ij}} a_{ji}(-y)\ud y > 0. 
$
We may write $V(\bar{x},v)$ in vector form as
$
 V(\bar{x},v) = \int_0^{\bar{x}}  {1}_{1 \times n} \mathcal{B}^Ta(y) \ud y + \int_{v^*{1}_{n\times 1}}^v \tilde{y}^T \Gamma^{-1}(y) {1}_{n \times 1} \ud y,
$
where
$ \tilde{y}= [
y_1{-}v^*, 
{\hdots} 
,y_n{-}v^*
]^T
$. Differentiating $V(\bar{x},v)$ with respect to time yields:
\begin{align*}
\dot{V} &=  a(\bar{x})^T \mathcal{B}^T v{-}(v- v^* {1}_{n\times 1})^T \mathcal{B}\left[ a(\bar{x}) + b(\mathcal{B}^T{v}) \right]  \\
&= - v^T  \mathcal{B}b(\mathcal{B}^T{v}) +  v^*1_{1\times n}  \mathcal{B}b(\mathcal{B}^T{v})  = - {v}^T\mathcal{B}b(\mathcal{B}^T{v}) \le 0
\end{align*}
due to Assumption \ref{ass:alpha}, with equality if and only if $\mathcal{B}^T{v}=0$.
We now invoke LaSalle's invariance principle to show that the agreement point satisfies $\dot{v}=0$.
The subspace where $\dot{ V}(\bar{x},v)=0$ is given by $S_1=\left\{ (\bar{x},v)| v=c{1}_{n\times 1} \right\}$. We note that on $S_1$, 
$\dot{v} =  - \Gamma(v) \left[ \mathcal{B}a(\bar{x}) + \mathcal{B}b(\mathcal{B}^Tv) \right] = - \Gamma(v)  \mathcal{B}a(\bar{x})  \ne k(t)  {1}_{n\times 1}.
$
To see this, suppose that 
$\dot{v}(t) = - \Gamma(v)  \mathcal{B}a(\bar{x}) = k(t) {1}_{n\times 1} \Leftrightarrow  \mathcal{B}a(\bar{x}) = \Gamma^{-1}(v)k(t) {1}_{n\times 1},
$ where $k(t) \ne 0 \forall t$. Premultiplying the above equation with ${1}_{1\times n}$ yields
$ 0 = {1}_{1\times n}    \mathcal{B}a(\bar{x}) =k(t) 1_{1\times n} \Gamma^{-1}(v) 1_{n\times 1} \ne 0,
$ which is a contradiction since $k(t)\ne 0 \forall t$ by assumption. Hence the only trajectories contained in $S_1$ are those where $v=v^*{1}_{n\times 1}, \dot{v}=0$. 
It can also be shown that no trajectories where $\bar{x}\ne 0$ are contained in $S_1$. Assume for the sake of contradiction that $\bar{x}\ne 0$ in $S_1$. Let $i^- = \min_{ j\in \mathcal{V}}  x_j \text{ s.t. } \exists  k\in \mathcal{N}_{i^-}:  x_k> x_{i^-}$. It is clear that such an $i^-$ exists, since otherwise $\bar{x}=0$. Consequently 
\begin{align*}
\dot{v}_{i^-} &= -\gamma_{i^-}(v_{i^-}) \sum_{j \in \mathcal{N}_{i^-}} \left[ a_{{i^-}j} \left( x_{i^-}-x_j \right) + b_{ij} \left( v_{i^-}-v_j \right) \right] \\
 &=  -\gamma_{i^-}(v_{i^-}) \sum_{j \in \mathcal{N}_{i^-}} \left[ a_{{i^-}j} \left( x_{i^-}-x_j \right)  \right] \\ 
 &>  -\gamma_{i^-}(v_{i^-}) a_{{i^-}k} \left( x_{i^-}-x_k \right) > 0,
\end{align*} 
by the assumption that $ x_k> x_{i^-}$. Thus, any trajectory in $S_1$ where  $\bar{x}\ne 0$ cannot remain in $S_1$, implying 
that $|x_i{-}x_j|\rightarrow 0, |v_i{-}v_j|\rightarrow 0  \; \forall  i,j \in \mathcal{G} \text{ as } t\rightarrow \infty$ and furthermore $\dot{v}(t)=0_{n\times 1}$. 
Next we show that
$ p(v)=\sum_{i\in\mathcal{V}} \int_0^{v_i} {1}/{\gamma_i(y)} \ud y = \int_0^{v} 1_{1\times n} \Gamma^{-1}(v)1_{n\times 1} \ud y
$
is invariant under \eqref{nnlsecorder3}. Consider:
\begin{align*}
{ \dot{p}(v(t))} &= - 1_{1\times n} \Gamma^{-1}(v)  \Gamma(v) \left[ \mathcal{B}a(\bar{x}) + \mathcal{B}b(\mathcal{B}^T{v}) \right] \\ 
&= -1_{1\times n} \mathcal{B}a(\bar{x}) -  1_{1\times n} \mathcal{B}b(\mathcal{B}^T{v}) =0 .
\end{align*}
Thus we conclude that $\lim_{t\rightarrow \infty} x(t) = x^*(t)1_{n\times 1}$ and $\lim_{t\rightarrow \infty} v(t) = v^*1_{n\times 1}$ with $v^*$ given by \eqref{eq:consensus_velocity_double_integrator}. 
The existence and uniqueness of the solution to \eqref{eq:consensus_velocity_double_integrator} follows from Assumption \ref{ass:gamma}, and by the proof of Theorem \ref{prop:1}, replacing $x_i$ and $x^*$ with $v_i$ and $v^*$, respectively.
\end{proof}

\begin{remark}
Theorem \ref{prop:nlsecorder} has a physical interpretation. If we regard ${1}/{\gamma_i(v_i)}$ as the velocity-dependent mass of agent $i$, e.g., due to the agents' masses scaling with the Lorentz factor $1/\sqrt{1-v_i^2/v_c^2}$, where $v_c$ is the speed of light, then the invariant quantity 
$p(v) =$ $\sum_{i\in\mathcal{V}} \int_0^{v_i}$ ${1}/{\gamma_i(y)} \ud y$ 
is the total (relativistic) momentum of the mechanical system. 
\end{remark}

\subsection{Consensus for double-integrator dynamics with state-dependent damping}
\label{sec:consensus_damping}
Consider agents with double-integrator dynamics \eqref{eq:dynamics_double_2}, where $d_i=0 \; \forall i \in \mathcal{V}$, with $u_i$ given by:
\begin{align}
u_i &=  - \kappa_i(x_i)v_i  - \sum_{j\in \mathcal{N}_i} a_{ij}(x_i-x_j).  \label{eq:consensus_2_gains}
\end{align}
The following theorem generalizes the results in \cite{Xie2007} to also include nonlinear state-dependent damping, as well as nonlinear interaction functions. With this framework, we are able to generalize the simple average consensus to a much broader class of consensus functions. 

\begin{theorem}
\label{th:1_nl}
Consider agents with dynamics \eqref{eq:dynamics_double_2} and $u_i$ given by \eqref{eq:consensus_2_gains}, where $\kappa_i(\cdot)$ satisfies Assumption \ref{ass:gamma}, replacing $\gamma_i$ with $\kappa_i$ , and $ a_{ij}(\cdot)$ satisfies Assumption \ref{ass:alpha} for all $i\in \mathcal{V}$ and for all $(i,j)\in \mathcal{E}$, respectively. Furthermore, the interaction functions $a_{ij(\cdot)}$ satisfy $\lim_{x\to\infty} \int_0^x a_{ij}(y) \text{d} y = \infty \; \forall (i,j) \in \mathcal{E}$ \footnote{We would like to thank the anonymous reviewer for pointing out the necessity of this condition.}. 
 Then the agents converge to a common point for any initial position $x_i(0)$.
Furthermore, the consensus point is uniquely determined by
\begin{equation}
\label{eq:consensus_point}
\sum_{i\in \mathcal{V}} \Big( \int_0^{x^0_i}{\kappa_i(y)} \ud y + v_i(0) \Big) = \int_0^{x^*} \sum_{i\in \mathcal{V}} {\kappa_i(y)} \ud y.
\end{equation}
\end{theorem}

\begin{proof}
The proof of convergence also relies on Lyapunov techniques, and LaSalle's invariance principle. A time-invariant function is introduced to characterize the equilibrium set. This function is given by
$
E(x,v)= \sum_{i\in \mathcal{V}} \left( \int_0^{x_i} \kappa_i(y) \ud y + {v_i} \right).
$
Differentiating $E(x,v)$ along trajectories of \eqref{eq:consensus_2_gains} yields
$
\dot{E}(x,v) 
= - \sum_{i \in \mathcal{V}} \sum_{j\in \mathcal{N}_i} a_{ij}(x_i-x_j) = 0,
$
by Assumption \ref{ass:alpha}. 
We first note that by Assumptions \ref{ass:gamma} and \ref{ass:alpha}, a unique continuous solution of \eqref{eq:dynamics_double_2} with $u_i$ given by \eqref{eq:consensus_2_gains} exists for all $t \ge 0$.
Consider the candidate Lyapunov function 
$V({x},v) = \sum_{i\in \mathcal{V}} [ {v^2_i}/{2} + \sum_{j\in \mathcal{N}_i} \int_0^{x_i-x_j} a_{ij}(y)\text{d} y  ] \ge 0.$ 
Differentiating $V({x},v)$ along trajectories of \eqref{eq:consensus_2_gains} yields
\begin{align*}
\dot{V}({x},v) &=\sum_{i \in \mathcal{V}}  v_i \big({-}\kappa_i(x_i) v_i{-}\sum_{j\in \mathcal{N}_i} a_{ij}(x_i{-}x_j) \big)\\ 
&+ \sum_{i \in \mathcal{V}} \big( \sum_{j\in \mathcal{N}_i} a_{ij}(x_i{-}x_j) \big) v_i ={-} \sum_{i \in \mathcal{V}} \kappa_i(x_i) v_i^2 \le 0.
\end{align*}

Since $V({x},v)$ is non-increasing under the dynamics  \eqref{eq:consensus_2_gains}, it is clear that for $\Omega=\{[{x}(t), v(t)]: V({x},v)\le V({x}_0,v_0)\}$, it holds that $[\bar{x}(t), v(t)] \in \Omega \; \forall t\ge 0$, where $\bar{x}=\mathcal{B}^Tx$ and ${x}_0={x}(0), v_0=v(0)$. Clearly $\Omega$ is compact. The following lemma is needed to use LaSalle's invariance principle. 
\begin{lemma}
\label{lem:1_nl}
Given the requirements in Theorem \ref{th:1_nl}, $[x(t),v(t)]$ evolve in a compact set, denoted $\Omega'$.
\end{lemma}
\begin{proof}
See Appendix.
\end{proof}
Now let $H=\{(x,y)| \dot{V}({x},v)=0\}=\{(x,y)| v=0\}$. Consider any trajectory of \eqref{eq:consensus_2_gains} with $x\ne x^*(t)1_{n\times 1}$. By \eqref{eq:consensus_2_gains} and the assumption that $\mathcal{G}$ is connected, $\dot{v}_i \ne 0$ for at least one index $i$. Thus the largest invariant manifold of $E$ is $\{(x,v)| x=x^*1_{n\times 1}, v=0 \}$. Since $\Omega'$ by Lemma \ref{lem:1_nl} is compact and positively invariant, then by LaSalle's invariance principle, the agents converge to a common point $x_i=x^* \; \forall i\in \mathcal{G}$, with $v_i=0 \; \forall i\in \mathcal{G}$.

It remains to show that the common point to which the agents converge to is the point given by \eqref{eq:consensus_point}, and that the solution is unique. Indeed, consider again the function
$E(x,v)$. Since $\dot{E}(x,v)=0$, and the agents converge to a point $x^*$ with $v_i = 0 \; \forall i\in \mathcal{V}$. It follows that $x^*$ is given by \eqref{eq:consensus_point}. Since $\kappa_i(y)>0$ by assumption, \eqref{eq:consensus_point} admits a unique solution.
\end{proof}

The following corollary follows directly from Theorem \ref{th:1_nl}.
\begin{corollary}
\label{cor:consensus_point_damping_cor}
Given $n$ agents starting from rest, i.e., $v_i(0)=0 \; \forall i \in \mathcal{V}$, with dynamics \eqref{eq:dynamics_double_2}, where $u_i$ is given by \eqref{eq:consensus_2_gains}, the agents converge to a common point $x^*$ for any initial position $x_i(0)$, which is uniquely determined by
\begin{eqnarray}
\label{eq:consensus_point_damping_cor}
\sum_{i\in \mathcal{V}}  \int_0^{x^0_i}{\kappa_i(y)} \ud y  = \int_0^{x^*} \sum_{i\in \mathcal{V}} {\kappa_i(y)} \ud y.
\end{eqnarray}
\end{corollary}

\begin{remark}
In Theorem \ref{prop:1}, $x^*$ is given by $
\sum_{i\in \mathcal{V}}  \int_0^{x^0_i} {1}/{\gamma_i(y)} \ud y = \int_0^{x^*} \sum_{i\in \mathcal{V}}  {1}/{\gamma_i(y)} \ud y,
$
as opposed to \eqref{eq:consensus_point_damping_cor} in Corollary \ref{cor:consensus_point_damping_cor}. The intuition behind this is that in \eqref{eq:ag}, $\gamma_i(x_i)$ acts as a gain of agent $i$, where an increased $\gamma_i(x_i)$ will increase the speed of agent $i$. In \eqref{eq:consensus_2_gains} however, $\kappa_i(x_i)$ acts as damping of agent $i$, where an increased $\kappa_i(x_i)$ will decrease the speed of agent $i$. 
\end{remark}

\section{Distributed control with integral action}
\label{sec:integral}

{M}{ulti}-agent systems are, like most control processes, sensitive to disturbances. 
Generally, static distributed control protocols, such as e.g., P-controllers, cannot reject even constant disturbances. 
In this section we propose a control protocol for single- and double-integrator dynamics that drives the agents to a common state under static disturbances. By using distributed integral action, we are able to compensate for the disturbances in a distributed setting.  Moreover, with the proposed control algorithm, the agents reach the average of their initial positions for arbitrary initial velocities in the absence of disturbances. We study the properties of the control protocols and derive necessary and sufficient conditions under which the multi-agent system is stable. 

\subsection{Consensus by distributed integral action for single-integrator dynamics with damping}
\label{sec:consensus_integral_single}
Consider agents with single-integrator dynamics \eqref{eq:dynamics_single}, with $u_i$ given by:
\begin{align}
\begin{aligned}
u_i &= -\sum_{j\in \mathcal{N}_i} \left( b(x_i-x_j) +  a  \int_0^t (x_i(\tau){-}x_j(\tau)) \ud \tau  \right) \\ &-\delta (x_i-x_i(0)) 
\end{aligned} 
\label{eq:consensus_single_absolute_2}
\end{align}
where $a \in \mathbb{R}^+$, $b \in \mathbb{R}^+$, $\delta\in  \bar{\mathbb{R}}^+$ are fixed parameters, and $d_i \in \mathbb{R}$ is an unknown disturbance.

\begin{theorem}
\label{th:single_abs}
Under the dynamics  \eqref{eq:consensus_single_absolute_2}, with $u_i$ given by \eqref{eq:consensus_single_absolute_2}, 
the agents converge to a common value $x^*$ for any constant disturbance $d_i$ and any initial condition. If $d_i=0\;\forall i\in\mathcal{V}$, the agents converge to 
$x^*= \frac 1n \sum_{i\in \mathcal{V}} x_i(0) \text{ given } v_i(0) \quad \forall i\in \mathcal{V}$.  
If absolute position measurements are not present, i.e., $\delta=0$, it still holds that $\lim_{t \rightarrow \infty } |x_i(t)-x_j(t)| =0 \; \forall i,j \in \mathcal{V}$ for any $d_i\in \mathbb{R}$, $a,b \in \mathbb{R}^+$, while the absolute states diverge, i.e., $\lim_{t\rightarrow \infty} |x_i(t)| = \infty \; \forall i \in \mathcal{V}$, unless ${1_{1\times n}} d =0$. 
\end{theorem}

\begin{proof}
First consider the case when $\delta=0$ and $d_i=0 \; \forall i\in \mathcal{V}$. 
By introducing the integral states $ z=[z_1, \hdots , z_n]^T$ we may rewrite the dynamics  \eqref{eq:consensus_single_absolute_2} in vector form as
\begin{align}
\label{eq:integral_vector_1}
\left[\begin{smallmatrix}
\dot{z}\\\dot{x}
\end{smallmatrix}\right] = \underbrace{\left[\begin{smallmatrix}
 {0}_{n\times n} & I_{n} \\
-a \mathcal{L} & -b \mathcal{L}
\end{smallmatrix}\right]}_{ \triangleq A} 
\left[\begin{smallmatrix}
{z}\\{x}
\end{smallmatrix}\right],
\end{align}
together with the initial condition $z(0)=0_{n\times 1}$. 
By elementary column operations we note that the characteristic equation of $A$ is given by
$0=\det \left( (a+b s) \mathcal{L} + s^2 I_{n} \right)
$. By comparing the characteristic polynomial with the characteristic equation of $\mathcal{L}$, being
$ 0=\det\left( \mathcal{L} - \kappa I_{n} \right)
$, with solutions $\kappa=\lambda_i \ge 0$, we obtain the equation
$0=s^2 + \lambda_ib s + \lambda_ia$. This equation has a double root $s=0$ if $\lambda_i=0$, and the remaining solutions $s\in \mathbb{C}^-$ if $\lambda_i>0$. Since the above equation has exactly two solutions for every $\lambda_i$, it follows that the algebraic multiplicity of the eigenvalue $0$ must be equal to two. It is well-known that for connected graphs $\mathcal{G}$, $\lambda_1$ is the only zero-eigenvalue of the Laplacian $\mathcal{L}$. By straightforward calculations we obtain that $e_1= [
1_{1\times n} ,\; 0_{1\times n}
]^T$ is an eigenvector and $e_2= [
0_{1\times n} ,\; 1_{1\times n} 
]^T$ is a generalized eigenvector of $A$ corresponding to the eigenvalue $0$. It can also be verified that $v_1= \frac 1n [
{1}_{1\times n}  ,\;  {0}_{1\times n} 
]$ and $ v_2 = \frac 1n [
{0}_{1\times n}  ,\; {1}_{1\times n} 
]$ are a left eigenvector and a generalized left eigenvector of $A$, respectively, corresponding to the eigenvalue 0. It is easily verified that $v_1 e_1 =1, \; v_2 e_2 = 1$ and $v_2 e_1 =0, \; v_1 e_2 = 0$. If we let $P$ be a matrix consisting of the normalized eigenvectors of $A$, we can chose the first columns of $P$ to be $e_1$ and $e_2$, and the first rows of $P^{-1}$ to be $v_1$ and $v_2$, respectively.
Since all remaining eigenvalues of $A$ have strictly negative real part we obtain 
\begin{align*}
\lim_{t\rightarrow \infty} e^{At} &{=}\lim_{t\rightarrow \infty} P e^{Jt} P^{-1} \\
&=P \lim_{t\rightarrow \infty}  \left[\begin{smallmatrix}
1 & t & 0_{1\times (2n-2)} \\
0 & 1 & 0_{1\times (2n-2)} \\
0_{ (2n-2)\times 1} & 0_{ (2n-2)\times 1} & e^{J't}
\end{smallmatrix}\right] P^{-1}  \\
&=\lim_{t\rightarrow \infty} \frac 1n \left[\begin{smallmatrix}
1_{n\times n} & t1_{n\times n} \\
0_{n\times n} & 1_{n\times n} 
\end{smallmatrix}\right],
\end{align*}
where $J$ is a Jordan matrix. 
Thus, given an initial position $x(0)=x_0$, we obtain 
$
\lim_{t\rightarrow \infty} x_i(t) =  \frac{1}{n} \sum_{i\in\mathcal{V}}  x_{0,i} \quad \forall \; i \in \mathcal{V}
$, 
i.e., the agents converge to the average of their initial positions.

We now consider the case where $\delta=0$ and $d_i \ne 0$ for at least one $i\in\mathcal{V}$. 
Define the output
\begin{align*}
\left[\begin{smallmatrix}
 y_x
\end{smallmatrix}\right] &= {\left[\begin{smallmatrix}
0_{m\times n}  &  \mathcal{B}^T
\end{smallmatrix}\right]}   \left[\begin{smallmatrix}
{z}\\{x}
\end{smallmatrix}\right],
\end{align*}
and consider the linear coordinate change:
\begin{align}
\begin{aligned}
\label{eq:proof_pi_single_coordinate_change}
x &= \left[\begin{smallmatrix}
 \frac{1}{\sqrt{n}} 1_{n\times 1} & S 
\end{smallmatrix}\right] u \quad 
u = \left[\begin{smallmatrix}
 \frac{1}{\sqrt{n}} 1_{1\times n} \\ S^T 
\end{smallmatrix}\right] x \\ 
z &= \left[\begin{smallmatrix}
 \frac{1}{\sqrt{n}} 1_{n\times 1} & S
\end{smallmatrix}\right] w 
 \quad
w = \left[\begin{smallmatrix}
 \frac{1}{\sqrt{n}} 1_{1\times n} \\ S^T
\end{smallmatrix}\right] z 
\end{aligned}
\end{align}
where $S$ is a matrix such that $[
 \frac{1}{\sqrt{n}} 1_{n\times 1}, S 
]$ is an orthonormal matrix.
In the new coordinates, the system dynamics become:
\begin{align}
\dot{w} &= u \nonumber \\
\dot{u} &= \left[\begin{smallmatrix}
0 & 0_{1\times (n-1)} \\
0_{(n-1)\times 1} & -a S^T\mathcal{L}S 
\end{smallmatrix}\right] w + \left[\begin{smallmatrix}
0 & 0_{1\times (n-1)} \\
0_{(n-1)\times 1} & -b S^T\mathcal{L}S 
\end{smallmatrix}\right] u \label{eq:proof_pi_single_2} \\
&+ \left[\begin{smallmatrix}
\frac {1}{\sqrt{n}}{1_{1\times n}} \\ S^T
\end{smallmatrix}\right] d. \nonumber
\end{align}
We note that the states $u_1$ and $w_1$ are both unobservable and uncontrollable. We thus omit these states to obtain a minimal realization by defining the new coordinates $u'=[
u_2, \hdots , u_n
]^T$ and $w'=[
w_2, \hdots , w_n
]^T$, thus obtaining the system dynamics
\begin{align*}
\left[\begin{smallmatrix}
\dot{w}'\\ \dot{u}' 
\end{smallmatrix}\right]
 &= \left[\begin{smallmatrix}
0_{(n-1)\times (n-1)} & I_{(n-1)} \\
 -a S^T\mathcal{L}S &  -b S^T\mathcal{L}S
\end{smallmatrix}\right]  
\left[\begin{smallmatrix}
{w}'\\ {u}' 
\end{smallmatrix}\right]
+ \left[\begin{smallmatrix}
 0_{(n-1)\times 1} \\  S^T d
\end{smallmatrix}\right].
\end{align*}
Clearly $x^TS^T\mathcal{L}Sx\ge 0$, with equality only if $Sx=k1_{n\times 1}$. However, since $[\frac{1}{\sqrt{n}} 1_{n\times 1}, S ]$ is orthonormal, $1_{1\times n}Sx=0_{1\times n}x=0=k1_{1\times n}1_{n\times 1} = kn$, which implies $k=0$. 
Hence $S^T\mathcal{L}S$ is positive definite and thus invertible, and we may define
\begin{align*}
\left[\begin{smallmatrix}
{w}''\\ {u}'' 
\end{smallmatrix}\right]
&= \left[\begin{smallmatrix}
{w}'\\ {u}' 
\end{smallmatrix}\right] 
- \left[\begin{smallmatrix}
\frac 1a (S^T\mathcal{L}S)^{-1} S^T d  \\ 0_{(n-1)\times 1} 
\end{smallmatrix}\right].
\end{align*}
It is easily verified that the origin is the only equilibrium of the system dynamics, which in the new coordinates are given by
\begin{align*}
\left[\begin{smallmatrix}
\dot{w}''\\ \dot{u}''
\end{smallmatrix}\right]
 &= \underbrace{\left[\begin{smallmatrix}
0_{(n-1)\times (n-1)} & I_{(n-1)} \\
 -a S^T\mathcal{L}S &  -b S^T\mathcal{L}S
\end{smallmatrix}\right]}_{\triangleq A}  
\left[\begin{smallmatrix}
{w}''\\ {u}''
\end{smallmatrix}\right].
\end{align*}
By elementary column operations, the characteristic polynomial in $\kappa$ of $A$ is given by
$ \det( \kappa^2 I_{(n-1)} + (b \kappa + a) S^T\mathcal{L}S )
$. By comparing this polynomial with the characteristic polynomial $\det( s I + S^T \mathcal{L} S ) $, which since $S^T \mathcal{L} S$ is positive definite has solutions $-s_i<0$, we know that the eigenvalues of $A$ must satisfy
$ \kappa^2 + s_ib \kappa + s_ia =0
$, with solutions $\kappa \in \mathbb{C}^-$. Thus $A$ is Hurwitz.
From the dynamics \eqref{eq:proof_pi_single_2}, it is clear that $\dot{u}_1= \frac 1{\sqrt{n}}{1_{1\times n}} d $. Hence $\lim_{t \rightarrow \infty} u_1(t) = \pm \infty$ unless ${1_{1\times n}} d =0$. Since ${u}'$ is bounded, by the coordinate change \eqref{eq:proof_pi_single_coordinate_change}, $x$ is bounded iff ${1_{1\times n }} d =0$. 

Now consider the case where $\delta> 0$ and $d_i = 0 \; \forall i\in\mathcal{V}$. The dynamics can be written as:
\begin{align}
\left[\begin{smallmatrix}
\dot{z}\\\dot{x}
\end{smallmatrix}\right] = {\left[\begin{smallmatrix}
 0_{n\times n} & I_{n } \\
-a \mathcal{L} & -b \mathcal{L} -\delta I
\end{smallmatrix}\right]} 
\left[\begin{smallmatrix}
{z}\\{x}
\end{smallmatrix}\right] +
\left[
\begin{smallmatrix}
0_{n\times 1} \\
\delta x(0)
\end{smallmatrix}
\right]
.
\label{eq:proof_pi_single_delta_ne_0}
\end{align}
Define the output of the system
\begin{align*}
\left[\begin{smallmatrix}
 y_x
\end{smallmatrix}\right] &= {\left[\begin{smallmatrix}
0_{m\times n}  &  \mathcal{B}^T
\end{smallmatrix}\right]}   \left[\begin{smallmatrix}
{z}\\{x}
\end{smallmatrix}\right],
\end{align*}
and consider the linear coordinate change:
\begin{align}
\label{eq:proof_pi_single_coordinate_change_2}
z = \left[\begin{smallmatrix}
 \frac{1}{\sqrt{n}} 1_{n\times 1} & S
\end{smallmatrix}\right] w 
 \quad
w = \left[\begin{smallmatrix}
 \frac{1}{\sqrt{n}} 1_{1\times n} \\ S^T
\end{smallmatrix}\right] z ,
\end{align}
where $S$ is a matrix such that $[
 \frac{1}{\sqrt{n}} 1_{n\times 1}, S 
]$ is an orthonormal matrix.
In the new coordinates, the system dynamics \eqref{eq:proof_pi_single_delta_ne_0} become:
\begin{align*}
\left[\begin{smallmatrix}
\dot{w}\\\dot{x}
\end{smallmatrix}\right] = {\left[\begin{smallmatrix}
 0_{n\times n} & \left[\begin{smallmatrix}
 \frac{1}{\sqrt{n}} 1_{1\times n} \\ S^T
\end{smallmatrix} \right]  \\
\left[\begin{smallmatrix}
0_{n\times 1} & -a \mathcal{L} S
\end{smallmatrix} \right]
 & -b \mathcal{L} -\delta I
\end{smallmatrix}\right]} 
\left[\begin{smallmatrix}
{w}\\{x}
\end{smallmatrix}\right] +
\left[
\begin{smallmatrix}
0_{n\times 1} \\
\delta x(0)
\end{smallmatrix}
\right]
.
\end{align*}
We note that the state $w_1$ is unobservable and uncontrollable. We thus omit this state to obtain a minimal realization by defining the new coordinates $w'=[
w_2, \hdots , w_n
]^T$, thus obtaining the system dynamics
\begin{align}
\left[\begin{smallmatrix}
\dot{w}'\\\dot{x}
\end{smallmatrix}\right] =
\underbrace{\left[
\begin{smallmatrix}
0_{(n-1)\times (n-1)} & S^T \\
-a\mathcal{L}S & -b \mathcal{L}  -\delta I_n
\end{smallmatrix}
\right]}_{\triangleq A}
\left[\begin{smallmatrix}
{w'}\\{x}
\end{smallmatrix}\right]
 +
\left[
\begin{smallmatrix}
0_{n\times 1} \\
\delta x(0)
\end{smallmatrix}
\right] 
.
\label{eq:proof_pi_single_delta_ne_0_prime}
\end{align}

By elementary column operations, the characteristic polynomial of $A$ may be written as
$0= \det \left( s^2 I_n + s(b\mathcal{L}+\delta I_n) +a\mathcal{L}SS^T \right) \triangleq \det(Q(s))
$. For a given $s\in \mathbb{C}$, the previous equation has a solution only if $x^TQ(s)x=0$ for some $x$ satisfying $x^Tx=1$. This equation becomes
\begin{align*}
0=   s^2 \underbrace{x^Tx}_{a_2} + s\underbrace{x^T(b\mathcal{L}+\delta I_n)x}_{a_1} +\underbrace{ax^T\mathcal{L}SS^Tx}_{a_0}.
\end{align*}
By the Routh-Hurwitz stability criterion, the above equation has all its solutions $s\in\mathbb{C}^-$ if and only if $a_i>0$ for $i=0,1,2$. Clearly $a_2=1$ by definition, while $a_1>0$ by the positive definiteness of $\mathcal{L}+\delta I_n$. It is easily shown that $SS^T=I_n-1/n1_{n\times n}$, by the orthonormality of $[
 \frac{1}{\sqrt{n}} 1_{n\times 1}, S 
]$. This implies that $\mathcal{L}SS^T=\mathcal{L}$, implying that $a_0\ge 0$. If $a_0>0$, all solutions are stable, whereas if $a_0=0$, $s=0$ is also a solution.  Thus, any eigenvalue of $A$ must be either zero, or have negative real part. However, it is easily verified that $A$ is full rank. Thus $0$ cannot be an eigenvalue of $A$, implying that $A$ is Hurwitz. 
 The first $n-1$ rows of the equilibrium of \eqref{eq:proof_pi_single_delta_ne_0_prime} imply $S^Tx=0_{(n-1)\times 1}$, which implies $x=x^* 1_{n\times 1}$. Finally, the last $n$ rows of \eqref{eq:proof_pi_single_delta_ne_0_prime} imply
$ 1_{1\times n} \left( -a \mathcal{L} Sw' -b\mathcal{L}x^* 1_{n\times 1} -\delta x^* 1_{n\times 1} + \delta x(0) \right) = 0    
$, so $ n x^* = \sum_{i\in\mathcal{V}} x_i(0)$


The case when $\delta> 0$ and $d_i\ne 0$ for at least one $i\in \mathcal{V}$ is analogous to the case when $\delta= 0$ and $d_i\ne 0$, and the proof is thus omitted.
\end{proof}

\subsection{Consensus by distributed integral action for double-integrator dynamics with damping}
\label{sec:consensus_integral_double}
Consider agents with double-integrator dynamics \eqref{eq:dynamics_double_2}, with input given by the velocity-damped PI-controller:
\begin{align}
\begin{aligned}
u_i &= -\sum_{j\in \mathcal{N}_i} \left( b(x_i{-}x_j) +  a \int_0^t (x_i(\tau){-}x_j(\tau)) \ud \tau \right) \\
&-\gamma v_i -\delta (x_i{-}x_i^0)  
\end{aligned}
 \label{eq:consensus_absolute_3}
\end{align}
where $a \in \mathbb{R}^+$, $b \in \mathbb{R}^+$, $\gamma \in \mathbb{R}^+$, $\delta \in  \mathbb{R}^+$ and $d_i \in \mathbb{R}$ is an unknown scalar disturbance. The above protocol does not require communication of the integral state between the agents, as it suffices for each agent to measure its neighbors' states and integrate the relative differences. 
\begin{theorem}
\label{th:double_abs}
Under the dynamics \eqref{eq:dynamics_double_2} with $u_i$ given by \eqref{eq:consensus_absolute_3},
the agents converge to a common value $x^*$ for any constant disturbance $d_i$, provided that $a<b\gamma$. If $d_i=0\;\forall i\in\mathcal{V}$, the agents converge to $x^*= \frac 1n \sum_{i\in \mathcal{V}} x_i(0)$ for arbitrary $v_i(0)$.
If absolute position measurements are not present, i.e., $\delta=0$, it still holds that $\lim_{t \rightarrow \infty } |x_i(t)-x_j(t)| =0 \; \forall i,j \in \mathcal{V}$ for any set of disturbances $d_i$. However the absolute states are in general unbounded, i.e., $\lim_{t\rightarrow \infty} |x_i(t)| = \infty \; \forall i \in \mathcal{V}$, unless ${1_{1\times n}} d =0$. Also, in this case the system is stable if and only if $a<b\gamma$. 
\end{theorem}

\begin{proof} 
The proof follows the same principle ideas as the proof of Theorem \ref{th:single_abs}. However, as second-order dynamics are considered, the problem is inherently different from first-order dynamics. 
First consider the case where $\delta = 0$. 
 Let also $d_i=0\; \forall i\in \mathcal{V}$.
By introducing the state vector $ z=[z_1, \hdots , z_n]^T$ we may rewrite the dynamics as:
\begin{align*}
\left[\begin{smallmatrix}
\dot{z}\\\dot{x} \\ \dot{v}
\end{smallmatrix}\right] = \underbrace{\left[\begin{smallmatrix}
 0_{n\times n} & I_n &  0_{n\times n}  \\
 0_{n\times n} &  0_{n\times n} & I_n  \\
-a \mathcal{L} & -b \mathcal{L} & -\gamma I_n 
\end{smallmatrix}\right]}_{ \triangleq A}
\left[\begin{smallmatrix}
{z}\\{x} \\ v
\end{smallmatrix}\right],
\end{align*}
together with the additional initial condition $z(0)=0_{n\times 1}$. 
By elementary column operations it is easily shown that the characteristic polynomial of $A$ can be written as
$0  =  \det ( (a + b s) \mathcal{L} +  s^2 (s+\gamma) I )
$, where $I$ is the identity matrix of appropriate dimensions. Comparing the above equation with the characteristic polynomial of $\mathcal{L}$, we get that
$0=s^3 + \gamma s^2 + \lambda_i b s + \lambda_i a
$, where $\lambda_i$ is an eigenvalue of $\mathcal{L}$. If $\lambda_i>0$, the above equation has all its solutions $s \in \mathbb{C}^{-}$ if and only if $a< b \gamma$, and $a,b,\gamma>0$ by the Routh-Hurwitz stability criterion. Since $\mathcal{G}$ by assumption is connected, $\lambda_1=0$ and $\lambda_i>0 \; \forall i=2, \hdots, n$. For $\lambda_1=0$, the above equation has the solutions $s=0, s=-\gamma$. By straightforward calculations it can be shown that 
 $e_1^1= [
1_{1\times n} ,\; 0_{1\times n} ,\ 0_{1\times n}
]^T$ and $e_1^2= [
0_{1\times n} ,\; 1_{1\times n} ,\; 0_{1\times n}
]^T$ are an eigenvector and a generalized eigenvector of $A$, respectively,  corresponding to the eigenvalue $0$. 
Furthermore $v_1={1}/{(\gamma^2 n)} [
\gamma^2 1_{1\times n} ,\; 0_{1\times n} ,\;  - 1_{1\times n} 
]$ and  $v_2={1}/{(\gamma n)} [
0_{1\times n} ,\; \gamma 1_{1\times n} ,\;   1_{1\times n} 
]$ are a generalized left eigenvector and a left eigenvector of $A$ corresponding to the eigenvalue $0$.
Furthermore $v_1 e_1^1 =1, \; v_2 e_1^2 = 1$ and $v_2 e_1^1 =0, \; v_1 e_1^2 = 0$.
Hence the first columns of $P$ can be chosen as $e_1^1$ and $e_1^2$, and the first rows of $P^{-1}$ can be chosen to be $v_1$ and $v_2$.
 Since all other eigenvalues of $A$ have strictly negative real part we obtain 
\begin{align*}
&\lim_{t\rightarrow \infty} e^{At} = \lim_{t\rightarrow \infty} P e^{Jt} P^{-1} \\ 
&= P \lim_{t\rightarrow \infty} \left[\begin{smallmatrix}
1 & t & 0_{1\times (3n-2)} \\
0 & 1 & 0_{1\times (3n-2)} \\
0_{ (3n-2)\times 1} & 0_{ (3n-2)\times 1} & e^{J't}
\end{smallmatrix}\right] P^{-1}  \\
&= \lim_{t\rightarrow \infty} P \left[\begin{smallmatrix}
1 & t & 0_{1\times (3n-2)} \\
0 & 1 & 0_{1\times (3n-2)} \\
0_{ (3n-2)\times 1} & 0_{ (3n-2)\times 1} & 0_{ (3n-2)\times (3n-2)}
\end{smallmatrix}\right] P^{-1}  \\
 &= \lim_{t\rightarrow \infty} \frac 1n \left[\begin{smallmatrix}
1_{n\times n} & t1_{n\times n} & \frac{t\gamma - 1}{\gamma^2}  1_{n\times n}  \\
0_{n\times n} & 1_{n\times n} &\frac{1}{\gamma}1_{n\times n}  \\
0_{n\times n} &0_{n\times n} &0_{n\times n} 
\end{smallmatrix}\right]
\end{align*}
where $J$ is a Jordan matrix. 
Given any initial position $x(0)=x_0, z(0)=0_{n\times 1}, \; v(0)=v_0$, we obtain 
$\lim_{t\rightarrow \infty} x_i(t) =  \frac{1}{n} \sum_{i\in\mathcal{V}}  x_{0,i} +  \frac{1}{\gamma n} \sum_{i\in\mathcal{V}}  v_{0,i} \quad \forall \; i \in \mathcal{V}
$. 

Now let us turn our attention to the case where $\delta> 0$ and $d_i\ne 0$ for at least one $i\in \mathcal{V}$.  We define the output of the system as 
\begin{align*}
\left[\begin{smallmatrix}
 y_x \\ y_v
\end{smallmatrix}\right] &= {\left[\begin{smallmatrix}
0_{m\times n}  &  \mathcal{B}^T & 0_{m\times n}\\
0_{m\times n}  & 0_{m\times n}  &  \mathcal{B}^T
\end{smallmatrix}\right]}
\left[\begin{smallmatrix}
{z}\\{x} \\ v
\end{smallmatrix}\right],
\end{align*}
and consider the same linear coordinate change of $z$, $x$ and $v$ as applied to $z$ and $x$ in the proof of Theorem \ref{th:single_abs}. In the new coordinates the system dynamics are
\begin{align}
\begin{aligned}
\dot{z}' &= x'  \\
\dot{x}' &= v'  \\
\dot{v}' &= \left[\begin{smallmatrix}
0 & 0_{1\times (n-1)}  \\
0_{(n-1)\times 1} & -a S^T\mathcal{L}S 
\end{smallmatrix}\right] z' +  \left[\begin{smallmatrix}
0 & 0_{1\times (n-1)} \\
0_{(n-1)\times 1} & -b S^T\mathcal{L}S 
\end{smallmatrix}\right] x' \\
&- \gamma v'
+ \left[\begin{smallmatrix}
\frac 1n{1_{1\times(n)}} \\ S^T
\end{smallmatrix}\right] d.
\end{aligned}
\label{eq:proof_pi_double_2}
\end{align}
We note that the states $z'_1$, $x'_1$ and $v'_1$ are unobservable and uncontrollable. We thus omit these states to obtain a minimal realization by defining the new coordinates $z''=[
z'_2, \hdots , z'_n
]^T$, 
$x''=[
x'_2, \hdots , x'_n
]^T$ and 
$v''=[
v'_2, \hdots , v'_n
]^T$
 we obtain the system dynamics
\begin{align*}
\left[\begin{smallmatrix}
\dot{z}''\\ \dot{x}'' \\ \dot{v}'' 
\end{smallmatrix}\right]
 &= \underbrace{\left[\begin{smallmatrix}
0_{(n-1)^2} & I_{(n-1)^2} & 0_{(n-1)^2} \\
0_{(n-1)^2} & 0_{(n-1)^2} &  I_{(n-1)^2} \\
 -a S^T\mathcal{L}S &  -b S^T\mathcal{L}S & -\gamma  I_{(n-1)^2}
\end{smallmatrix}\right]}_{ \triangleq A'} 
\left[\begin{smallmatrix}
{z}''\\ {x}'' \\ {v}'' 
\end{smallmatrix}\right] + \left[\begin{smallmatrix}
 0_{(n-1)\times 1} \\ 0_{(n-1)\times 1} \\ S^T d 
\end{smallmatrix}\right].
\end{align*}
We now shift the state space by defining
\begin{align*}
\left[\begin{smallmatrix}
{z}^{(3)}\\ {x}^{(3)} \\ {v}^{(3)}
\end{smallmatrix}\right]
&=
\left[\begin{smallmatrix}
{z}''\\ {x}'' \\ {v}'' 
\end{smallmatrix}\right]
- \left[\begin{smallmatrix}
\frac 1a (S^T\mathcal{L}S)^{-1} S^T d  \\ 0_{(n-1)\times 1} \\ 0_{(n-1)\times 1} 
\end{smallmatrix}\right].
\end{align*}
It is easily verified that the origin is the only equilibrium of the system dynamics, and that the stability in the new coordinates is characterized by the matrix $A'$.
By a similar argument used when showing that $A$ has eigenvalues with non-positive real part, we may show that $A'$ has eigenvalues with non-positive real part. But since $S^T\mathcal{L}S$ is full-rank, $A'$ must also be full-rank, and hence $A'$ is Hurwitz. 
Thus $\lim_{t\rightarrow \infty} x^{(3)} = \lim_{t\rightarrow \infty} x'' = 0_{(n-1)\times 1}$, which implies that $\lim_{t\rightarrow \infty} |x_i(t)-x_j(t)| = 0 \; \forall i,j \in \mathcal{V} $, even in the presence of disturbances $d_i$. It is also clear that whenever $a \ge b \gamma$, at least one eigenvalue will have non-negative real part, and that its (generalized) eigenvector will be distinct from $e_1$ and $e_2$. 
From the dynamics \eqref{eq:proof_pi_double_2}, it is clear that $\dot{x}'_1= \frac 1n{1_{1\times n}} d $. Hence $\lim_{t \rightarrow \infty} x'_1(t) = \pm \infty$ unless ${1_{1\times n}} d =0$. Since ${x}''$ is bounded, by the coordinate change \eqref{eq:proof_pi_single_coordinate_change}, $x$ is bounded if and only if ${1_{1\times n }} d =0$. 

The stability analysis of the case when $\delta>0$ is analogous to the corresponding part of the proof of Theorem \ref{th:single_abs}, and hence omitted. 
If $d_i=0 \; \forall i\in\mathcal{V}$, stationarity of $v(t)$ implies: \\
$\lim_{t\rightarrow\infty} 1_{1\times n} \left( -a \mathcal{L} z -b\mathcal{L}x -\delta x + \delta x(0) -\gamma v \right) = 0$, so $n x^* = \sum_{i\in\mathcal{V}} x_i(0)$.
\end{proof}

\section{Motivating applications revisited}
\label{sec:applications}

\subsection{Thermal energy storage in smart buildings}
We here return to the example of thermal energy storage in smart buildings, introduced in Section \ref{subsec:green_intro}. Recall that the  temperatures dynamics in the rooms can be described by \eqref{eq:ex_building_temperature}. 

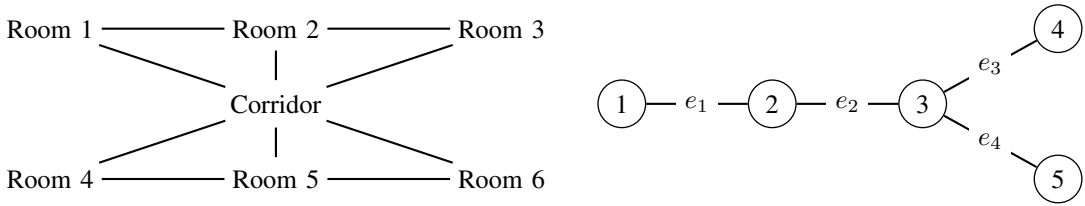
\begin{figure*}[t]
\begin{center}
$
\begin{array}{cc}
\raisebox{-0.5\height}{
\begin{tikzpicture}
\GraphInit[vstyle=empty]
\Vertex[x=-3, y=1] {Room 1}
\Vertex[x=0, y=1] {Room 2}
\Vertex[x=3, y=1] {Room 3}
\Vertex[x=0, y=0] {Corridor}
\Vertex[x=-3, y=-1] {Room 4}
\Vertex[x=0, y=-1] {Room 5}
\Vertex[x=3, y=-1] {Room 6}
\Edge[](Room 1)(Room 2) 
\Edge[](Room 2)(Room 3) 
\Edge[](Room 1)(Corridor) 
\Edge[](Room 2)(Corridor) 
\Edge[](Room 3)(Corridor) 
\Edge[](Room 4)(Corridor) 
\Edge[](Room 5)(Corridor) 
\Edge[](Room 6)(Corridor) 
\Edge[](Room 4)(Room 5) 
\Edge[](Room 5)(Room 6) 
\end{tikzpicture}
}
&
\raisebox{-0.5\height}{
\begin{tikzpicture}
\GraphInit[vstyle=Normal]
\Vertex[x=-2, y=0] {1} 
\Vertex[x=0, y=0] {2}
\Vertex[x=2, y=0] {3}
\Vertex[x=3.8, y=1] {4}
\Vertex[x=3.8, y=-1] {5}
\Edge[label=$e_1$](1)(2) 
\Edge[label=$e_3$](3)(4) 
\Edge[label=$e_2$](3)(2)
\Edge[label=$e_4$](3)(5)   
\end{tikzpicture}
}
\end{array}
$
\end{center}
\caption{The leftmost figure illustrates the floor topology. The rightmost figure shows the communication topology of the space satellites.}
\label{fig:graph_corridor}
\end{figure*}
\begin{figure}[t]
\begin{center}
$
\begin{array}{c}
\input{Simulations/hc.tikz} \\ \input{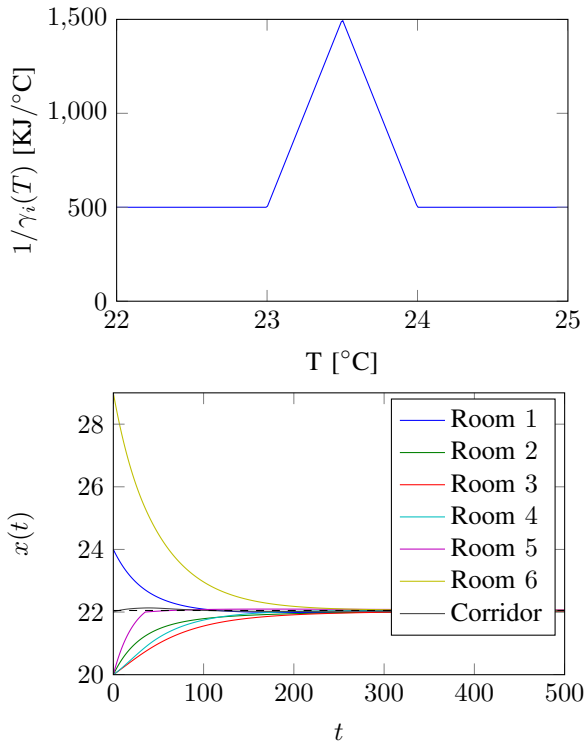} 
\end{array}
$
\end{center}
\caption{The top figure shows the heat capacities of Room 2 and 5. The bottom figure shows the temperatures in the building floor.}
\label{fig:smartbuilding}
\end{figure}

The heat conductivity $a$ is assumed to be constant and uniform, implying $a_{ij}(x)=a x \; \forall (i,j) \in \mathcal{E}$, where $a=0.5 \text{W/$^{\circ}$C}$. 
Consider the floor topology in Figure \ref{fig:graph_corridor}. 
Let the desired maximum temperature be given by $t_b=23^{\circ}$C. The heat capacity is assumed to be given by Figure \ref{fig:smartbuilding}  
 for $i\in \{\text{Room 2},\text{Room 5}\}$ due to thermal energy storage installations, and  ${1}/{\gamma_i(T)}= 50 \text{kJ/$^{\circ}$C}$ for $i \in \{ \text{Room 1},\text{Room 3}, \text{Room 4}, \text{Room 6}, \text{Corridor} \}$ where no thermal energy storage is installed. 
 The initial temperatures was assumed to be $29 ^{\circ}$C for room $6$, $24 ^{\circ}$C for room $1$, $22 ^{\circ}$C for the corridor and $20 ^{\circ}$C for the other rooms. 
The temperatures as a function of time are shown in Figure \ref{fig:smartbuilding} for given initial temperatures. We note that the temperature in room $2$ and $5$ never exceeds the desired maximum temperature $t_b=23 ^{\circ}$C due to the thermal energy storage, and that the temperatures converge to a temperature below $t_b$ in all rooms, in accordance with Theorem \ref{prop:1}. 


	\setlength\fheight{3.75cm} 
	\setlength\fwidth{6.0cm}

\subsection{Autonomous space satellites}
\label{sec:example_underwater}
\begin{figure*}[t]
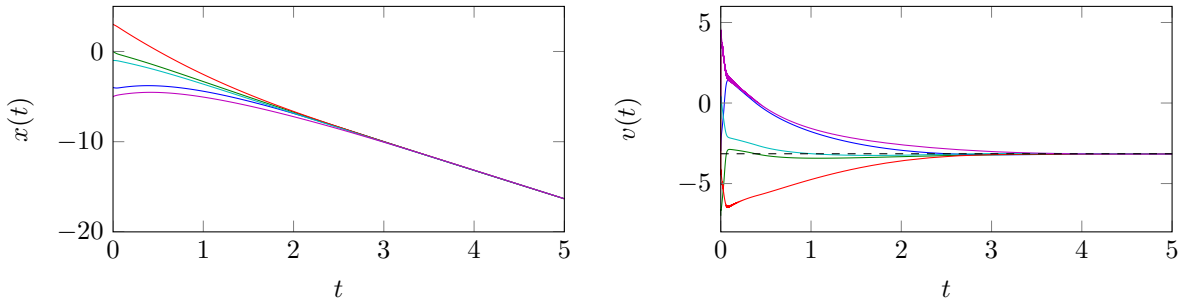

	\centering
$
\begin{array}{cc}
	\input{Simulations/secondorder_sim1_1.tikz}
 & \input{Simulations/secondorder_sim1_2.tikz}
\end{array}
$
\caption{The leftmost figure shows the positions of the satellites for the initial conditions  $x(0) = [-4, 0, 3, -1, -5]^T, v(0)=[-3,    -7,     3,    -1,     0 ]^T$, while their velocities.} 
\label{fig:secondorder_sim_1}
\end{figure*}

Consider a group of autonomous space satellites with unitary masses.  
The satellites are denoted as $1, \hdots, 5$, and their communication topology is given by the undirected graph in Figure \ref{fig:graph_corridor}. 
The objective is to reach consensus in one dimension by a distributed control law using only relative position and velocity measurements. 
 The raw control signal is the power applied to the agent's engine, $P_i$. However, the acceleration in an inertial reference frame is $a_i={P_i}/{|v_i|}$ due to $P_i= \langle F_i,v_i \rangle $ and $P_i$ being parallel to the agent's velocity $v_i$. We assume that the agents only have access to relative measurements, and hence are unaware of their absolute positions. 
 Assuming that $a_i={P_i}/{(|v_i|+c)}$, $c>0$, to ensure that the accelerations remain bounded, this scenario can be modelled by \eqref{nnlsecorder3}. This is clearly a special case of the dynamics \eqref{eq:dynamics_double_2} with $u_i$ given by \eqref{nnlsecorder3}. 
The interaction functions in this example are chosen to be
$ a_{ij}(y) = 2 b_{ij}(y) =  20 ( e^{|y|} -1 )\sgn{(y)} \quad \forall  (i,j) \in \mathcal{E} 
$, which satisfies Assumption \ref{ass:alpha}. It is clear that this situation cannot be modelled by any previously proposed linear consensus protocols.
Figure \ref{fig:secondorder_sim_1} shows the state trajectories for different initial conditions. 
As predicted by Theorem \ref{prop:nlsecorder}, consensus is reached, and the final consensus velocity, as seen from an observer, can be calculated by \eqref{eq:consensus_point}. 

	\setlength\fheight{3.0cm} 
	\setlength\fwidth{6.0cm}

\subsection{Mobile robot coordination under disturbances}
\label{ex:1}
In this section we revisit the example of mobile robots from Section \ref{subsec:ex_mobile_robots}. The dynamics of the robots are given by \eqref{eq:ex_mobile_robots}. 
Let the damping coefficient be given by $\gamma=3$, and the static gain $b=5$. 
 We consider the system with a constant disturbance  $d=[1,0,0,0,0]^T$, and for the different integral gains $a=0$, $a=1$, and $a=15$. The initial conditions are given by $x(0)=[5, -6, 8, 4, 5]^T$, $v(0)=[0, 0, 0, 0, 0 ]^T$. The setup we will consider consists of a string of 5 mobile robots, whose communication topology is a string graph. 

\setlength\fheight{3.0cm} 
\setlength\fwidth{6.0cm}

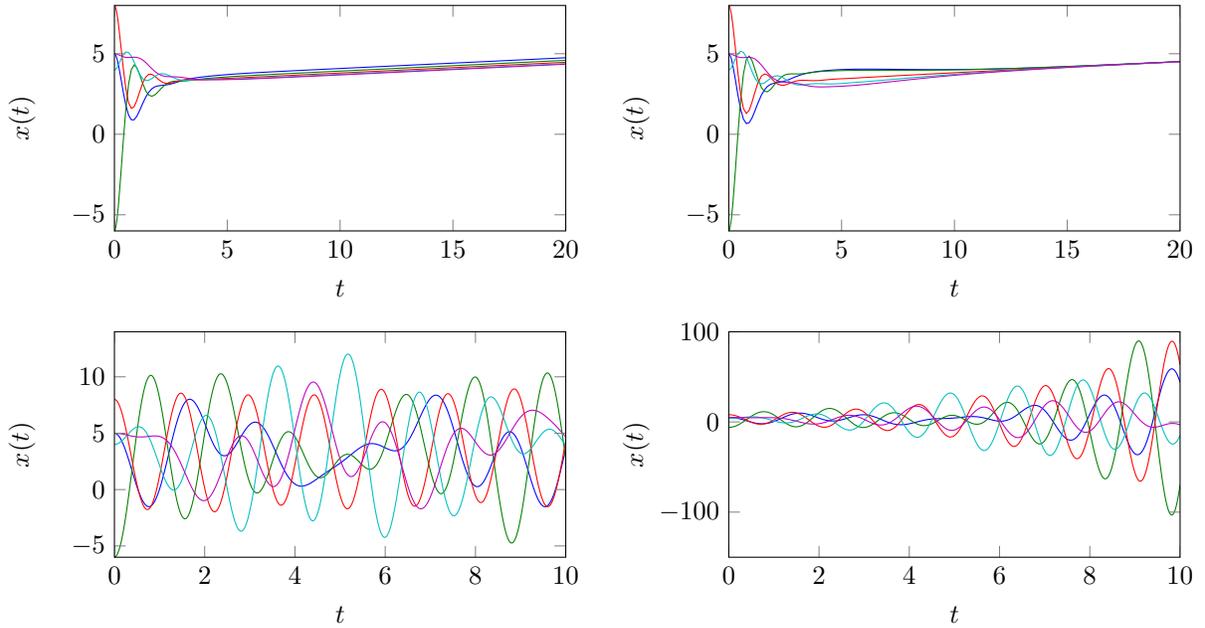
\begin{figure*}[t]
	\centering
$
\begin{array}{cc}
	\input{Simulations/consensus_without_I.tikz} & 	\input{Simulations/consensus_with_I.tikz} \\
	\input{Simulations/consensus_with_I_too_large.tikz} & 	\input{Simulations/consensus_with_way_I_too_large.tikz} 
\end{array} 
$
\caption{The upper left figure shows the state trajectories of \eqref{eq:ex_mobile_robots} when $a=0$, the upper right figure shows the state trajectories when $a=1$, the lower left figure shows the state trajectories when $a=15$, and the lower right figure shows the state trajectories when $a=20$.} 
\label{fig:secondorder_sim_1_no_disturbance}
\end{figure*}

\setlength{\fwidth}{\columnwidth-24mm}
\setlength{\fheight}{0.618\fwidth} 
By Theorem \ref{th:double_abs} stability is guaranteed if and only if $a < b\gamma$. In Figure~\ref{fig:secondorder_sim_1_no_disturbance}, the state trajectories are shown for different choices of $a$. We observe that asymptotic consensus amongst the mobile robots is only reached for $a=1$. For $a=0$, consensus is not reached due to the presence of a static disturbance. When $a=1$, the disturbance is attenuated by the integrators, and asymptotic consensus is reached. However, as $a$ is increased to $15=b\gamma$, the system becomes marginally stable, i.e., stable but not asymptotically stable. By increasing $a$ further to $20$, the system becomes unstable, in accordance with Theorem \ref{th:double_abs}.

\subsection{Frequency control of power systems}
\label{sec:pow}
In this section we demonstrate that a similar protocol to the one proposed in Section \ref{sec:consensus_integral_double}, can be employed for frequency control of power systems. Let us consider a power system, whose topology modeled by a graph $\mathcal{G}=(\mathcal{V}, \mathcal{E})$. Each node, here referred to as a bus, is assumed to obey the linearized swing equation \eqref{eq:swing_intro}. 
 By defining $\delta = [
\delta_1 \hdots , \delta_n
]^T$, we may rewrite \eqref{eq:swing_intro} as:
\begin{eqnarray}
\label{eq:swing_vector}
\left[\begin{smallmatrix}
\dot{\delta} \\ \dot{\omega}
\end{smallmatrix}\right] = \left[\begin{smallmatrix}
0_{n\times n} & I_{n} \\
-M \mathcal{L}_k & - M D
\end{smallmatrix}\right] 
\left[\begin{smallmatrix}
{\delta} \\ {\omega}
\end{smallmatrix}\right] +
\left[\begin{smallmatrix}
0_{n\times 1} \\ M p^m
\end{smallmatrix}\right] +
\left[\begin{smallmatrix}
0_{n\times 1} \\ M u
\end{smallmatrix}\right],
\end{eqnarray}
where $M=\diag({1}/{m_1}, \hdots , {1}/{m_n})$, $D=\diag(d_1, \hdots, d_n)$, $\mathcal{L}_k$ is the weighted Laplacian with edge weights $k_{ij}$,
 $p^m=[
p^m_1,\hdots, p^m_n
]^T$ and $u=[
u_i,\hdots, u_n
]^T$. 

\subsubsection{Centralized PI control}
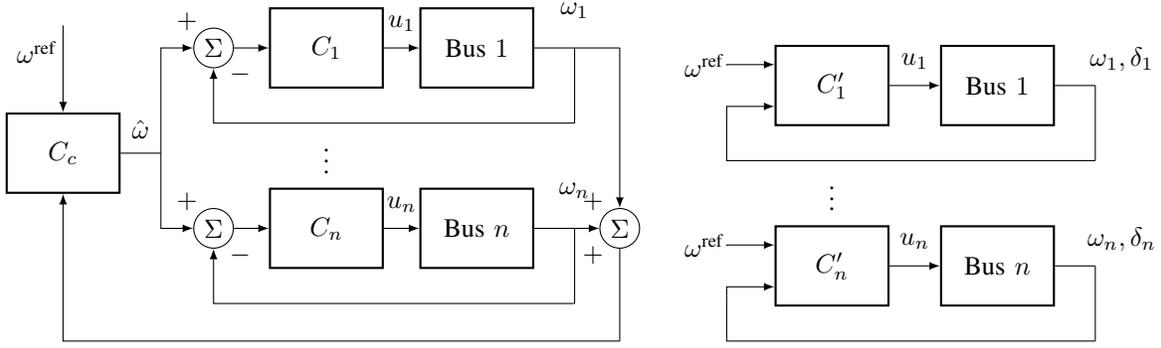
\begin{figure*}[t]
\begin{center}
$
\begin{array}{cc}
\begin{tikzpicture}[auto, >=latex] 
    	
   	\node [input, name=ref] {};
   	\node [block, above of=ref, node distance =1.0cm, minimum width=1.5cm] (p1) {$C_1$};
	\node [draw=none,fill=none, below of =p1,node distance = 1.4cm, minimum width=1.0cm] (dots) {$\vdots$};
   	\node [block, below of =dots,node distance = 1.0cm, minimum width=1.5cm] (pn) {$C_n$};
   	\node [block, right of=p1, node distance =2.0cm, minimum width=1.5cm] (b1) {Bus $1$};
   	\node [block, right of=pn, node distance =2.0cm, minimum width=1.5cm] (bn) {Bus $n$};
   	\node [block, left of=dots, node distance = 3.5cm, minimum width=1.5cm] (pc) {$C_c$};
   	\coordinate [input, right of=pc, node distance =1.3cm, minimum width=1.5cm] (connector1) {};
	\node [output, name=output1, right of=b1, node distance=1.3cm]{};
    	\node [output, name=outputn, right of=bn, node distance=1.3cm]{};
    	\node [output, name=node3, below of=output1, node distance=1.0cm]{};
    	\node [output, name=node4, below of=outputn, node distance=1.0cm]{};
    	\node [output, name=node5, below of=node3, node distance=1.0cm]{};	
    	\node [output, name=node6, right of=output1, node distance=0.6cm]{};	
	\node [sum, name=s2, right of=outputn, node distance=0.6cm] {$\Sigma$};
    	\node [output, name=node8, below of=s2, node distance=1.5cm]{};	
    	\node [output, name=node9, above of=pc, node distance=1.7cm]{};	

	\node [sum, name=s1, left of = p1, node distance = 1.5cm] {$\Sigma$};
	\node [sum, name=sn, left of = pn, node distance = 1.5cm] {$\Sigma$};

	\draw [connector] (s1) --  (p1) {};
	\draw [connector] (sn) --  (pn) {};
	\draw [connector] (p1) -- node (u1) {} (b1) {};
	\draw [connector] (pn) -- node (un) {} (bn) {};
	\draw [line] (pc) -- node  (omegahat) {} (connector1) {};
	\draw [connector] (connector1) |- (s1) {};
	\draw [connector] (connector1) |- (sn) {};
	\draw [line] (b1) -- (output1) {};
	\draw [line] (bn) -- (outputn) {};
	\draw [line] (output1) -- (node3) {};
	\draw [line] (outputn) -- (node4) {};
	\draw [connector] (node3) -| node (feedb1) {} (s1) {};
	\draw [connector] (node4) -| node (feedbn) {} (sn) {};
	\draw [line] (output1) -- (node6) {};
	\draw [connector] (outputn) -- (s2) {};
	\draw [connector] (node6) -- (s2) {};
	\draw [line] (s2) -- (node8) {};
	\draw [connector] (node8) -| (pc) {};
	\draw [connector] (node9) -| (pc) {};

	\node [draw=none,fill=none, above of =omegahat,node distance = 0.2cm, minimum width=1.0cm] {$\hat{\omega}$};
	\node [draw=none,fill=none, above of =u1,node distance = 0.2cm, minimum width=1.0cm] {$u_1$};
	\node [draw=none,fill=none, above of =un,node distance = 0.2cm, minimum width=1.0cm] {$u_n$};
	\node [draw=none,fill=none, above of =output1,node distance = 0.5cm, minimum width=1.0cm] {$\omega_1$};
	\node [draw=none,fill=none, above of =outputn,node distance = 0.5cm, minimum width=1.0cm] {$\omega_n$};
	\node [draw=none,fill=none, above left of =s1,node distance = 0.5cm, minimum width=1.0cm] {$+$};
	\node [draw=none,fill=none, below right of =s1,node distance = 0.5cm, minimum width=1.0cm] {$-$};
	\node [draw=none,fill=none, above left of =sn,node distance = 0.5cm, minimum width=1.0cm] {$+$};
	\node [draw=none,fill=none, below right of =sn,node distance = 0.5cm, minimum width=1.0cm] {$-$};
	\node [draw=none,fill=none, above left of =s2,node distance = 0.5cm, minimum width=1.0cm] {$+$};
	\node [draw=none,fill=none, below left of =s2,node distance = 0.5cm, minimum width=1.0cm] {$+$};
	\node [draw=none,fill=none, below left of =node9,node distance = 0.5cm, minimum width=1.0cm] {$\omega^{\text{ref}}$};

\end{tikzpicture}
&
\begin{tikzpicture}[auto, >=latex]
    	
   	\node [input, name=ref] {};
   	\node [block, above of=ref, node distance =1.0cm, minimum width=1.5cm] (p1) {$C'_1$};
	\node [draw=none,fill=none, below of =p1,node distance = 1.4cm, minimum width=1.0cm] (dots) {$\vdots$};
   	\node [block, below of =dots,node distance = 1.0cm, minimum width=1.5cm] (pn) {$C'_n$};
    	\node [output, name=l1c,left of =p1,node distance = 1.4cm]{};	
    	\node [output, name=l1,below of =l1c,node distance = 0.4cm]{};	
    	\node [output, name=lnc,left of =pn,node distance = 1.4cm]{};	
    	\node [output, name=ln,below of =lnc,node distance = 0.4cm]{};	
    	\node [output, name=l12,above of =l1,node distance = 0.67cm]{};	
    	\node [output, name=ln2,above of =ln,node distance = 0.67cm]{};

   	\node [block, right of=p1, node distance =2.2cm, minimum width=1.5cm] (b1) {Bus $1$};
   	\node [block, right of=pn, node distance =2.2cm, minimum width=1.5cm] (bn) {Bus $n$};
	\node [output, name=output1, right of=b1, node distance=1.3cm]{};
    	\node [output, name=outputn, right of=bn, node distance=1.3cm]{};
    	\node [output, name=node3, below of=output1, node distance=1.0cm]{};
    	\node [output, name=node4, below of=outputn, node distance=1.0cm]{};
    	\node [output, name=node5, below of=node3, node distance=1.0cm]{};	

	\draw [connector] (l1) |-  (p1.-160)  {};
	\draw [connector] (l12) |-  (p1.160)  {};
	\draw [connector] (ln) |-  (pn.-160)  {};
	\draw [connector] (ln2) |-  (pn.160)  {};
	\draw [connector] (p1) -- node (u1) {} (b1) {};
	\draw [connector] (pn) -- node (un) {} (bn) {};
	\draw [line] (b1) -- (output1) {};
	\draw [line] (bn) -- (outputn) {};
	\draw [line] (output1) -- (node3) {};
	\draw [line] (outputn) -- (node4) {};
	\draw [line] (node3) -| (l1) {};
	\draw [line] (node4) -| (ln) {};

	\node [draw=none,fill=none, above of =u1,node distance = 0.2cm, minimum width=1.0cm] {$u_1$};
	\node [draw=none,fill=none, above of =un,node distance = 0.2cm, minimum width=1.0cm] {$u_n$};
	\node [draw=none,fill=none, above right of =output1,node distance = 0.5cm, minimum width=1.0cm] {$\omega_1, \delta_1$};
	\node [draw=none,fill=none, above right of =outputn,node distance = 0.5cm, minimum width=1.0cm] {$\omega_n, \delta_n$};
	\node [draw=none,fill=none, left of =l12,node distance = 0.3cm, minimum width=1.0cm] {$\omega^{\text{ref}}$};
	\node [draw=none,fill=none, left of =ln2,node distance = 0.3cm, minimum width=1.0cm] {$\omega^{\text{ref}}$};

\end{tikzpicture}
\end{array}
$
\caption{The left figure illustrates the centralized controller, while the right figure illustrates the decentralized controller.}
\label{fig:power_decentralized}
\end{center}
\end{figure*}

We will here present a a centralized frequency control protocol for power systems and analyze its stability properties. 
Traditionally, frequency control of a power systems is carried out at two levels, see e.g., \cite{machowski2008power}. In the first level, the frequency is controlled with a proportional controller against a dynamic reference frequency. At the second level, the dynamic reference frequency is controlled with a proportional controller to eliminate static errors. We model the first level, proportional controller of an arbitrary bus $i$ as:
\begin{align}
u_i &= a (\hat{\omega} - \omega_i)
 \label{eq:u_central}
\end{align}
The second level proportional controller, regulating $\hat{\omega}$ is assumed to be given by:
\begin{align}
\dot{\hat{\omega}} &= b \Bigg(\omega^{\text{ref}} - \frac 1n \sum_{i\in\mathcal{V}} \omega_i \Bigg), 
 \label{eq:centralized_control_2}
\end{align}
where we have assumed that the average frequency of the buses is measured by the central controller\footnote{In reality the frequency is often measured at a specific bus. This will typically lead to longer delays, since disturbances need to propagate through the system before control action can be taken.}. Note that the second level controller integrates $\hat{\omega}$, thus acting as an integral controller. 
The centralized controller architecture is illustrated in Figure \ref{fig:power_decentralized}. 

\begin{proposition}
The power system described by \eqref{eq:swing_vector} where $u_i$ is given by  \eqref{eq:u_central}--\eqref{eq:centralized_control_2} satisfies $\lim_{t\rightarrow \infty} \omega_i(t) =  \omega^{\text{ref}}$ for any set of initial conditions, given that $a, b > 0$.
\end{proposition}

\begin{proof}
We may write \eqref{eq:swing_vector} with $u$ given by \eqref{eq:u_central}--\eqref{eq:centralized_control_2} as:
\begin{align}
\left[\begin{smallmatrix}
\dot{\hat{\omega}}\\ \dot{\delta} \\ \dot{\omega}
\end{smallmatrix}\right] &= 
{\left[\begin{smallmatrix}
0 &  0_{1\times n} & -\frac{b}{n} 1_{1\times n} \\
0_{n\times 1} & 0_{n\times n} & I_{n} \\
a M 1_{n\times 1} & -M \mathcal{L}_k & - M D -a M
\end{smallmatrix}\right]}
\left[\begin{smallmatrix}
\hat{\omega} \\ {\delta} \\ {\omega}
\end{smallmatrix}\right]  +
\left[\begin{smallmatrix}
 b \omega^{\text{ref}} \\ 0_{n\times 1} \\ Mp^m
\end{smallmatrix}\right].  
\label{eq:swing_vector_decentralized}
\end{align}
Define the output of the system as
\begin{align*}
\left[\begin{smallmatrix}
 y_\omega
\end{smallmatrix}\right] &= {\left[\begin{smallmatrix}
0_{n\times 1}  &  0_{n\times n} & I_n 
\end{smallmatrix}\right]}
\left[\begin{smallmatrix}
\hat{\omega} \\ {\delta} \\ {\omega}
\end{smallmatrix}\right],
\end{align*}
and consider the linear coordinate change:
\begin{align}
\label{eq:proof_pi_single_coordinate_change_2}
\delta = \left[\begin{smallmatrix}
 \frac{1}{\sqrt{n}} 1_{n\times 1} & S
\end{smallmatrix}\right] \delta' 
 \quad
\delta' = \left[\begin{smallmatrix}
 \frac{1}{\sqrt{n}} 1_{1\times n} \\ S^T
\end{smallmatrix}\right] \delta
\end{align}
where $S$ is a matrix such that $[
 \frac{1}{\sqrt{n}} 1_{n\times 1}, S 
]$ is an orthonormal matrix.
In the new coordinates, the system dynamics \eqref{eq:swing_vector_decentralized} become:
\begin{align*}
\left[\begin{smallmatrix}
\dot{\hat{\omega}}\\ \dot{\delta}' \\ \dot{\omega}
\end{smallmatrix}\right] &= 
{\left[\begin{smallmatrix}
0 &  0_{1\times n} & -\frac{b}{n} 1_{1\times n} \\
0_{n\times 1} & 0_{n\times n} & \left[\begin{smallmatrix}
 \frac{1}{\sqrt{n}} 1_{1\times n} \\ S^T
\end{smallmatrix}\right] \\
a M 1_{n\times 1} & \left[\begin{smallmatrix}
 0_{n\times 1} & -M \mathcal{L}_k S
\end{smallmatrix}\right]  & - M D -a M
\end{smallmatrix}\right]}
\left[\begin{smallmatrix}
\hat{\omega} \\ {\delta'} \\ {\omega}
\end{smallmatrix}\right]  +
\left[\begin{smallmatrix}
 b \omega^{\text{ref}} \\ 0_{n\times 1} \\ Mp^m
\end{smallmatrix}\right]. 
\end{align*}
We note that $\delta'_1$ is unobservable, and hence omit this state by defining $\delta'' = [\delta'_2, \dots, \delta'_n]^T$. In these coordinates the system dynamics \eqref{eq:swing_vector_decentralized} become:
\begin{align}
\left[\begin{smallmatrix}
\dot{\hat{\omega}}\\ \dot{\delta}'' \\ \dot{\omega}
\end{smallmatrix}\right] &= 
\underbrace{\left[\begin{smallmatrix}
0 &  0_{1\times (n-1)} & -\frac{b}{n} 1_{1\times n} \\
0_{(n-1)\times 1} & 0_{(n-1)\times (n-1)} &  S^T \\
a M 1_{n\times 1} & -M \mathcal{L}_k S  & - M D -a M
\end{smallmatrix}\right]}_{\triangleq A}
\left[\begin{smallmatrix}
\hat{\omega} \\ {\delta}'' \\ {\omega}
\end{smallmatrix}\right]  +
\left[\begin{smallmatrix}
 b \omega^{\text{ref}} \\ 0_{n\times 1} \\ Mp^m
\end{smallmatrix}\right].   \label{eq:swing_vector_decentralized_doubleprime}
\end{align}

By elementary row and column operations, we may rewrite the characteristic equation of $A$ as
$
 \det\left( {M}\mathcal{L}_kSS^T + \frac{a b}{n}M1_{n\times n} + s(MD + a M) + s^2 I_n  \right) = 0.
$
This equation is equivalent to
\begin{align}
\label{eq:proof_suboptimal_centralized_1}
\begin{aligned}
 \det\left( \mathcal{L}_kSS^T + \frac{a b}{n}1_{n\times n} + s(D + a I_n) + s^2 M^{-1}  \right) = 0.
\end{aligned}
\end{align}
For a given $s$, this equation has one solution only if \\
$
x^T(\mathcal{L}_kSS^T + \frac{a b}{n}1_{n\times n} + s(D + a I_n) + s^2 M^{-1}   )x=0
$
has a solution for some $\norm{x}=1$. Hence, if the previous equation has all its solutions in $\mathbb{C}^-$ for all $\norm{x}=1$, then \eqref{eq:proof_suboptimal_centralized_1} has all its solutions in $\mathbb{C}^-$. By similar arguments used in the proof of Theorem \ref{th:single_abs}, we can show that $\mathcal{L}_kSS^T=\mathcal{L}_k$. 
Thus, if the equation 
\begin{align*}
\underbrace{x^T( \mathcal{L}_k +  \frac{a b}{n} 1_{n\times n} )x}_{a_0} + s \underbrace{x^T (D+a I_n) x}_{a_1} + s^2 \underbrace{x^T M^{-1} x}_{a_2} =0 
\end{align*}
has all its solutions in $\mathbb{C}^-$, then $A$ has only one zero eigenvalue, and all other eigenvalues in $\mathbb{C}^-$. By the Routh-Hurwitz stability criterion, the above equation has all its solutions in $\mathbb{C}^-$ iff $a_i>0, i=0,1,2$. Clearly $a_0>0$ since $x^T\mathcal{L}x > 0$ for $x\ne c 1_{n\times 1}$, where $c\in \mathbb{R}$, and $1_{n\times n}1_{n\times 1} \ne 0_{n\times 1}$. Also $a_1, a_2>0$ since $D+a I_n$ and $M^{-1}$ are diagonal with positive elements. 
We conclude thus that $A$ is Hurwitz. 
Stationarity of \eqref{eq:swing_vector_decentralized_doubleprime} implies  $\omega = \omega^{\text{ref}} 1_{n\times 1} $.
\end{proof}

\subsubsection{Decentralized PI control}
In this section we analyze a decentralized frequency controller, where each bus controls its own frequency based only on local phase and frequency measurements. 
Thus, there is no need to send control signals or reference values to the buses. 
This controller architecture might be favorable due to security concerns when sending unencrypted data over large areas.  
Another benefit is improved performance when the tripping of one or several power lines causes the network to be split up into two or more sub-networks, so called islanding. 
 The controller of node $i$ is assumed to be given by \eqref{eq:decentralized_control_22}, here written as
\begin{align}
\dot{z}_i &=   \omega_i - \omega^{\text{ref}}
 \label{eq:decentralized_control_1} \\
u_i &= a ( \omega^{\text{ref}} - \omega_i) - b z_i.
 \label{eq:decentralized_control_2}
\end{align}
The controller architecture is illustrated in Figure~\ref{fig:power_decentralized}.
The decentralized controller \eqref{eq:decentralized_control_1}--\eqref{eq:decentralized_control_2} is typically not practically feasible with only frequency measurements available at the generation buses. Even the slightest measurement error will be integrated and cause instability, see, e.g., \cite{machowski2008power}. 
 However, with recent advances in phasor measurement unit (PMU) technology however, phase measurements are becoming increasingly available \cite{Phadke1993}. 
 By integrating \eqref{eq:decentralized_control_1} we obtain
$z_i =  \omega^{\text{ref}} t - \delta_i.
$ This implies that in order to accurately estimate the integral state $z_i$, each generator bus needs access only to time and phase measurements, both provided by PMU's with high accuracy. 
\begin{proposition}
\label{th:power_decentralized_stability}
The power system described by \eqref{eq:swing_vector} where $u_i$ is given by  \eqref{eq:decentralized_control_1}--\eqref{eq:decentralized_control_2} satisfies $\lim_{t\rightarrow \infty} \omega_i(t) =  \omega^{\text{ref}}$ for any set of initial conditions, given that $a, b > 0$.
\end{proposition}
\begin{proof}
If we consider $[\mathcal{B}^T \delta, \omega]$ to be the output, the dynamics of \eqref{eq:swing_vector} may be modified as long as the dynamics of $[\mathcal{B}^T \delta, \omega]$ are left unchanged. We thus may rewrite \eqref{eq:swing_vector} with $u$ given by \eqref{eq:decentralized_control_1}--\eqref{eq:decentralized_control_2} as:
\begin{eqnarray*}
\left[\begin{smallmatrix}
\dot{z} \\ \dot{\omega}
\end{smallmatrix}\right] = 
\underbrace{\left[\begin{smallmatrix}
0_{n\times n} & I_{n} \\
-M \mathcal{L}_k{-}b M & - M D{-}a M
\end{smallmatrix}\right]}_{\triangleq A} 
\left[\begin{smallmatrix}
{z} \\ {\omega}
\end{smallmatrix}\right] +
\left[\begin{smallmatrix}
- \omega^{\text{ref}}1_{n\times 1} \\ M (p^m{+}a \omega^{\text{ref}} 1_{n\times 1} )
\end{smallmatrix}\right],
\end{eqnarray*}
since $\dot{\delta}-\dot{z} =  \omega^{\text{ref}}1_{n\times 1}$, implying that ${\delta}-{z} = \delta_0- t \omega^{\text{ref}}1_{n\times 1}$. Since $\mathcal{L}_k 1_{n\times 1}=0_{n\times 1}$, the output dynamics of the above equation is equivalent to that of \eqref{eq:swing_vector} with respect to the output  $[\mathcal{B}^T \delta, \omega]$. 
By elementary column operations, we may rewrite the characteristic equation of $A$ as
\begin{align}
\label{eq:proof_suboptimal_decentralized_1}
\begin{aligned}
&\det\left(  s^2 I_n +  s  MD + M\mathcal{L} + b I _{n} \right) = 0 \Leftrightarrow \\ 
&\det\left(  s^2 M^{-1} +  s D + \mathcal{L} +  b M^{-1} \right) = 0. 
\end{aligned}
\end{align}
For a given $s$, the above equation has a solution only if 
$
x^T\left(  s^2 M^{-1} +  s  D + \mathcal{L} +  b M^{-1} \right)x=0
$
has a solution. Hence, if the previous equation has all its solutions in $\mathbb{C}^-$ for all $\norm{x}=1$, then \eqref{eq:proof_suboptimal_decentralized_1} has all its solutions in $\mathbb{C}^-$. Thus, if the equation 
\begin{align*}
\underbrace{x^T(\mathcal{L} +  b M^{-1})x}_{a_0} + s \underbrace{x^T D x}_{a_1} + s^2 \underbrace{x^TM^{-1} x}_{a_2} =0 
\end{align*}
has all its solutions in $\mathbb{C}^-$, then $A$ is stable. By the Routh-Hurwitz stability criterion, the above equation has all its solutions in $\mathbb{C}^-$ iff $a_i>0, i=0,1,2$. Clearly $a_0>0$ since $x^T\mathcal{L}x > 0$ for $x\ne c 1_{n\times 1}$ for any $c\in \mathbb{R}$, and $M^{-1}1_{n\times 1} \ne 0_{n\times 1}$. Also $a_1, a_2>0$ since $D$ and $M^{-1}$ are diagonal with positive elements. We conclude thus that $A$ is Hurwitz. 
 Now consider the coordinate shift
\begin{align*}
\left[\begin{smallmatrix}
{z'} \\ {\omega'}
\end{smallmatrix}\right] = 
\left[\begin{smallmatrix}
{z} \\ {\omega} 
\end{smallmatrix}\right] -
\left[\begin{smallmatrix}
{z_o} \\ {\omega_o}
\end{smallmatrix}\right]
\end{align*}
where 
\begin{align*}
z_0=(b I_{n} +  \mathcal{L}_k )^{-1}(D  \omega^{\text{ref}} 1_{n\times 1} - p^m),\qquad
\omega_0= \omega^{\text{ref}} 1_{n\times 1}.
\end{align*}
In the translated coordinates, the origin is the only equilibrium of the system. Hence
$ \lim_{t\rightarrow \infty} \omega_i(t) =  \omega^{\text{ref}} \; \forall i \;\in \mathcal{V} 
$. \end{proof}


	\setlength\fheight{3.0cm} 
	\setlength\fwidth{6.0cm} 
\begin{figure*}[t]
	\centering
$
\begin{array}{cc}
\quad \quad\text{Centralized} & \quad \quad\text{Decentralized} \\

	\input{Simulations/Powersystems/omega_cen.tikz} & \input{Simulations/Powersystems/omega_dec.tikz} 	 \\ 	\input{Simulations/Powersystems/u_cen.tikz} & \input{Simulations/Powersystems/u_dec.tikz} 
\end{array}
$
\caption{The upper left figure shows the bus frequencies with centralized frequency control, while the lower left figure shows the control signals at all buses. 
The upper right figure shows the bus frequencies with decentralized frequency control, while the lower right figure shows the control signals at all buses. The controller parameters were  $a=0.8$, $b=0.04$, for both the centralized and the decentralized controller. } 
\label{fig:powersystems_sim_1_disturbance}
\end{figure*}

\begin{figure*}[t]
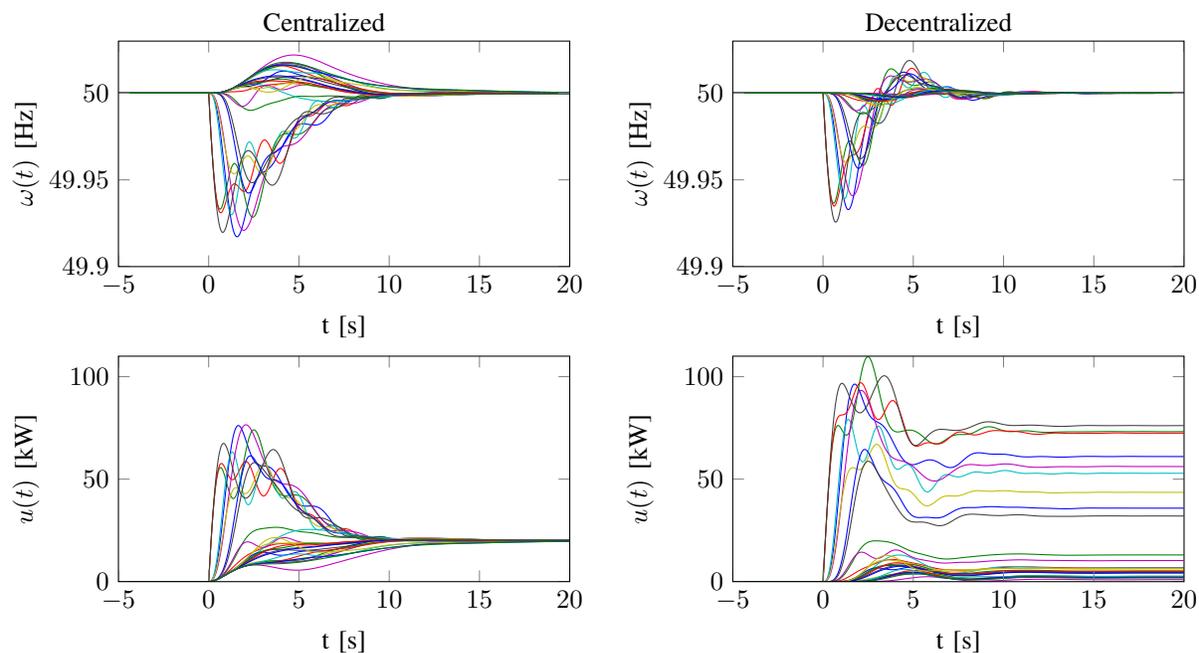

	\centering
$
\begin{array}{cc}
\quad \quad\text{Centralized} & \quad \quad\text{Decentralized} \\
	\input{Simulations/Powersystems/omega_cen_I.tikz} & \input{Simulations/Powersystems/omega_dec_I.tikz} 	 \\ 	\input{Simulations/Powersystems/u_cen_I.tikz} & \input{Simulations/Powersystems/u_dec_I.tikz} 
\end{array}
$
\caption{The upper left figure shows the bus frequencies with centralized frequency control, while the lower left figure shows the control signals at all buses. 
The upper right figure shows the bus frequencies with decentralized frequency control, while the lower right figure shows the control signals at all buses. The controller parameters were  $a=0.8$, $b=0.8$, for both the centralized and the decentralized controller. } 
\label{fig:powersystems_sim_2_disturbance} 
\end{figure*}

\subsubsection{Simulations}
The centralized and decentralized frequency controllers were tested on the IEEE 30 bus test system \cite{IEEE30}. 
The line admittances were extracted from  \cite{IEEE30} and the voltages were assumed to be 132 kV for all buses.
The values of $M$ and $D$ were assumed to be given by $m_i = 10^5\; \text{kg}\,\text{m}^2$ and $d_i = 1 \; s^{-1} $, respectively, $ \forall i \in \mathcal{V}$.  
The power system was assumed to be initially in an operational equilibrium, until the power load is increased by a step of $200$ kW in the buses $2,3$ and $7$. This will immediately result in decreased frequencies at these buses. The frequency controllers at the buses will then control the frequencies towards the reference frequency of $\omega^{\text{ref}}=50$ Hz. 
The controller parameters were set to  $a=0.8$, $b=0.04$, for both the centralized and the decentralized controller. 
The step responses of the frequencies are plotted in figure~\ref{fig:powersystems_sim_1_disturbance}. We note that for the centralized PI controller, the generation is increased uniformly among the generators. If however the integral action is distributed amongst the generators, some generators will increase their generation more than others.
Figure \ref{fig:powersystems_sim_2_disturbance} shows the step response under the significantly larger integral action, $\beta = 0.8$ for both the centralized and the decentralized controller. We notice that the step response of the decentralized controller shows better performance compared to the centralized controller. This is due to that the centralized controller can only measure the average frequency in the power system, as opposed to the individual frequencies which the decentralized controller can measure.  

Note that the controller parameters are assumed to be identical for both the centralized and the decentralized controller. This might not be restrictive when the generators are homogeneous. However, when generators are more heterogeneous, this might be restrictive.

\section{Conclusions}
\label{sec:conclusions}
In this paper we have studied a class of nonlinear consensus protocols for single and double-integrator dynamics. Necessary and sufficient conditions for consensus were derived for static communication topologies under single and double-integrator dynamics. In all cases, expressions for the convergence points were given.
We have also studied consensus controllers with integral action for agents with single integrator dynamics and agents with damped double-integrator dynamics. We proved that with the proposed consensus controllers, the agents reach asymptotic consensus even in the presence of constant disturbances. If we allow for absolute position measurements the agents in addition converge asymptotically to a common state. In the absence of disturbances, the proposed consensus protocols asymptotically solve the initial average consensus problem. 
We have demonstrated by simulations that the proposed controllers have applications in controlling autonomous satellites in space, control of mobile robots, building temperature control and frequency control of electrical power systems.

\bibliography{references}
\bibliographystyle{plain}

\appendix

\begin{proof}[Proof of Lemma \ref{lem:1_nl}]
 Since $\Omega$ is compact, the relative states $\bar{x}$ are bounded. 
 Then clearly $x$ is bounded if and only if the average $x'=\frac 1n \sum_{i\in \mathcal{V}} x_i $ is also bounded, as seen by the following inequalities
 \begin{align*}
 \norm{x}_\infty < n\norm{\bar{x}}_\infty + |x'|  \\
 |x'| = \left|\frac 1n \sum_{i\in \mathcal{V}} x_i\right| < \frac 1n \norm{x}_\infty.
 \end{align*}
Let
$
E_0 = E(x_0,v_0) =  \sum_{i\in \mathcal{V}} \left( \int_0^{x_i} \kappa_i(y) \ud y + {v_i} \right),
$
where  $[x_0,v_0]$ denotes the initial condition. 
Since $[\bar{x}(t), v(t)] $ evolve in the compact set $\Omega$, $v_i(t)$ is also bounded. Hence $\forall i \in \mathcal{V}$ $\exists M\in \mathbb{R}^+: |v_i(t)|\le M \; \forall t\ge0, \; \forall i\in \mathcal{V}$. By Assumption~\ref{ass:gamma} $\kappa_i(x)\ge\underline{\kappa}>0 \quad \forall i \in\mathcal{V},\; \forall x \in \mathbb{R}$. Using these inequalities we obtain
\begin{align}
\Bigg| \sum_{i \in \mathcal{V}} \int_0^{x_i} \kappa_i(y) \ud y \Bigg| \le nM + |E_0|. 
\label{eq_bounded_gamma_function}
\end{align}
Assume for the sake of contradiction that $x'(t)$ is unbounded. Let us consider the case when $x'(t)\rightarrow + \infty$. Since $\bar{x}$ is bounded by say $M'>0$ in the $\infty$-norm, and $\mathcal{G}$ is connected, $|x_i(t){-}x_j(t)|$ is bounded by $(n-1)M'$ $\forall i,j \in \mathcal{V}$. Thus $x_i(t)>0 \; \forall i \in \mathcal{V}$ whenever $x'(t)> (n-1)M'$. Thus, if $x'(t)> (n-1)M'$, we obtain that
$
\sum_{i\in\mathcal{V}} \int_0^{x_i} \kappa_i(y) \ud y \ge \sum_{i\in\mathcal{V}} \underline{\kappa} x_i
$. 
By assumption, $x'(t)$ is unbounded, implying that also $\sum_{i \in \mathcal{V}} x_i(t)$ is unbounded. Thus $\exists t_1:  \sum_{i\in\mathcal{V}}  x_i(t_1) > \max\{ {1}/{\underline{\kappa}} \left( nM + |E_0| \right),   M'\}$. But this contradicts \eqref{eq_bounded_gamma_function}. Hence $x'(t)$ must be bounded. The cases when $x'(t) \rightarrow -\infty$ as well as the case when no limit of $x'(t)$ exists are treated analogously. 
We conclude that $x$ must be bounded, and thus $\Omega'$ is compact by the Heine-Borel Theorem. 
\end{proof}

\section*{Acknowledgment}
\addcontentsline{toc}{section}{Acknowledgment}
The authors would like to thank Guodong Shi and Tao Yang, as well as the anonymous reviewers for their valuable feedback.

\begin{IEEEbiography}[{\includegraphics[width=1in,height=1.25in,clip,keepaspectratio]{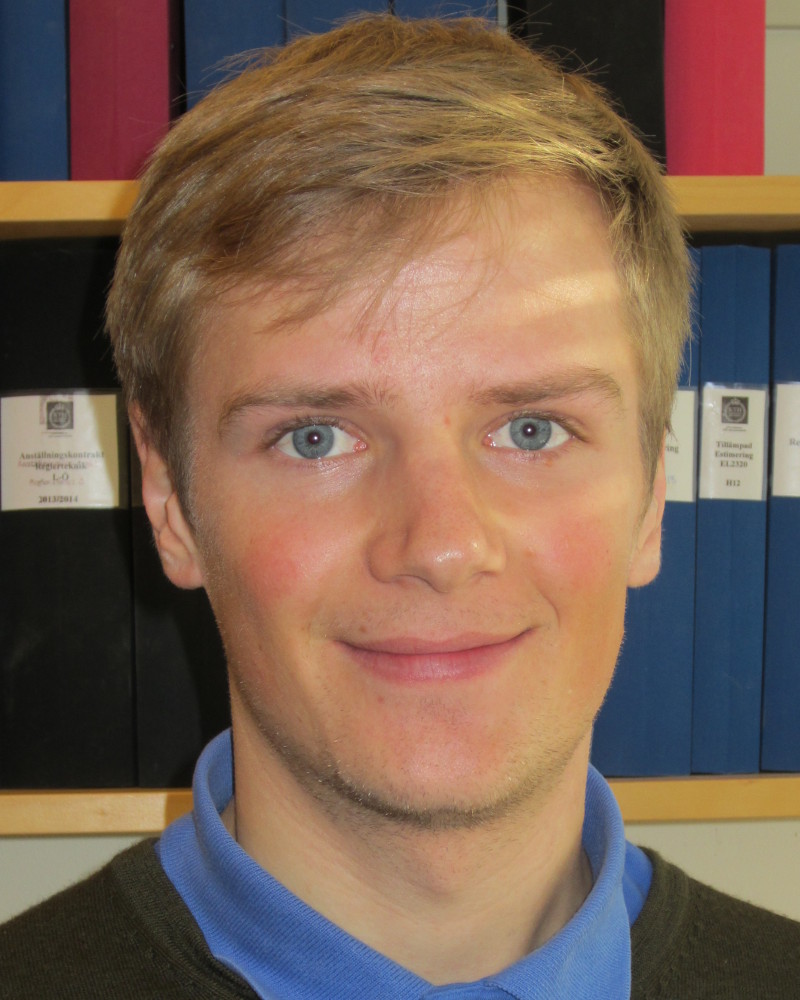}}]{Martin Andreasson} 
Martin Andreasson received the M.Sc. degree in engineering physics KTH Royal Institute of Technology, Stockholm, Sweden, in 2011. He is currently a PhD Student at the Automatic Control Laboratory, KTH Royal Institute of Technology, Stockholm, Sweden. His research interests include distributed control of multi-agent systems, and control of power systems. 
\end{IEEEbiography}

\begin{IEEEbiography}[{\includegraphics[width=1in,height=1.25in,clip,keepaspectratio]{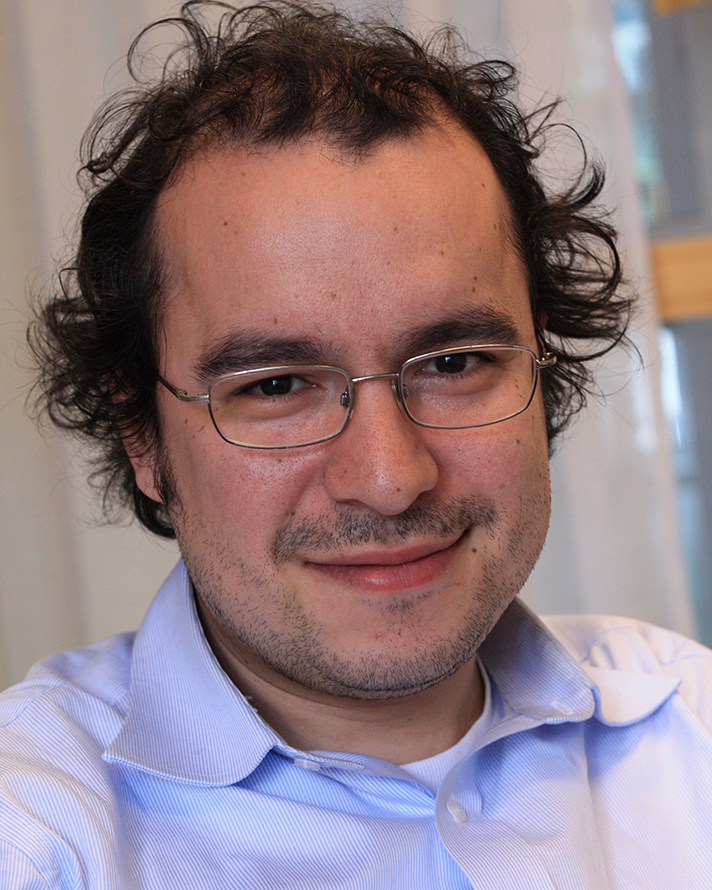}}]{Dimos V. Dimarogonas} 
Dimos V. Dimarogonas was born in Athens, Greece,
in 1978. He received the Diploma in Electrical and
Computer Engineering in 2001 and the Ph.D. in Mechanical
Engineering in 2006, both from the National Technical
University of Athens (NTUA), Greece. Between May 2007
and February 2009, he was a Postdoctoral Researcher
at the Automatic Control Laboratory, School of Electrical
Engineering, ACCESS Linnaeus Center, Royal Institute of
Technology (KTH), Stockholm, Sweden. Between February
2009 and March 2010, he was a Postdoctoral Associate
at the Laboratory for Information and Decision Systems
(LIDS) at the Massachusetts Institute of Technology (MIT), Boston, MA, USA. He
is currently an Assistant Professor at the Automatic Control Laboratory, ACCESS
Linnaeus Center, Royal Institute of Technology (KTH), Stockholm, Sweden.
His current research interests include Multi-Agent Systems, Hybrid Systems
and Control, Robot Navigation and Networked Control. He was awarded a Docent in
Automatic Control from KTH in 2012. He serves in the Editorial Board of Automatica
and the IET Control Theory and Applications and is a member of IEEE and the
Technical Chamber of Greece.
\end{IEEEbiography}

\begin{IEEEbiography}[{\includegraphics[width=1in,height=1.25in,clip,keepaspectratio]{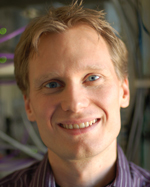}}]{Henrik Sandberg} 
Henrik Sandberg received the M.Sc. degree in engineering physics and the Ph.D. degree in automatic control from Lund University, Lund, Sweden, in 1999 and 2004, respectively. He is an Associate Professor with the Automatic Control Laboratory, KTH Royal Institute of Technology, Stockholm, Sweden. From 2005 to 2007, he was a Post-Doctoral Scholar with the California Institute of Technology, Pasadena, USA. In 2013, he was a visiting scholar at the Laboratory for Information and Decision Systems (LIDS) at MIT, Cambridge, USA. He has also held visiting appointments with the Australian National University and the University of Melbourne, Australia. His current research interests include secure networked control, power systems, model reduction, and fundamental limitations in control. Dr. Sandberg was a recipient of the Best Student Paper Award from the IEEE Conference on Decision and Control in 2004 and an Ingvar Carlsson Award from the Swedish Foundation for Strategic Research in 2007. He is currently an Associate Editor of the IFAC Journal Automatica.
\end{IEEEbiography}

\begin{IEEEbiography}[{\includegraphics[width=1in,height=1.25in,clip,keepaspectratio]{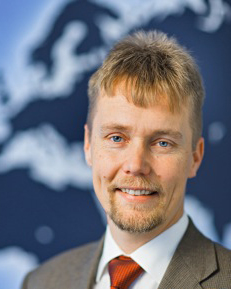}}]{Karl Henrik Johansson} 
Karl Henrik Johansson is Director of the KTH ACCESS Linnaeus Centre and Professor at the School of Electrical Engineering, Royal Institute of Technology, Sweden. He is a Wallenberg Scholar and has held a six-year Senior Researcher Position with the Swedish Research Council. He is Director of the Stockholm Strategic Research Area ICT The Next Generation. He received MSc and PhD degrees in Electrical Engineering from Lund University. He has held visiting positions at UC Berkeley (1998-2000) and California Institute of Technology (2006-2007). His research interests are in networked control systems, hybrid and embedded system, and applications in transportation, energy, and automation systems. He has been a member of the IEEE Control Systems Society Board of Governors and the Chair of the IFAC Technical Committee on Networked Systems. He has been on the Editorial Boards of several journals, including Automatica, IEEE Transactions on Automatic Control, and IET Control Theory and Applications. He is currently on the Editorial Board of IEEE Transactions on Control of Network Systems and the European Journal of Control. He has been Guest Editor for special issues, including the one on “Wireless Sensor and Actuator Networks” of IEEE Transactions on Automatic Control 2011. He was the General Chair of the ACM/IEEE Cyber-Physical Systems Week 2010 in Stockholm and IPC Chair of many conferences. He has served on the Executive Committees of several European research projects in the area of networked embedded systems. In 2009, he received the Best Paper Award of the IEEE International Conference on Mobile Ad-hoc and Sensor Systems. In 2009, he was also awarded Wallenberg Scholar, as one of the first ten scholars from all sciences by the Knut and Alice Wallenberg Foundation. He was awarded an Individual Grant for the Advancement of Research Leaders from the Swedish Foundation for Strategic Research in 2005. He received the triennial Young Author Prize from IFAC in 1996 and the Peccei Award from the International Institute of System Analysis, Austria, in 1993. He received Young Researcher Awards from Scania in 1996 and from Ericsson in 1998 and 1999. He is a Fellow of the IEEE.
\end{IEEEbiography}

\end{document}

%% file: Simulations/hc.tikz
%
%
\begin{tikzpicture}

\begin{axis}[%
scale only axis,
width=\fwidth,
height=\fheight,
xmin=22, xmax=25,
ymin=0, ymax=1500,
xlabel={$\text{T [}^{\circ}\text{C]}$},
ylabel={$1/\gamma{}_i(T)\text{ [KJ}/^{\circ}\text{C]}$},
axis on top]
\addplot [
color=blue,
solid
]
coordinates{
 (21.992,500)(22.002,500)(22.012,500)(22.022,500)(22.032,500)(22.042,500)(22.0521,500)(22.0621,500)(22.0721,500)(22.0821,500)(22.0921,500)(22.1021,500)(22.1121,500)(22.1221,500)(22.1321,500)(22.1421,500)(22.1522,500)(22.1622,500)(22.1722,500)(22.1822,500)(22.1922,500)(22.2022,500)(22.2122,500)(22.2222,500)(22.2322,500)(22.2422,500)(22.2523,500)(22.2623,500)(22.2723,500)(22.2823,500)(22.2923,500)(22.3023,500)(22.3123,500)(22.3223,500)(22.3323,500)(22.3423,500)(22.3524,500)(22.3624,500)(22.3724,500)(22.3824,500)(22.3924,500)(22.4024,500)(22.4124,500)(22.4224,500)(22.4324,500)(22.4424,500)(22.4525,500)(22.4625,500)(22.4725,500)(22.4825,500)(22.4925,500)(22.5025,500)(22.5125,500)(22.5225,500)(22.5325,500)(22.5425,500)(22.5526,500)(22.5626,500)(22.5726,500)(22.5826,500)(22.5926,500)(22.6026,500)(22.6126,500)(22.6226,500)(22.6326,500)(22.6426,500)(22.6527,500)(22.6627,500)(22.6727,500)(22.6827,500)(22.6927,500)(22.7027,500)(22.7127,500)(22.7227,500)(22.7327,500)(22.7427,500)(22.7528,500)(22.7628,500)(22.7728,500)(22.7828,500)(22.7928,500)(22.8028,500)(22.8128,500)(22.8228,500)(22.8328,500)(22.8428,500)(22.8529,500)(22.8629,500)(22.8729,500)(22.8829,500)(22.8929,500)(22.9029,500)(22.9129,500)(22.9229,500)(22.9329,500)(22.9429,500)(22.953,500)(22.963,500)(22.973,500)(22.983,500)(22.993,500)(23.003,506.006)(23.013,526.026)(23.023,546.046)(23.033,566.066)(23.043,586.086)(23.0531,606.106)(23.0631,626.126)(23.0731,646.146)(23.0831,666.166)(23.0931,686.186)(23.1031,706.206)(23.1131,726.226)(23.1231,746.246)(23.1331,766.266)(23.1431,786.286)(23.1532,806.306)(23.1632,826.326)(23.1732,846.346)(23.1832,866.366)(23.1932,886.386)(23.2032,906.406)(23.2132,926.426)(23.2232,946.446)(23.2332,966.466)(23.2432,986.486)(23.2533,1006.51)(23.2633,1026.53)(23.2733,1046.55)(23.2833,1066.57)(23.2933,1086.59)(23.3033,1106.61)(23.3133,1126.63)(23.3233,1146.65)(23.3333,1166.67)(23.3433,1186.69)(23.3534,1206.71)(23.3634,1226.73)(23.3734,1246.75)(23.3834,1266.77)(23.3934,1286.79)(23.4034,1306.81)(23.4134,1326.83)(23.4234,1346.85)(23.4334,1366.87)(23.4434,1386.89)(23.4535,1406.91)(23.4635,1426.93)(23.4735,1446.95)(23.4835,1466.97)(23.4935,1486.99)(23.5035,1492.99)(23.5135,1472.97)(23.5235,1452.95)(23.5335,1432.93)(23.5435,1412.91)(23.5536,1392.89)(23.5636,1372.87)(23.5736,1352.85)(23.5836,1332.83)(23.5936,1312.81)(23.6036,1292.79)(23.6136,1272.77)(23.6236,1252.75)(23.6336,1232.73)(23.6436,1212.71)(23.6537,1192.69)(23.6637,1172.67)(23.6737,1152.65)(23.6837,1132.63)(23.6937,1112.61)(23.7037,1092.59)(23.7137,1072.57)(23.7237,1052.55)(23.7337,1032.53)(23.7437,1012.51)(23.7538,992.492)(23.7638,972.472)(23.7738,952.452)(23.7838,932.432)(23.7938,912.412)(23.8038,892.392)(23.8138,872.372)(23.8238,852.352)(23.8338,832.332)(23.8438,812.312)(23.8539,792.292)(23.8639,772.272)(23.8739,752.252)(23.8839,732.232)(23.8939,712.212)(23.9039,692.192)(23.9139,672.172)(23.9239,652.152)(23.9339,632.132)(23.9439,612.112)(23.954,592.092)(23.964,572.072)(23.974,552.052)(23.984,532.032)(23.994,512.012)(24.004,500)(24.014,500)(24.024,500)(24.034,500)(24.044,500)(24.0541,500)(24.0641,500)(24.0741,500)(24.0841,500)(24.0941,500)(24.1041,500)(24.1141,500)(24.1241,500)(24.1341,500)(24.1441,500)(24.1542,500)(24.1642,500)(24.1742,500)(24.1842,500)(24.1942,500)(24.2042,500)(24.2142,500)(24.2242,500)(24.2342,500)(24.2442,500)(24.2543,500)(24.2643,500)(24.2743,500)(24.2843,500)(24.2943,500)(24.3043,500)(24.3143,500)(24.3243,500)(24.3343,500)(24.3443,500)(24.3544,500)(24.3644,500)(24.3744,500)(24.3844,500)(24.3944,500)(24.4044,500)(24.4144,500)(24.4244,500)(24.4344,500)(24.4444,500)(24.4545,500)(24.4645,500)(24.4745,500)(24.4845,500)(24.4945,500)(24.5045,500)(24.5145,500)(24.5245,500)(24.5345,500)(24.5445,500)(24.5546,500)(24.5646,500)(24.5746,500)(24.5846,500)(24.5946,500)(24.6046,500)(24.6146,500)(24.6246,500)(24.6346,500)(24.6446,500)(24.6547,500)(24.6647,500)(24.6747,500)(24.6847,500)(24.6947,500)(24.7047,500)(24.7147,500)(24.7247,500)(24.7347,500)(24.7447,500)(24.7548,500)(24.7648,500)(24.7748,500)(24.7848,500)(24.7948,500)(24.8048,500)(24.8148,500)(24.8248,500)(24.8348,500)(24.8448,500)(24.8549,500)(24.8649,500)(24.8749,500)(24.8849,500)(24.8949,500)(24.9049,500)(24.9149,500)(24.9249,500)(24.9349,500)(24.9449,500)(24.955,500)(24.965,500)(24.975,500)(24.985,500)(24.995,500)(25.005,500) 
};

\end{axis}
\end{tikzpicture}

%% file: Simulations/consensus_without_I.tikz
%
%
\begin{tikzpicture}

\definecolor{mycolor1}{rgb}{0,0.5,0}
\definecolor{mycolor2}{rgb}{0,0.75,0.75}
\definecolor{mycolor3}{rgb}{0.75,0,0.75}

\begin{axis}[%
scale only axis,
width=\fwidth,
height=\fheight,
xmin=0, xmax=20,
ymin=-6, ymax=8,
xlabel={$t$},
ylabel={$x(t)$},
axis on top]
\addplot [
color=blue,
solid
]
coordinates{
 (0,5)(0.000119652,5)(0.000777737,4.99998)(0.00197426,4.99989)(0.00855511,4.99805)(0.0205203,4.98888)(0.0557841,4.92087)(0.114347,4.68969)(0.172909,4.34277)(0.283485,3.50036)(0.39406,2.59886)(0.504636,1.81027)(0.626408,1.19747)(0.74818,0.900921)(0.869953,0.89707)(0.991725,1.11463)(1.12526,1.50007)(1.25879,1.92976)(1.39233,2.31947)(1.52586,2.62247)(1.66253,2.83042)(1.7992,2.94903)(1.93587,3.00619)(2.0896,3.0342)(2.26038,3.05873)(2.39868,3.0899)(2.53699,3.13484)(2.67529,3.19042)(2.8149,3.25142)(2.95583,3.31199)(3.09676,3.36717)(3.23769,3.41488)(3.38919,3.45773)(3.54069,3.49331)(3.66231,3.51797)(3.78392,3.54023)(3.90553,3.56084)(4.03526,3.58151)(4.16499,3.60114)(4.29471,3.61982)(4.46981,3.64358)(4.61562,3.66212)(4.76143,3.6796)(4.90724,3.69616)(5.0797,3.71472)(5.2788,3.73503)(5.4779,3.75429)(5.677,3.77259)(5.96088,3.79713)(6.24475,3.82016)(6.52863,3.84214)(6.85863,3.86689)(7.18863,3.89104)(7.51863,3.91465)(7.97361,3.94648)(8.42859,3.97778)(8.88357,4.00877)(9.51072,4.05113)(10.1379,4.09324)(11.2085,4.16487)(12.7226,4.26592)(14.2367,4.36689)(17.4686,4.58236)(20,4.75112) 
};

\addplot [
color=mycolor1,
solid
]
coordinates{
 (0,-6)(0.000119652,-6)(0.000777737,-5.99996)(0.00197426,-5.99976)(0.00855511,-5.99548)(0.0205203,-5.97427)(0.0557841,-5.81678)(0.114347,-5.28105)(0.172909,-4.47555)(0.283485,-2.50926)(0.39406,-0.378772)(0.504636,1.53113)(0.626408,3.10156)(0.74818,4.00771)(0.869953,4.30257)(0.991725,4.14245)(1.12526,3.67584)(1.25879,3.13511)(1.39233,2.69064)(1.52586,2.42862)(1.66253,2.36069)(1.7992,2.45239)(1.93587,2.63595)(2.0896,2.87007)(2.26038,3.08772)(2.39868,3.20493)(2.53699,3.27016)(2.67529,3.29722)(2.8149,3.30361)(2.95583,3.30447)(3.09676,3.30992)(3.23769,3.32399)(3.38919,3.34789)(3.54069,3.37689)(3.66231,3.40078)(3.78392,3.42329)(3.90553,3.44349)(4.03526,3.46226)(4.16499,3.47858)(4.29471,3.49319)(4.46981,3.51146)(4.61562,3.5262)(4.76143,3.54078)(4.90724,3.55512)(5.0797,3.57153)(5.2788,3.58936)(5.4779,3.606)(5.677,3.62174)(5.96088,3.64347)(6.24475,3.66498)(6.52863,3.68623)(6.85863,3.71022)(7.18863,3.73345)(7.51863,3.75629)(7.97361,3.78754)(8.42859,3.81853)(8.88357,3.84925)(9.51072,3.89139)(10.1379,3.93341)(11.2085,4.00492)(12.7226,4.10594)(14.2367,4.20689)(17.4686,4.42236)(20,4.59112) 
};

\addplot [
color=red,
solid
]
coordinates{
 (0,8)(0.000119652,8)(0.000777737,7.99997)(0.00197426,7.99982)(0.00855511,7.99674)(0.0205203,7.98148)(0.0557841,7.86818)(0.114347,7.48381)(0.172909,6.90926)(0.283485,5.52828)(0.39406,4.08163)(0.504636,2.86009)(0.626408,1.97132)(0.74818,1.61207)(0.869953,1.70508)(0.991725,2.09681)(1.12526,2.66495)(1.25879,3.19243)(1.39233,3.55999)(1.52586,3.72823)(1.66253,3.71993)(1.7992,3.59447)(1.93587,3.42678)(2.0896,3.2627)(2.26038,3.16441)(2.39868,3.15474)(2.53699,3.19099)(2.67529,3.24978)(2.8149,3.3101)(2.95583,3.35712)(3.09676,3.38432)(3.23769,3.39352)(3.38919,3.39061)(3.54069,3.38347)(3.66231,3.3797)(3.78392,3.37997)(3.90553,3.38481)(4.03526,3.39427)(4.16499,3.40658)(4.29471,3.41993)(4.46981,3.43706)(4.61562,3.4493)(4.76143,3.45956)(4.90724,3.46836)(5.0797,3.47795)(5.2788,3.48931)(5.4779,3.50173)(5.677,3.51517)(5.96088,3.53501)(6.24475,3.55444)(6.52863,3.57326)(6.85863,3.59498)(7.18863,3.61694)(7.51863,3.63903)(7.97361,3.66938)(8.42859,3.69968)(8.88357,3.73003)(9.51072,3.77184)(10.1379,3.81364)(11.2085,3.88503)(12.7226,3.98595)(14.2367,4.0869)(17.4686,4.30236)(20,4.47112) 
};

\addplot [
color=mycolor2,
solid
]
coordinates{
 (0,4)(0.000119652,4)(0.000777737,4.00001)(0.00197426,4.00005)(0.00855511,4.0009)(0.0205203,4.00514)(0.0557841,4.0365)(0.114347,4.14171)(0.172909,4.29509)(0.283485,4.63882)(0.39406,4.94085)(0.504636,5.10556)(0.626408,5.07981)(0.74818,4.85365)(0.869953,4.49532)(0.991725,4.09536)(1.12526,3.70853)(1.25879,3.44755)(1.39233,3.33759)(1.52586,3.3591)(1.66253,3.46702)(1.7992,3.60154)(1.93587,3.7113)(2.0896,3.76785)(2.26038,3.73934)(2.39868,3.66434)(2.53699,3.56906)(2.67529,3.47734)(2.8149,3.40565)(2.95583,3.36251)(3.09676,3.34766)(3.23769,3.35376)(3.38919,3.37247)(3.54069,3.39287)(3.66231,3.40534)(3.78392,3.41249)(3.90553,3.41438)(4.03526,3.41202)(4.16499,3.40733)(4.29471,3.40252)(4.46981,3.39875)(4.61562,3.39927)(4.76143,3.40313)(4.90724,3.40957)(5.0797,3.41908)(5.2788,3.43075)(5.4779,3.44195)(5.677,3.45245)(5.96088,3.4671)(6.24475,3.48264)(6.52863,3.49944)(6.85863,3.51992)(7.18863,3.54076)(7.51863,3.56179)(7.97361,3.59115)(8.42859,3.62088)(8.88357,3.6508)(9.51072,3.69227)(10.1379,3.7339)(11.2085,3.80511)(12.7226,3.90599)(14.2367,4.0069)(17.4686,4.22236)(20,4.39112) 
};

\addplot [
color=mycolor3,
solid
]
coordinates{
 (0,5)(0.000119652,5)(0.000777737,5)(0.00197426,4.99999)(0.00855511,4.99982)(0.0205203,4.99897)(0.0557841,4.9927)(0.114347,4.97169)(0.172909,4.9411)(0.283485,4.87265)(0.39406,4.81174)(0.504636,4.77451)(0.626408,4.76452)(0.74818,4.77575)(0.869953,4.78706)(0.991725,4.77594)(1.12526,4.71842)(1.25879,4.60622)(1.39233,4.44707)(1.52586,4.26027)(1.66253,4.06581)(1.7992,3.89175)(1.93587,3.75435)(2.0896,3.65088)(2.26038,3.59234)(2.39868,3.57469)(2.53699,3.56961)(2.67529,3.56598)(2.8149,3.55648)(2.95583,3.53815)(3.09676,3.51214)(3.23769,3.48203)(3.38919,3.44998)(3.54069,3.42264)(3.66231,3.40592)(3.78392,3.39427)(3.90553,3.38728)(4.03526,3.38398)(4.16499,3.38365)(4.29471,3.38506)(4.46981,3.38803)(4.61562,3.3906)(4.76143,3.39302)(4.90724,3.39549)(5.0797,3.3989)(5.2788,3.40409)(5.4779,3.41094)(5.677,3.41933)(5.96088,3.4332)(6.24475,3.44831)(6.52863,3.46409)(6.85863,3.48314)(7.18863,3.50296)(7.51863,3.5234)(7.97361,3.55226)(8.42859,3.5816)(8.88357,3.61129)(9.51072,3.65254)(10.1379,3.69404)(11.2085,3.76518)(12.7226,3.866)(14.2367,3.96691)(17.4686,4.18236)(20,4.35112) 
};

\end{axis}
\end{tikzpicture}

%% file: Simulations/consensus_with_I.tikz
%
%
\begin{tikzpicture}

\definecolor{mycolor1}{rgb}{0,0.5,0}
\definecolor{mycolor2}{rgb}{0,0.75,0.75}
\definecolor{mycolor3}{rgb}{0.75,0,0.75}

\begin{axis}[%
scale only axis,
width=\fwidth,
height=\fheight,
xmin=0, xmax=20,
ymin=-6, ymax=8,
xlabel={$t$},
ylabel={$x(t)$},
axis on top]
\addplot [
color=blue,
solid
]
coordinates{
 (0,5)(0.000123851,5)(0.000805034,4.99998)(0.00204355,4.99989)(0.00885537,4.9979)(0.0212405,4.98808)(0.0570874,4.91692)(0.116396,4.67671)(0.175705,4.31563)(0.291346,3.40343)(0.406988,2.42555)(0.52263,1.57898)(0.65114,0.937411)(0.77965,0.664288)(0.903714,0.719173)(1.02333,0.99233)(1.14295,1.3905)(1.28413,1.90561)(1.44686,2.43265)(1.58211,2.75936)(1.71736,2.97581)(1.85261,3.10068)(2.00441,3.17018)(2.17275,3.21085)(2.3411,3.25435)(2.40404,3.27591)(2.46697,3.30095)(2.5299,3.32933)(2.62191,3.37597)(2.74299,3.44377)(2.82556,3.49182)(2.8696,3.51741)(2.91365,3.54271)(2.99394,3.58747)(3.11049,3.64782)(3.19971,3.68941)(3.28893,3.7267)(3.40118,3.76766)(3.53646,3.80918)(3.67175,3.84379)(3.80924,3.87366)(3.92377,3.89551)(3.99726,3.90836)(4.07075,3.92044)(4.20989,3.94134)(4.3861,3.96439)(4.53372,3.98086)(4.68134,3.99488)(4.84535,4.00784)(5.02575,4.01934)(5.20615,4.02842)(5.41978,4.03652)(5.63584,4.04212)(5.8211,4.04502)(6.00636,4.04633)(6.19162,4.0463)(6.38763,4.04513)(6.56071,4.04342)(6.73379,4.0413)(6.90688,4.03893)(7.11926,4.03584)(7.33164,4.03265)(7.54402,4.02946)(7.77245,4.02611)(8.01695,4.02278)(8.26144,4.01987)(8.50593,4.01751)(8.93251,4.01499)(9.31926,4.01463)(9.70602,4.01619)(10.0928,4.01971)(10.7595,4.03038)(11.4263,4.04652)(12.093,4.06763)(13.3076,4.11683)(14.5222,4.17682)(15.7368,4.24442)(17.3434,4.34148)(18.9501,4.44366)(20,4.512) 
};

\addplot [
color=mycolor1,
solid
]
coordinates{
 (0,-6)(0.000123851,-6)(0.000805034,-5.99996)(0.00204355,-5.99974)(0.00885537,-5.99515)(0.0212405,-5.97241)(0.0570874,-5.80766)(0.116396,-5.25106)(0.175705,-4.41284)(0.291346,-2.28386)(0.406988,0.0285773)(0.52263,2.08502)(0.65114,3.74924)(0.77965,4.64448)(0.903714,4.8501)(1.02333,4.60317)(1.14295,4.11594)(1.28413,3.4748)(1.44686,2.90265)(1.58211,2.66473)(1.71736,2.64325)(1.85261,2.78357)(2.00441,3.04059)(2.17275,3.33623)(2.3411,3.56246)(2.40404,3.62227)(2.46697,3.66864)(2.5299,3.70269)(2.62191,3.73365)(2.74299,3.74958)(2.82556,3.75085)(2.8696,3.74999)(2.91365,3.74872)(2.99394,3.74671)(3.11049,3.74766)(3.19971,3.75316)(3.28893,3.76302)(3.40118,3.78067)(3.53646,3.80679)(3.67175,3.83417)(3.80924,3.85979)(3.92377,3.87791)(3.99726,3.88779)(4.07075,3.89636)(4.20989,3.90951)(4.3861,3.92208)(4.53372,3.93078)(4.68134,3.9387)(4.84535,3.9469)(5.02575,3.95486)(5.20615,3.96117)(5.41978,3.96615)(5.63584,3.96884)(5.8211,3.97)(6.00636,3.97073)(6.19162,3.97136)(6.38763,3.97201)(6.56071,3.97249)(6.73379,3.97276)(6.90688,3.97278)(7.11926,3.97255)(7.33164,3.97221)(7.54402,3.972)(7.77245,3.97211)(8.01695,3.97271)(8.26144,3.97382)(8.50593,3.97544)(8.93251,3.9795)(9.31926,3.9847)(9.70602,3.99146)(10.0928,3.99977)(10.7595,4.01757)(11.4263,4.0395)(12.093,4.06515)(13.3076,4.11992)(14.5222,4.18282)(15.7368,4.25156)(17.3434,4.34854)(18.9501,4.44979)(20,4.51738) 
};

\addplot [
color=red,
solid
]
coordinates{
 (0,8)(0.000123851,8)(0.000805034,7.99997)(0.00204355,7.99981)(0.00885537,7.99651)(0.0212405,7.98014)(0.0570874,7.86162)(0.116396,7.46233)(0.175705,6.86464)(0.291346,5.37132)(0.406988,3.80771)(0.52263,2.50676)(0.65114,1.59465)(0.77965,1.29476)(0.903714,1.48815)(1.02333,1.95942)(1.14295,2.53608)(1.28413,3.1561)(1.44686,3.61302)(1.58211,3.74256)(1.71736,3.68545)(1.85261,3.51802)(2.00441,3.2982)(2.17275,3.11316)(2.3411,3.03898)(2.40404,3.03927)(2.46697,3.0521)(2.5299,3.07492)(2.62191,3.12053)(2.74299,3.19003)(2.82556,3.23588)(2.8696,3.25819)(2.91365,3.27847)(2.99394,3.30932)(3.11049,3.33896)(3.19971,3.35013)(3.28893,3.35322)(3.40118,3.34932)(3.53646,3.33972)(3.67175,3.33175)(3.80924,3.33009)(3.92377,3.33487)(3.99726,3.3407)(4.07075,3.34827)(4.20989,3.36585)(4.3861,3.39001)(4.53372,3.40853)(4.68134,3.424)(4.84535,3.4376)(5.02575,3.44971)(5.20615,3.46099)(5.41978,3.47548)(5.63584,3.49212)(5.8211,3.50746)(6.00636,3.523)(6.19162,3.53816)(6.38763,3.55358)(6.56071,3.56676)(6.73379,3.57972)(6.90688,3.59266)(7.11926,3.60864)(7.33164,3.62473)(7.54402,3.64081)(7.77245,3.658)(8.01695,3.67621)(8.26144,3.69426)(8.50593,3.71222)(8.93251,3.74342)(9.31926,3.77153)(9.70602,3.79945)(10.0928,3.82719)(10.7595,3.87468)(11.4263,3.92177)(12.093,3.96853)(13.3076,4.05295)(14.5222,4.1366)(15.7368,4.21966)(17.3434,4.32883)(18.9501,4.43741)(20,4.50813) 
};

\addplot [
color=mycolor2,
solid
]
coordinates{
 (0,4)(0.000123851,4)(0.000805034,4.00001)(0.00204355,4.00005)(0.00885537,4.00097)(0.0212405,4.00551)(0.0570874,4.03831)(0.116396,4.14754)(0.175705,4.30687)(0.291346,4.67629)(0.406988,4.99478)(0.52263,5.15198)(0.65114,5.08328)(0.77965,4.78436)(0.903714,4.36054)(1.02333,3.9229)(1.14295,3.54183)(1.28413,3.24058)(1.44686,3.12614)(1.58211,3.18767)(1.71736,3.32828)(1.85261,3.48181)(2.00441,3.6044)(2.17275,3.63508)(2.3411,3.55809)(2.40404,3.50951)(2.46697,3.45457)(2.5299,3.39627)(2.62191,3.31082)(2.74299,3.2093)(2.82556,3.15297)(2.8696,3.12816)(2.91365,3.10721)(2.99394,3.07893)(3.11049,3.05915)(3.19971,3.05824)(3.28893,3.06625)(3.40118,3.08387)(3.53646,3.10828)(3.67175,3.12869)(3.80924,3.14103)(3.92377,3.14433)(3.99726,3.14362)(4.07075,3.1413)(4.20989,3.13455)(4.3861,3.12661)(4.53372,3.12432)(4.68134,3.1276)(4.84535,3.13751)(5.02575,3.154)(5.20615,3.17346)(5.41978,3.1972)(5.63584,3.22038)(5.8211,3.2397)(6.00636,3.25915)(6.19162,3.2793)(6.38763,3.30164)(6.56071,3.32217)(6.73379,3.34321)(6.90688,3.36452)(7.11926,3.39073)(7.33164,3.41684)(7.54402,3.44279)(7.77245,3.47052)(8.01695,3.49998)(8.26144,3.52915)(8.50593,3.55795)(8.93251,3.60711)(9.31926,3.65038)(9.70602,3.69236)(10.0928,3.73301)(10.7595,3.80007)(11.4263,3.86352)(12.093,3.9237)(13.3076,4.02625)(14.5222,4.1217)(15.7368,4.2121)(17.3434,4.32657)(18.9501,4.43754)(20,4.50896) 
};

\addplot [
color=mycolor3,
solid
]
coordinates{
 (0,5)(0.000123851,5)(0.000805034,5)(0.00204355,4.99999)(0.00885537,4.99981)(0.0212405,4.9989)(0.0570874,4.99234)(0.116396,4.97052)(0.175705,4.93875)(0.291346,4.86517)(0.406988,4.8007)(0.52263,4.76355)(0.65114,4.75715)(0.77965,4.77166)(0.903714,4.7796)(1.02333,4.75739)(1.14295,4.68918)(1.28413,4.54224)(1.44686,4.2982)(1.58211,4.06294)(1.71736,3.82923)(1.85261,3.62281)(2.00441,3.44398)(2.17275,3.31803)(2.3411,3.25553)(2.40404,3.2434)(2.46697,3.23507)(2.5299,3.2291)(2.62191,3.22196)(2.74299,3.21061)(2.82556,3.19929)(2.8696,3.19173)(2.91365,3.18306)(2.99394,3.16449)(3.11049,3.13218)(3.19971,3.10457)(3.28893,3.07605)(3.40118,3.04113)(3.53646,3.00377)(3.67175,2.97443)(3.80924,2.9541)(3.92377,2.94423)(3.99726,2.94087)(4.07075,2.93947)(4.20989,2.94098)(4.3861,2.94787)(4.53372,2.95567)(4.68134,2.96419)(4.84535,2.97419)(5.02575,2.98626)(5.20615,3.00027)(5.41978,3.02017)(5.63584,3.04407)(5.8211,3.06711)(6.00636,3.09183)(6.19162,3.11767)(6.38763,3.14576)(6.56071,3.17099)(6.73379,3.19652)(6.90688,3.22232)(7.11926,3.25424)(7.33164,3.28637)(7.54402,3.31853)(7.77245,3.353)(8.01695,3.38956)(8.26144,3.42563)(8.50593,3.46111)(8.93251,3.52141)(9.31926,3.5741)(9.70602,3.62481)(10.0928,3.67349)(10.7595,3.75273)(11.4263,3.82637)(12.093,3.89493)(13.3076,4.00884)(14.5222,4.11171)(15.7368,4.20677)(17.3434,4.32465)(18.9501,4.43721)(20,4.50911) 
};

\end{axis}
\end{tikzpicture}

%% file: Simulations/consensus_with_I_too_large.tikz
%
%
\begin{tikzpicture}

\definecolor{mycolor1}{rgb}{0,0.5,0}
\definecolor{mycolor2}{rgb}{0,0.75,0.75}
\definecolor{mycolor3}{rgb}{0.75,0,0.75}

\begin{axis}[%
scale only axis,
width=\fwidth,
height=\fheight,
xmin=0, xmax=10,
ymin=-6, ymax=14,
xlabel={$t$},
ylabel={$x(t)$},
axis on top]
\addplot [
color=blue,
solid
]
coordinates{
 (0,5)(0.00046341,4.99999)(0.00301216,4.99965)(0.00764626,4.99814)(0.0258359,4.9788)(0.0449307,4.94047)(0.0649307,4.88151)(0.0849307,4.80187)(0.104931,4.70184)(0.124931,4.58209)(0.144931,4.44339)(0.164931,4.28665)(0.184931,4.11289)(0.204931,3.92325)(0.224931,3.71896)(0.244931,3.50139)(0.264931,3.27194)(0.284931,3.03215)(0.304931,2.78359)(0.324931,2.52791)(0.344931,2.2668)(0.364931,2.00201)(0.384931,1.73528)(0.404931,1.46841)(0.424931,1.20317)(0.444931,0.941352)(0.464931,0.684717)(0.484931,0.434999)(0.504931,0.193894)(0.524931,-0.0369501)(0.544931,-0.255948)(0.564931,-0.461583)(0.584931,-0.652418)(0.604931,-0.827107)(0.624931,-0.984397)(0.644931,-1.12314)(0.664931,-1.24231)(0.684931,-1.34098)(0.704931,-1.41837)(0.724931,-1.47381)(0.744931,-1.50676)(0.764931,-1.51684)(0.784931,-1.50379)(0.804931,-1.46748)(0.824931,-1.40794)(0.844931,-1.32533)(0.864931,-1.21995)(0.884931,-1.09223)(0.904931,-0.942744)(0.924931,-0.772181)(0.944931,-0.581364)(0.964931,-0.371232)(0.984931,-0.142836)(1.00493,0.102668)(1.02493,0.364022)(1.04493,0.639883)(1.06493,0.928824)(1.08493,1.22934)(1.10493,1.53988)(1.12493,1.85883)(1.14493,2.18451)(1.16493,2.51526)(1.18493,2.84935)(1.20493,3.18506)(1.22493,3.52067)(1.24493,3.85446)(1.26493,4.18475)(1.28493,4.50986)(1.30493,4.82819)(1.32493,5.13816)(1.34493,5.43826)(1.36493,5.72705)(1.38493,6.00317)(1.40493,6.26533)(1.42493,6.51235)(1.44493,6.74315)(1.46493,6.95674)(1.48493,7.15226)(1.50493,7.32894)(1.52493,7.48614)(1.54493,7.62335)(1.56493,7.74018)(1.58493,7.83636)(1.60493,7.91175)(1.62493,7.96633)(1.64493,8.0002)(1.66493,8.0136)(1.68493,8.00687)(1.70493,7.98049)(1.72493,7.93503)(1.74493,7.87117)(1.76493,7.78971)(1.78493,7.69153)(1.80493,7.57759)(1.82493,7.44896)(1.84493,7.30677)(1.86493,7.15221)(1.88493,6.98654)(1.90493,6.81107)(1.92493,6.62713)(1.94493,6.4361)(1.96493,6.23939)(1.98493,6.03841)(2.00493,5.83457)(2.02493,5.62929)(2.04493,5.42395)(2.06493,5.21994)(2.08493,5.01859)(2.10493,4.82121)(2.12493,4.62904)(2.14493,4.44328)(2.16493,4.26507)(2.18493,4.09546)(2.20493,3.93545)(2.22493,3.78594)(2.24493,3.64775)(2.26493,3.52161)(2.28493,3.40814)(2.30493,3.3079)(2.32493,3.2213)(2.34493,3.14869)(2.36493,3.0903)(2.38493,3.04626)(2.40493,3.01659)(2.42493,3.00122)(2.44493,2.99997)(2.46493,3.01257)(2.48493,3.03867)(2.50493,3.07779)(2.52493,3.1294)(2.54493,3.19288)(2.56493,3.26751)(2.58493,3.35254)(2.60493,3.44713)(2.62493,3.55037)(2.64493,3.66133)(2.66493,3.77901)(2.68493,3.90239)(2.70493,4.03041)(2.72493,4.16198)(2.74493,4.29601)(2.76493,4.4314)(2.78493,4.56705)(2.80493,4.70185)(2.82493,4.83474)(2.84493,4.96464)(2.86493,5.09054)(2.88493,5.21143)(2.90493,5.32637)(2.92493,5.43444)(2.94493,5.5348)(2.96493,5.62665)(2.98493,5.70925)(3.00493,5.78195)(3.02493,5.84414)(3.04493,5.8953)(3.06493,5.93498)(3.08493,5.96282)(3.10493,5.97852)(3.12493,5.98188)(3.14493,5.97276)(3.16493,5.95111)(3.18493,5.91696)(3.20493,5.87041)(3.22493,5.81166)(3.24493,5.74095)(3.26493,5.65861)(3.28493,5.56503)(3.30493,5.46068)(3.32493,5.34607)(3.34493,5.22177)(3.36493,5.08842)(3.38493,4.94667)(3.40493,4.79725)(3.42493,4.6409)(3.44493,4.4784)(3.46493,4.31054)(3.48493,4.13816)(3.50493,3.96208)(3.52493,3.78315)(3.54493,3.6022)(3.56493,3.42008)(3.58493,3.23762)(3.60493,3.05564)(3.62493,2.87492)(3.64493,2.69625)(3.66493,2.52037)(3.68493,2.34799)(3.70493,2.17979)(3.72493,2.0164)(3.74493,1.85841)(3.76493,1.70637)(3.78493,1.56079)(3.80493,1.4221)(3.82493,1.2907)(3.84493,1.16695)(3.86493,1.05112)(3.88493,0.943463)(3.90493,0.844151)(3.92493,0.753318)(3.94493,0.671041)(3.96493,0.597346)(3.98493,0.532211)(4.00493,0.475564)(4.02493,0.42729)(4.04493,0.387231)(4.06493,0.355188)(4.08493,0.330926)(4.10493,0.314176)(4.12493,0.304639)(4.14493,0.301988)(4.16493,0.305875)(4.18493,0.31593)(4.20493,0.33177)(4.22493,0.352997)(4.24493,0.379208)(4.26493,0.409992)(4.28493,0.444939)(4.30493,0.483643)(4.32493,0.525699)(4.34493,0.570716)(4.36493,0.618312)(4.38493,0.66812)(4.40493,0.719791)(4.42493,0.772995)(4.44493,0.827422)(4.46493,0.882789)(4.48493,0.938832)(4.50493,0.995318)(4.52493,1.05204)(4.54493,1.10881)(4.56493,1.16548)(4.58493,1.22192)(4.60493,1.27803)(4.62493,1.33373)(4.64493,1.38899)(4.66493,1.44376)(4.68493,1.49806)(4.70493,1.5519)(4.72493,1.60532)(4.74493,1.65838)(4.76493,1.71114)(4.78493,1.7637)(4.80493,1.81614)(4.82493,1.86858)(4.84493,1.92111)(4.86493,1.97385)(4.88493,2.02691)(4.90493,2.08039)(4.92493,2.13441)(4.94493,2.18906)(4.96493,2.24443)(4.98493,2.30059)(5.00493,2.35761)(5.02493,2.41553)(5.04493,2.47438)(5.06493,2.53417)(5.08493,2.5949)(5.10493,2.65653)(5.12493,2.71902)(5.14493,2.78228)(5.16493,2.84624)(5.18493,2.91076)(5.20493,2.97571)(5.22493,3.04094)(5.24493,3.10625)(5.26493,3.17146)(5.28493,3.23635)(5.30493,3.30069)(5.32493,3.36423)(5.34493,3.42671)(5.36493,3.48787)(5.38493,3.54744)(5.40493,3.60514)(5.42493,3.66069)(5.44493,3.71382)(5.46493,3.76425)(5.48493,3.81172)(5.50493,3.85597)(5.52493,3.89678)(5.54493,3.9339)(5.56493,3.96714)(5.58493,3.99633)(5.60493,4.02129)(5.62493,4.04191)(5.64493,4.05809)(5.66493,4.06976)(5.68493,4.07689)(5.70493,4.07949)(5.72493,4.07759)(5.74493,4.07129)(5.76493,4.0607)(5.78493,4.04598)(5.80493,4.02733)(5.82493,4.005)(5.84493,3.97927)(5.86493,3.95044)(5.88493,3.91889)(5.90493,3.885)(5.92493,3.8492)(5.94493,3.81194)(5.96493,3.77371)(5.98493,3.73503)(6.00493,3.69642)(6.02493,3.65845)(6.04493,3.62168)(6.06493,3.58669)(6.08493,3.55407)(6.10493,3.5244)(6.12493,3.49829)(6.14493,3.47631)(6.16493,3.45903)(6.18493,3.44701)(6.20493,3.44078)(6.22493,3.44085)(6.24493,3.44771)(6.26493,3.46179)(6.28493,3.4835)(6.30493,3.51321)(6.32493,3.55122)(6.34493,3.59781)(6.36493,3.65319)(6.38493,3.7175)(6.40493,3.79083)(6.42493,3.87323)(6.44493,3.96464)(6.46493,4.06497)(6.48493,4.17404)(6.50493,4.29163)(6.52493,4.41741)(6.54493,4.55103)(6.56493,4.69202)(6.58493,4.83989)(6.60493,4.99406)(6.62493,5.1539)(6.64493,5.31871)(6.66493,5.48773)(6.68493,5.66018)(6.70493,5.83518)(6.72493,6.01186)(6.74493,6.18927)(6.76493,6.36644)(6.78493,6.54239)(6.80493,6.71609)(6.82493,6.88651)(6.84493,7.05261)(6.86493,7.21334)(6.88493,7.36767)(6.90493,7.51457)(6.92493,7.65302)(6.94493,7.78204)(6.96493,7.90067)(6.98493,8.00801)(7.00493,8.10318)(7.02493,8.18537)(7.04493,8.25381)(7.06493,8.30781)(7.08493,8.34675)(7.10493,8.37007)(7.12493,8.3773)(7.14493,8.36806)(7.16493,8.34205)(7.18493,8.29905)(7.20493,8.23895)(7.22493,8.16175)(7.24493,8.0675)(7.26493,7.9564)(7.28493,7.82872)(7.30493,7.68484)(7.32493,7.52523)(7.34493,7.35047)(7.36493,7.16121)(7.38493,6.95821)(7.40493,6.74232)(7.42493,6.51446)(7.44493,6.27564)(7.46493,6.02694)(7.48493,5.76951)(7.50493,5.50454)(7.52493,5.23331)(7.54493,4.95713)(7.56493,4.67734)(7.58493,4.39533)(7.60493,4.11249)(7.62493,3.83026)(7.64493,3.55005)(7.66493,3.27329)(7.68493,3.00139)(7.70493,2.73575)(7.72493,2.47772)(7.74493,2.22862)(7.76493,1.98975)(7.78493,1.76232)(7.80493,1.54749)(7.82493,1.34636)(7.84493,1.15993)(7.86493,0.989145)(7.88493,0.834829)(7.90493,0.697727)(7.92493,0.578478)(7.94493,0.477613)(7.96493,0.395554)(7.98493,0.332607)(8.00493,0.288963)(8.02493,0.264693)(8.04493,0.25975)(8.06493,0.273965)(8.08493,0.307052)(8.10493,0.358603)(8.12493,0.428098)(8.14493,0.514901)(8.16493,0.618265)(8.18493,0.73734)(8.20493,0.871172)(8.22493,1.01871)(8.24493,1.17883)(8.26493,1.3503)(8.28493,1.53182)(8.30493,1.72205)(8.32493,1.91956)(8.34493,2.12289)(8.36493,2.33052)(8.38493,2.54092)(8.40493,2.75253)(8.42493,2.96379)(8.44493,3.17311)(8.46493,3.37895)(8.48493,3.57977)(8.50493,3.77405)(8.52493,3.96032)(8.54493,4.13717)(8.56493,4.30323)(8.58493,4.45722)(8.60493,4.59791)(8.62493,4.72419)(8.64493,4.83501)(8.66493,4.92944)(8.68493,5.00664)(8.70493,5.06592)(8.72493,5.10665)(8.74493,5.12838)(8.76493,5.13074)(8.78493,5.11353)(8.80493,5.07665)(8.82493,5.02014)(8.84493,4.94418)(8.86493,4.84908)(8.88493,4.73528)(8.90493,4.60334)(8.92493,4.45395)(8.94493,4.28793)(8.96493,4.10621)(8.98493,3.90983)(9.00493,3.69993)(9.02493,3.47777)(9.04493,3.24468)(9.06493,3.00207)(9.08493,2.75145)(9.10493,2.49435)(9.12493,2.2324)(9.14493,1.96725)(9.16493,1.70058)(9.18493,1.43412)(9.20493,1.16958)(9.22493,0.908692)(9.24493,0.653173)(9.26493,0.404719)(9.28493,0.164994)(9.30493,-0.0643793)(9.32493,-0.281831)(9.34493,-0.485854)(9.36493,-0.675014)(9.38493,-0.847959)(9.40493,-1.00343)(9.42493,-1.14026)(9.44493,-1.2574)(9.46493,-1.3539)(9.48493,-1.42896)(9.50493,-1.48186)(9.52493,-1.51205)(9.54493,-1.5191)(9.56493,-1.50272)(9.58493,-1.46275)(9.60493,-1.39919)(9.62493,-1.31217)(9.64493,-1.20195)(9.66493,-1.06896)(9.68493,-0.913735)(9.70493,-0.736966)(9.72493,-0.539466)(9.74493,-0.322175)(9.76493,-0.0861528)(9.78493,0.167428)(9.80493,0.437288)(9.82493,0.722049)(9.84493,1.02025)(9.86493,1.33033)(9.88493,1.65068)(9.90493,1.97961)(9.92493,2.3154)(9.94493,2.65627)(9.96493,3.00042)(9.98493,3.34604)(10,3.60631) 
};

\addplot [
color=mycolor1,
solid
]
coordinates{
 (0,-6)(0.00046341,-5.99998)(0.00301216,-5.9992)(0.00764626,-5.99569)(0.0258359,-5.95095)(0.0449307,-5.86229)(0.0649307,-5.72597)(0.0849307,-5.54183)(0.104931,-5.31054)(0.124931,-5.03362)(0.144931,-4.71282)(0.164931,-4.35014)(0.184931,-3.94787)(0.204931,-3.50854)(0.224931,-3.0349)(0.244931,-2.52991)(0.264931,-1.99675)(0.284931,-1.43875)(0.304931,-0.859384)(0.324931,-0.262272)(0.344931,0.348871)(0.364931,0.97025)(0.384931,1.59802)(0.404931,2.22829)(0.424931,2.85718)(0.444931,3.48082)(0.464931,4.09541)(0.484931,4.69718)(0.504931,5.28249)(0.524931,5.84781)(0.544931,6.38972)(0.564931,6.90501)(0.584931,7.39063)(0.604931,7.84372)(0.624931,8.26167)(0.644931,8.64208)(0.664931,8.98282)(0.684931,9.28202)(0.704931,9.53808)(0.724931,9.74969)(0.744931,9.91585)(0.764931,10.0358)(0.784931,10.1092)(0.804931,10.1359)(0.824931,10.1161)(0.844931,10.0503)(0.864931,9.93932)(0.884931,9.78421)(0.904931,9.58636)(0.924931,9.34741)(0.944931,9.06927)(0.964931,8.7541)(0.984931,8.40427)(1.00493,8.0224)(1.02493,7.61128)(1.04493,7.17389)(1.06493,6.71337)(1.08493,6.23298)(1.10493,5.73612)(1.12493,5.22627)(1.14493,4.70697)(1.16493,4.18182)(1.18493,3.65443)(1.20493,3.12841)(1.22493,2.60735)(1.24493,2.09477)(1.26493,1.59413)(1.28493,1.1088)(1.30493,0.642017)(1.32493,0.196886)(1.34493,-0.223649)(1.36493,-0.616819)(1.38493,-0.980054)(1.40493,-1.31099)(1.42493,-1.6075)(1.44493,-1.86768)(1.46493,-2.08989)(1.48493,-2.27275)(1.50493,-2.41513)(1.52493,-2.5162)(1.54493,-2.57539)(1.56493,-2.59241)(1.58493,-2.56726)(1.60493,-2.50021)(1.62493,-2.39182)(1.64493,-2.2429)(1.66493,-2.05454)(1.68493,-1.82808)(1.70493,-1.5651)(1.72493,-1.26742)(1.74493,-0.93707)(1.76493,-0.576289)(1.78493,-0.187505)(1.80493,0.226686)(1.82493,0.663533)(1.84493,1.12016)(1.86493,1.59356)(1.88493,2.08065)(1.90493,2.57828)(1.92493,3.08323)(1.94493,3.59227)(1.96493,4.10215)(1.98493,4.60965)(2.00493,5.11158)(2.02493,5.60481)(2.04493,6.08627)(2.06493,6.55303)(2.08493,7.00222)(2.10493,7.43116)(2.12493,7.83728)(2.14493,8.2182)(2.16493,8.57173)(2.18493,8.89585)(2.20493,9.18877)(2.22493,9.4489)(2.24493,9.67489)(2.26493,9.86562)(2.28493,10.0202)(2.30493,10.138)(2.32493,10.2186)(2.34493,10.2619)(2.36493,10.268)(2.38493,10.2371)(2.40493,10.17)(2.42493,10.0673)(2.44493,9.93014)(2.46493,9.75972)(2.48493,9.5575)(2.50493,9.32509)(2.52493,9.06433)(2.54493,8.77719)(2.56493,8.4658)(2.58493,8.13241)(2.60493,7.7794)(2.62493,7.40926)(2.64493,7.02453)(2.66493,6.62783)(2.68493,6.22184)(2.70493,5.80923)(2.72493,5.39272)(2.74493,4.97497)(2.76493,4.55864)(2.78493,4.14634)(2.80493,3.74062)(2.82493,3.34392)(2.84493,2.95861)(2.86493,2.58694)(2.88493,2.23103)(2.90493,1.89285)(2.92493,1.57424)(2.94493,1.27686)(2.96493,1.0022)(2.98493,0.751565)(3.00493,0.526089)(3.02493,0.326697)(3.04493,0.154125)(3.06493,0.00890578)(3.08493,-0.108626)(3.10493,-0.198337)(3.12493,-0.26029)(3.14493,-0.294745)(3.16493,-0.302149)(3.18493,-0.283137)(3.20493,-0.238518)(3.22493,-0.169271)(3.24493,-0.0765337)(3.26493,0.0384074)(3.28493,0.174129)(3.30493,0.329083)(3.32493,0.50161)(3.34493,0.689954)(3.36493,0.892272)(3.38493,1.10666)(3.40493,1.33115)(3.42493,1.56374)(3.44493,1.80241)(3.46493,2.04513)(3.48493,2.28987)(3.50493,2.53464)(3.52493,2.77748)(3.54493,3.01646)(3.56493,3.24975)(3.58493,3.47558)(3.60493,3.69226)(3.62493,3.89823)(3.64493,4.09203)(3.66493,4.2723)(3.68493,4.43783)(3.70493,4.58756)(3.72493,4.72054)(3.74493,4.836)(3.76493,4.93331)(3.78493,5.01198)(3.80493,5.07169)(3.82493,5.11229)(3.84493,5.13375)(3.86493,5.13624)(3.88493,5.12002)(3.90493,5.08555)(3.92493,5.03339)(3.94493,4.96425)(3.96493,4.87897)(3.98493,4.7785)(4.00493,4.66387)(4.02493,4.53626)(4.04493,4.39689)(4.06493,4.24707)(4.08493,4.08818)(4.10493,3.92164)(4.12493,3.74891)(4.14493,3.57148)(4.16493,3.39085)(4.18493,3.20853)(4.20493,3.026)(4.22493,2.84472)(4.24493,2.66613)(4.26493,2.49161)(4.28493,2.32248)(4.30493,2.16)(4.32493,2.00534)(4.34493,1.85959)(4.36493,1.72375)(4.38493,1.5987)(4.40493,1.48524)(4.42493,1.38404)(4.44493,1.29563)(4.46493,1.22047)(4.48493,1.15884)(4.50493,1.11095)(4.52493,1.07684)(4.54493,1.05644)(4.56493,1.04958)(4.58493,1.05596)(4.60493,1.07514)(4.62493,1.1066)(4.64493,1.14972)(4.66493,1.20377)(4.68493,1.26792)(4.70493,1.3413)(4.72493,1.42292)(4.74493,1.51176)(4.76493,1.60675)(4.78493,1.70676)(4.80493,1.81064)(4.82493,1.91722)(4.84493,2.02533)(4.86493,2.13378)(4.88493,2.24141)(4.90493,2.34709)(4.92493,2.44972)(4.94493,2.54824)(4.96493,2.64165)(4.98493,2.72902)(5.00493,2.80949)(5.02493,2.88228)(5.04493,2.94671)(5.06493,3.00217)(5.08493,3.04818)(5.10493,3.08436)(5.12493,3.11043)(5.14493,3.12624)(5.16493,3.13175)(5.18493,3.12701)(5.20493,3.11223)(5.22493,3.08771)(5.24493,3.05386)(5.26493,3.01121)(5.28493,2.9604)(5.30493,2.90214)(5.32493,2.83728)(5.34493,2.76671)(5.36493,2.69144)(5.38493,2.61252)(5.40493,2.53108)(5.42493,2.44828)(5.44493,2.36535)(5.46493,2.28353)(5.48493,2.20409)(5.50493,2.12828)(5.52493,2.0574)(5.54493,1.99267)(5.56493,1.93533)(5.58493,1.88656)(5.60493,1.84749)(5.62493,1.8192)(5.64493,1.80268)(5.66493,1.79885)(5.68493,1.80854)(5.70493,1.83245)(5.72493,1.87122)(5.74493,1.92534)(5.76493,1.99517)(5.78493,2.08097)(5.80493,2.18285)(5.82493,2.30078)(5.84493,2.4346)(5.86493,2.584)(5.88493,2.74854)(5.90493,2.92764)(5.92493,3.12058)(5.94493,3.3265)(5.96493,3.54442)(5.98493,3.77322)(6.00493,4.0117)(6.02493,4.2585)(6.04493,4.51221)(6.06493,4.77128)(6.08493,5.03412)(6.10493,5.29905)(6.12493,5.56433)(6.14493,5.82818)(6.16493,6.0888)(6.18493,6.34436)(6.20493,6.59302)(6.22493,6.83295)(6.24493,7.06236)(6.26493,7.27948)(6.28493,7.48261)(6.30493,7.67009)(6.32493,7.84036)(6.34493,7.99195)(6.36493,8.12348)(6.38493,8.23372)(6.40493,8.32153)(6.42493,8.38594)(6.44493,8.4261)(6.46493,8.44135)(6.48493,8.43116)(6.50493,8.39519)(6.52493,8.33328)(6.54493,8.24544)(6.56493,8.13187)(6.58493,7.99295)(6.60493,7.82923)(6.62493,7.64147)(6.64493,7.43058)(6.66493,7.19768)(6.68493,6.94403)(6.70493,6.67106)(6.72493,6.38036)(6.74493,6.07366)(6.76493,5.75284)(6.78493,5.41989)(6.80493,5.07691)(6.82493,4.7261)(6.84493,4.36975)(6.86493,4.0102)(6.88493,3.64987)(6.90493,3.29119)(6.92493,2.93662)(6.94493,2.58862)(6.96493,2.24963)(6.98493,1.92206)(7.00493,1.60827)(7.02493,1.31055)(7.04493,1.03111)(7.06493,0.772048)(7.08493,0.535354)(7.10493,0.322879)(7.12493,0.13633)(7.14493,-0.0227497)(7.16493,-0.152988)(7.18493,-0.253201)(7.20493,-0.322397)(7.22493,-0.359792)(7.24493,-0.364812)(7.26493,-0.337103)(7.28493,-0.276536)(7.30493,-0.183204)(7.32493,-0.0574345)(7.34493,0.100221)(7.36493,0.288982)(7.38493,0.507844)(7.40493,0.755583)(7.42493,1.03076)(7.44493,1.33174)(7.46493,1.65668)(7.48493,2.00355)(7.50493,2.37016)(7.52493,2.75415)(7.54493,3.15302)(7.56493,3.56413)(7.58493,3.98475)(7.60493,4.41204)(7.62493,4.84307)(7.64493,5.27488)(7.66493,5.70446)(7.68493,6.12879)(7.70493,6.54483)(7.72493,6.9496)(7.74493,7.34014)(7.76493,7.71356)(7.78493,8.06706)(7.80493,8.39794)(7.82493,8.70364)(7.84493,8.98173)(7.86493,9.22992)(7.88493,9.44615)(7.90493,9.62851)(7.92493,9.77531)(7.94493,9.88509)(7.96493,9.9566)(7.98493,9.98887)(8.00493,9.98113)(8.02493,9.93293)(8.04493,9.84403)(8.06493,9.71448)(8.08493,9.5446)(8.10493,9.33498)(8.12493,9.08647)(8.14493,8.80019)(8.16493,8.47751)(8.18493,8.12006)(8.20493,7.72971)(8.22493,7.30854)(8.24493,6.85888)(8.26493,6.38324)(8.28493,5.88434)(8.30493,5.36506)(8.32493,4.82844)(8.34493,4.27763)(8.36493,3.71594)(8.38493,3.14674)(8.40493,2.57349)(8.42493,1.99969)(8.44493,1.42886)(8.46493,0.864556)(8.48493,0.310278)(8.50493,-0.230495)(8.52493,-0.754357)(8.54493,-1.25799)(8.56493,-1.73818)(8.58493,-2.19185)(8.60493,-2.61609)(8.62493,-3.00813)(8.64493,-3.36542)(8.66493,-3.68562)(8.68493,-3.96659)(8.70493,-4.20647)(8.72493,-4.40361)(8.74493,-4.55666)(8.76493,-4.66453)(8.78493,-4.72643)(8.80493,-4.74183)(8.82493,-4.71053)(8.84493,-4.63261)(8.86493,-4.50845)(8.88493,-4.33873)(8.90493,-4.12441)(8.92493,-3.86676)(8.94493,-3.56732)(8.96493,-3.22789)(8.98493,-2.85054)(9.00493,-2.43761)(9.02493,-1.99162)(9.04493,-1.51537)(9.06493,-1.01181)(9.08493,-0.484088)(9.10493,0.0644789)(9.12493,0.630445)(9.14493,1.21025)(9.16493,1.80022)(9.18493,2.39665)(9.20493,2.99574)(9.22493,3.59371)(9.24493,4.18677)(9.26493,4.77114)(9.28493,5.34313)(9.30493,5.89909)(9.32493,6.43551)(9.34493,6.94897)(9.36493,7.43621)(9.38493,7.89416)(9.40493,8.3199)(9.42493,8.71076)(9.44493,9.06426)(9.46493,9.37819)(9.48493,9.65058)(9.50493,9.87976)(9.52493,10.0643)(9.54493,10.2031)(9.56493,10.2954)(9.58493,10.3405)(9.60493,10.3384)(9.62493,10.2892)(9.64493,10.1932)(9.66493,10.0512)(9.68493,9.86427)(9.70493,9.63367)(9.72493,9.36103)(9.74493,9.04825)(9.76493,8.69749)(9.78493,8.31116)(9.80493,7.8919)(9.82493,7.44258)(9.84493,6.96624)(9.86493,6.46615)(9.88493,5.94568)(9.90493,5.40837)(9.92493,4.85787)(9.94493,4.29791)(9.96493,3.73229)(9.98493,3.16483)(10,2.73839) 
};

\addplot [
color=red,
solid
]
coordinates{
 (0,8)(0.00046341,7.99999)(0.00301216,7.99942)(0.00764626,7.99689)(0.0258359,7.96469)(0.0449307,7.90092)(0.0649307,7.80293)(0.0849307,7.67072)(0.104931,7.50489)(0.124931,7.30673)(0.144931,7.07766)(0.164931,6.81938)(0.184931,6.53378)(0.204931,6.22295)(0.224931,5.88917)(0.244931,5.5349)(0.264931,5.16273)(0.284931,4.7754)(0.304931,4.37575)(0.324931,3.96673)(0.344931,3.55133)(0.364931,3.13261)(0.384931,2.71367)(0.404931,2.29757)(0.424931,1.8874)(0.444931,1.48617)(0.464931,1.09686)(0.484931,0.722339)(0.504931,0.365383)(0.524931,0.0286408)(0.544931,-0.285385)(0.564931,-0.574355)(0.584931,-0.836109)(0.604931,-1.06868)(0.624931,-1.27033)(0.644931,-1.43951)(0.664931,-1.57494)(0.684931,-1.67556)(0.704931,-1.74059)(0.724931,-1.76946)(0.744931,-1.76192)(0.764931,-1.71793)(0.784931,-1.63776)(0.804931,-1.52192)(0.824931,-1.37118)(0.844931,-1.18656)(0.864931,-0.969339)(0.884931,-0.721024)(0.904931,-0.443354)(0.924931,-0.138279)(0.944931,0.19205)(0.964931,0.545296)(0.984931,0.918956)(1.00493,1.31037)(1.02493,1.71676)(1.04493,2.13523)(1.06493,2.56279)(1.08493,2.99638)(1.10493,3.43291)(1.12493,3.86924)(1.14493,4.30225)(1.16493,4.72884)(1.18493,5.14593)(1.20493,5.55054)(1.22493,5.93974)(1.24493,6.31075)(1.26493,6.66087)(1.28493,6.98759)(1.30493,7.28855)(1.32493,7.56156)(1.34493,7.80466)(1.36493,8.01606)(1.38493,8.19424)(1.40493,8.3379)(1.42493,8.44597)(1.44493,8.51765)(1.46493,8.55241)(1.48493,8.54996)(1.50493,8.51031)(1.52493,8.4337)(1.54493,8.32065)(1.56493,8.17195)(1.58493,7.98864)(1.60493,7.772)(1.62493,7.52356)(1.64493,7.24507)(1.66493,6.9385)(1.68493,6.60602)(1.70493,6.24999)(1.72493,5.87294)(1.74493,5.47753)(1.76493,5.06659)(1.78493,4.64303)(1.80493,4.20986)(1.82493,3.77016)(1.84493,3.32706)(1.86493,2.88372)(1.88493,2.44328)(1.90493,2.00888)(1.92493,1.58361)(1.94493,1.17049)(1.96493,0.77247)(1.98493,0.392368)(2.00493,0.0328888)(2.02493,-0.303414)(2.04493,-0.614154)(2.06493,-0.897126)(2.08493,-1.15033)(2.10493,-1.37196)(2.12493,-1.56047)(2.14493,-1.71452)(2.16493,-1.83303)(2.18493,-1.91518)(2.20493,-1.9604)(2.22493,-1.96839)(2.24493,-1.93911)(2.26493,-1.87279)(2.28493,-1.76994)(2.30493,-1.6313)(2.32493,-1.4579)(2.34493,-1.25099)(2.36493,-1.01207)(2.38493,-0.742884)(2.40493,-0.445365)(2.42493,-0.121661)(2.44493,0.225893)(2.46493,0.594796)(2.48493,0.982396)(2.50493,1.38591)(2.52493,1.80243)(2.54493,2.22898)(2.56493,2.6625)(2.58493,3.09989)(2.60493,3.538)(2.62493,3.97373)(2.64493,4.40394)(2.66493,4.82558)(2.68493,5.23563)(2.70493,5.63119)(2.72493,6.00944)(2.74493,6.36769)(2.76493,6.70342)(2.78493,7.01425)(2.80493,7.29799)(2.82493,7.55265)(2.84493,7.77645)(2.86493,7.96785)(2.88493,8.12551)(2.90493,8.24837)(2.92493,8.3356)(2.94493,8.38665)(2.96493,8.40122)(2.98493,8.37926)(3.00493,8.32102)(3.02493,8.22698)(3.04493,8.0979)(3.06493,7.93477)(3.08493,7.73886)(3.10493,7.51165)(3.12493,7.25484)(3.14493,6.97038)(3.16493,6.66039)(3.18493,6.32717)(3.20493,5.9732)(3.22493,5.6011)(3.24493,5.21364)(3.26493,4.81367)(3.28493,4.40415)(3.30493,3.9881)(3.32493,3.56859)(3.34493,3.14869)(3.36493,2.7315)(3.38493,2.32009)(3.40493,1.91747)(3.42493,1.52661)(3.44493,1.15036)(3.46493,0.791502)(3.48493,0.452651)(3.50493,0.136293)(3.52493,-0.155254)(3.54493,-0.419855)(3.56493,-0.655573)(3.58493,-0.860684)(3.60493,-1.03369)(3.62493,-1.17332)(3.64493,-1.27856)(3.66493,-1.34864)(3.68493,-1.38305)(3.70493,-1.38154)(3.72493,-1.34412)(3.74493,-1.27107)(3.76493,-1.16293)(3.78493,-1.02049)(3.80493,-0.844799)(3.82493,-0.637141)(3.84493,-0.399041)(3.86493,-0.132244)(3.88493,0.161295)(3.90493,0.479422)(3.92493,0.819807)(3.94493,1.17995)(3.96493,1.55722)(3.98493,1.94884)(4.00493,2.35195)(4.02493,2.76358)(4.04493,3.18072)(4.06493,3.60029)(4.08493,4.01923)(4.10493,4.43445)(4.12493,4.84291)(4.14493,5.24159)(4.16493,5.62757)(4.18493,5.99799)(4.20493,6.35014)(4.22493,6.68141)(4.24493,6.98935)(4.26493,7.27169)(4.28493,7.52633)(4.30493,7.75137)(4.32493,7.94514)(4.34493,8.10619)(4.36493,8.23331)(4.38493,8.32552)(4.40493,8.38211)(4.42493,8.40263)(4.44493,8.38688)(4.46493,8.33494)(4.48493,8.24714)(4.50493,8.12408)(4.52493,7.96661)(4.54493,7.77584)(4.56493,7.55311)(4.58493,7.3)(4.60493,7.01832)(4.62493,6.71009)(4.64493,6.3775)(4.66493,6.02295)(4.68493,5.64898)(4.70493,5.25829)(4.72493,4.85369)(4.74493,4.43811)(4.76493,4.01453)(4.78493,3.58604)(4.80493,3.15572)(4.82493,2.7267)(4.84493,2.30209)(4.86493,1.88498)(4.88493,1.47839)(4.90493,1.08529)(4.92493,0.708538)(4.94493,0.350889)(4.96493,0.0149588)(4.98493,-0.296793)(5.00493,-0.582076)(5.02493,-0.838789)(5.04493,-1.06504)(5.06493,-1.25913)(5.08493,-1.41964)(5.10493,-1.54533)(5.12493,-1.63526)(5.14493,-1.68872)(5.16493,-1.70527)(5.18493,-1.68472)(5.20493,-1.62717)(5.22493,-1.53298)(5.24493,-1.40275)(5.26493,-1.23737)(5.28493,-1.03796)(5.30493,-0.805911)(5.32493,-0.542816)(5.34493,-0.250515)(5.36493,0.0689504)(5.38493,0.41334)(5.40493,0.780235)(5.42493,1.16705)(5.44493,1.57107)(5.46493,1.98943)(5.48493,2.41919)(5.50493,2.8573)(5.52493,3.30067)(5.54493,3.74616)(5.56493,4.19062)(5.58493,4.6309)(5.60493,5.06388)(5.62493,5.4865)(5.64493,5.89576)(5.66493,6.28877)(5.68493,6.66274)(5.70493,7.01502)(5.72493,7.34312)(5.74493,7.64473)(5.76493,7.9177)(5.78493,8.16013)(5.80493,8.37029)(5.82493,8.54671)(5.84493,8.68815)(5.86493,8.79364)(5.88493,8.86244)(5.90493,8.89407)(5.92493,8.88836)(5.94493,8.84535)(5.96493,8.76538)(5.98493,8.64905)(6.00493,8.49722)(6.02493,8.31099)(6.04493,8.09173)(6.06493,7.84101)(6.08493,7.56066)(6.10493,7.25271)(6.12493,6.91937)(6.14493,6.56305)(6.16493,6.18631)(6.18493,5.79187)(6.20493,5.38257)(6.22493,4.96134)(6.24493,4.5312)(6.26493,4.09525)(6.28493,3.6566)(6.30493,3.21841)(6.32493,2.7838)(6.34493,2.35588)(6.36493,1.93772)(6.38493,1.53229)(6.40493,1.14249)(6.42493,0.771094)(6.44493,0.42074)(6.46493,0.0939181)(6.48493,-0.207056)(6.50493,-0.480055)(6.52493,-0.723154)(6.54493,-0.934646)(6.56493,-1.11305)(6.58493,-1.25713)(6.60493,-1.3659)(6.62493,-1.43862)(6.64493,-1.47481)(6.66493,-1.47428)(6.68493,-1.43707)(6.70493,-1.3635)(6.72493,-1.25416)(6.74493,-1.10989)(6.76493,-0.931776)(6.78493,-0.72116)(6.80493,-0.479608)(6.82493,-0.208914)(6.84493,0.0889234)(6.86493,0.411707)(6.88493,0.757065)(6.90493,1.12246)(6.92493,1.50522)(6.94493,1.90253)(6.96493,2.3115)(6.98493,2.72914)(7.00493,3.1524)(7.02493,3.57821)(7.04493,4.00345)(7.06493,4.42504)(7.08493,4.83992)(7.10493,5.24508)(7.12493,5.63757)(7.14493,6.01455)(7.16493,6.3733)(7.18493,6.71122)(7.20493,7.02588)(7.22493,7.315)(7.24493,7.57651)(7.26493,7.80853)(7.28493,8.0094)(7.30493,8.17769)(7.32493,8.31221)(7.34493,8.412)(7.36493,8.47639)(7.38493,8.50493)(7.40493,8.49747)(7.42493,8.45408)(7.44493,8.37513)(7.46493,8.26123)(7.48493,8.11326)(7.50493,7.93233)(7.52493,7.71981)(7.54493,7.47728)(7.56493,7.20656)(7.58493,6.90967)(7.60493,6.58881)(7.62493,6.24636)(7.64493,5.88487)(7.66493,5.50702)(7.68493,5.1156)(7.70493,4.71352)(7.72493,4.30375)(7.74493,3.88932)(7.76493,3.47331)(7.78493,3.05879)(7.80493,2.64883)(7.82493,2.24646)(7.84493,1.85468)(7.86493,1.47638)(7.88493,1.11437)(7.90493,0.771338)(7.92493,0.449835)(7.94493,0.152256)(7.96493,-0.11918)(7.98493,-0.362443)(8.00493,-0.57571)(8.02493,-0.757373)(8.04493,-0.906058)(8.06493,-1.02063)(8.08493,-1.1002)(8.10493,-1.14413)(8.12493,-1.15206)(8.14493,-1.12388)(8.16493,-1.05973)(8.18493,-0.960047)(8.20493,-0.825495)(8.22493,-0.657009)(8.24493,-0.455768)(8.26493,-0.223193)(8.28493,0.0390674)(8.30493,0.329146)(8.32493,0.644971)(8.34493,0.984283)(8.36493,1.34465)(8.38493,1.72347)(8.40493,2.11804)(8.42493,2.5255)(8.44493,2.94291)(8.46493,3.36727)(8.48493,3.79549)(8.50493,4.22449)(8.52493,4.65116)(8.54493,5.07239)(8.56493,5.48514)(8.58493,5.88639)(8.60493,6.27324)(8.62493,6.64287)(8.64493,6.99257)(8.66493,7.31979)(8.68493,7.62214)(8.70493,7.8974)(8.72493,8.14354)(8.74493,8.35875)(8.76493,8.54144)(8.78493,8.69024)(8.80493,8.80403)(8.82493,8.88195)(8.84493,8.92339)(8.86493,8.928)(8.88493,8.8957)(8.90493,8.82668)(8.92493,8.72138)(8.94493,8.58052)(8.96493,8.40506)(8.98493,8.19622)(9.00493,7.95545)(9.02493,7.68446)(9.04493,7.38513)(9.06493,7.05959)(9.08493,6.71015)(9.10493,6.33926)(9.12493,5.94957)(9.14493,5.54385)(9.16493,5.12498)(9.18493,4.69594)(9.20493,4.2598)(9.22493,3.81965)(9.24493,3.37865)(9.26493,2.93993)(9.28493,2.50664)(9.30493,2.08186)(9.32493,1.66865)(9.34493,1.26995)(9.36493,0.888621)(9.38493,0.527389)(9.40493,0.188844)(9.42493,-0.124586)(9.44493,-0.410655)(9.46493,-0.667308)(9.48493,-0.892702)(9.50493,-1.08522)(9.52493,-1.24347)(9.54493,-1.36631)(9.56493,-1.45286)(9.58493,-1.50248)(9.60493,-1.51482)(9.62493,-1.48977)(9.64493,-1.4275)(9.66493,-1.32844)(9.68493,-1.1933)(9.70493,-1.02302)(9.72493,-0.818808)(9.74493,-0.58212)(9.76493,-0.314631)(9.78493,-0.01824)(9.80493,0.304947)(9.82493,0.652633)(9.84493,1.02235)(9.86493,1.41146)(9.88493,1.81722)(9.90493,2.23672)(9.92493,2.667)(9.94493,3.105)(9.96493,3.5476)(9.98493,3.99167)(10,4.32532) 
};

\addplot [
color=mycolor2,
solid
]
coordinates{
 (0,4)(0.00046341,4)(0.00301216,4.00016)(0.00764626,4.00086)(0.0258359,4.00979)(0.0449307,4.02744)(0.0649307,4.05447)(0.0849307,4.09077)(0.104931,4.13602)(0.124931,4.18968)(0.144931,4.2511)(0.164931,4.31956)(0.184931,4.39426)(0.204931,4.47428)(0.224931,4.55865)(0.244931,4.64634)(0.264931,4.73625)(0.284931,4.82725)(0.304931,4.91816)(0.324931,5.00778)(0.344931,5.0949)(0.364931,5.1783)(0.384931,5.25679)(0.404931,5.32918)(0.424931,5.39432)(0.444931,5.4511)(0.464931,5.49846)(0.484931,5.53543)(0.504931,5.56107)(0.524931,5.57457)(0.544931,5.57518)(0.564931,5.56227)(0.584931,5.5353)(0.604931,5.49386)(0.624931,5.43765)(0.644931,5.3665)(0.664931,5.28036)(0.684931,5.17932)(0.704931,5.0636)(0.724931,4.93354)(0.744931,4.78963)(0.764931,4.63247)(0.784931,4.4628)(0.804931,4.28147)(0.824931,4.08946)(0.844931,3.88785)(0.864931,3.67783)(0.884931,3.46068)(0.904931,3.23777)(0.924931,3.01056)(0.944931,2.78054)(0.964931,2.54928)(0.984931,2.3184)(1.00493,2.08954)(1.02493,1.86435)(1.04493,1.64449)(1.06493,1.43161)(1.08493,1.22736)(1.10493,1.03333)(1.12493,0.851058)(1.14493,0.682041)(1.16493,0.52769)(1.18493,0.389332)(1.20493,0.268195)(1.22493,0.165399)(1.24493,0.0819436)(1.26493,0.0186998)(1.28493,-0.0236)(1.30493,-0.0443691)(1.32493,-0.0431735)(1.34493,-0.0197373)(1.36493,0.0260526)(1.38493,0.0941431)(1.40493,0.184313)(1.42493,0.296172)(1.44493,0.429161)(1.46493,0.582555)(1.48493,0.755463)(1.50493,0.946837)(1.52493,1.15547)(1.54493,1.38002)(1.56493,1.619)(1.58493,1.87077)(1.60493,2.13361)(1.62493,2.40567)(1.64493,2.685)(1.66493,2.96957)(1.68493,3.2573)(1.70493,3.54602)(1.72493,3.83355)(1.74493,4.11767)(1.76493,4.39616)(1.78493,4.6668)(1.80493,4.92741)(1.82493,5.17583)(1.84493,5.40997)(1.86493,5.62781)(1.88493,5.82743)(1.90493,6.00699)(1.92493,6.1648)(1.94493,6.29926)(1.96493,6.40896)(1.98493,6.49263)(2.00493,6.54916)(2.02493,6.57764)(2.04493,6.57732)(2.06493,6.54768)(2.08493,6.48838)(2.10493,6.3993)(2.12493,6.28051)(2.14493,6.13232)(2.16493,5.95524)(2.18493,5.74999)(2.20493,5.51752)(2.22493,5.25895)(2.24493,4.97564)(2.26493,4.66911)(2.28493,4.34111)(2.30493,3.99351)(2.32493,3.62839)(2.34493,3.24796)(2.36493,2.85457)(2.38493,2.4507)(2.40493,2.03893)(2.42493,1.62194)(2.44493,1.20249)(2.46493,0.783371)(2.48493,0.36743)(2.50493,-0.0424736)(2.52493,-0.443482)(2.54493,-0.832759)(2.56493,-1.20751)(2.58493,-1.565)(2.60493,-1.90258)(2.62493,-2.21769)(2.64493,-2.5079)(2.66493,-2.77091)(2.68493,-3.00456)(2.70493,-3.20688)(2.72493,-3.37607)(2.74493,-3.51054)(2.76493,-3.6089)(2.78493,-3.66999)(2.80493,-3.69286)(2.82493,-3.67684)(2.84493,-3.62147)(2.86493,-3.52657)(2.88493,-3.3922)(2.90493,-3.21868)(2.92493,-3.00659)(2.94493,-2.75678)(2.96493,-2.47032)(2.98493,-2.14857)(3.00493,-1.7931)(3.02493,-1.40573)(3.04493,-0.988482)(3.06493,-0.543613)(3.08493,-0.0735652)(3.10493,0.419036)(3.12493,0.931399)(3.14493,1.46058)(3.16493,2.00351)(3.18493,2.55699)(3.20493,3.11777)(3.22493,3.68248)(3.24493,4.24774)(3.26493,4.81012)(3.28493,5.36621)(3.30493,5.91259)(3.32493,6.4459)(3.34493,6.96285)(3.36493,7.46021)(3.38493,7.93488)(3.40493,8.38388)(3.42493,8.80438)(3.44493,9.19372)(3.46493,9.54943)(3.48493,9.86924)(3.50493,10.1511)(3.52493,10.3932)(3.54493,10.594)(3.56493,10.7521)(3.58493,10.8665)(3.60493,10.9365)(3.62493,10.9616)(3.64493,10.9417)(3.66493,10.8767)(3.68493,10.7671)(3.70493,10.6136)(3.72493,10.4171)(3.74493,10.179)(3.76493,9.90069)(3.78493,9.584)(3.80493,9.231)(3.82493,8.84395)(3.84493,8.42537)(3.86493,7.978)(3.88493,7.50476)(3.90493,7.00874)(3.92493,6.4932)(3.94493,5.96154)(3.96493,5.41725)(3.98493,4.86393)(4.00493,4.30524)(4.02493,3.74489)(4.04493,3.18661)(4.06493,2.63411)(4.08493,2.09108)(4.10493,1.56116)(4.12493,1.04792)(4.14493,0.55479)(4.16493,0.0851169)(4.18493,-0.35792)(4.20493,-0.771306)(4.22493,-1.15222)(4.24493,-1.49803)(4.26493,-1.80636)(4.28493,-2.07506)(4.30493,-2.30223)(4.32493,-2.48625)(4.34493,-2.62578)(4.36493,-2.71976)(4.38493,-2.76745)(4.40493,-2.76839)(4.42493,-2.72244)(4.44493,-2.62977)(4.46493,-2.49085)(4.48493,-2.30646)(4.50493,-2.07769)(4.52493,-1.8059)(4.54493,-1.49277)(4.56493,-1.14022)(4.58493,-0.750464)(4.60493,-0.325945)(4.62493,0.130654)(4.64493,0.616432)(4.66493,1.12829)(4.68493,1.66297)(4.70493,2.21702)(4.72493,2.78689)(4.74493,3.3689)(4.76493,3.95929)(4.78493,4.55425)(4.80493,5.1499)(4.82493,5.74238)(4.84493,6.32784)(4.86493,6.90244)(4.88493,7.46245)(4.90493,8.00419)(4.92493,8.52411)(4.94493,9.01879)(4.96493,9.48497)(4.98493,9.91957)(5.00493,10.3197)(5.02493,10.6827)(5.04493,11.0062)(5.06493,11.288)(5.08493,11.5261)(5.10493,11.719)(5.12493,11.8653)(5.14493,11.964)(5.16493,12.0143)(5.18493,12.0158)(5.20493,11.9684)(5.22493,11.8723)(5.24493,11.728)(5.26493,11.5363)(5.28493,11.2983)(5.30493,11.0155)(5.32493,10.6895)(5.34493,10.3223)(5.36493,9.91616)(5.38493,9.47353)(5.40493,8.99713)(5.42493,8.48986)(5.44493,7.95485)(5.46493,7.39538)(5.48493,6.81486)(5.50493,6.21687)(5.52493,5.60506)(5.54493,4.98319)(5.56493,4.35503)(5.58493,3.72443)(5.60493,3.09521)(5.62493,2.4712)(5.64493,1.85616)(5.66493,1.25379)(5.68493,0.667713)(5.70493,0.101412)(5.72493,-0.441754)(5.74493,-0.958588)(5.76493,-1.44607)(5.78493,-1.90139)(5.80493,-2.32194)(5.82493,-2.70535)(5.84493,-3.04949)(5.86493,-3.35252)(5.88493,-3.61284)(5.90493,-3.82914)(5.92493,-4.0004)(5.94493,-4.12589)(5.96493,-4.20519)(5.98493,-4.23816)(6.00493,-4.22498)(6.02493,-4.16611)(6.04493,-4.0623)(6.06493,-3.91461)(6.08493,-3.72434)(6.10493,-3.49309)(6.12493,-3.22271)(6.14493,-2.91526)(6.16493,-2.57308)(6.18493,-2.19869)(6.20493,-1.7948)(6.22493,-1.3643)(6.24493,-0.910255)(6.26493,-0.435834)(6.28493,0.0556674)(6.30493,0.560864)(6.32493,1.07631)(6.34493,1.5985)(6.36493,2.12393)(6.38493,2.64909)(6.40493,3.17049)(6.42493,3.6847)(6.44493,4.18834)(6.46493,4.67816)(6.48493,5.151)(6.50493,5.60384)(6.52493,6.0338)(6.54493,6.43821)(6.56493,6.81455)(6.58493,7.16053)(6.60493,7.47407)(6.62493,7.75332)(6.64493,7.99668)(6.66493,8.20278)(6.68493,8.37054)(6.70493,8.49912)(6.72493,8.58797)(6.74493,8.63679)(6.76493,8.64556)(6.78493,8.61453)(6.80493,8.54421)(6.82493,8.43537)(6.84493,8.28904)(6.86493,8.1065)(6.88493,7.88924)(6.90493,7.639)(6.92493,7.35771)(6.94493,7.0475)(6.96493,6.71069)(6.98493,6.34974)(7.00493,5.96728)(7.02493,5.56606)(7.04493,5.14891)(7.06493,4.71878)(7.08493,4.27866)(7.10493,3.83162)(7.12493,3.38071)(7.14493,2.92902)(7.16493,2.4796)(7.18493,2.03548)(7.20493,1.59961)(7.22493,1.17489)(7.24493,0.764093)(7.26493,0.369906)(7.28493,-0.00512969)(7.30493,-0.358619)(7.32493,-0.688326)(7.34493,-0.992194)(7.36493,-1.26835)(7.38493,-1.51513)(7.40493,-1.73107)(7.42493,-1.91492)(7.44493,-2.06566)(7.46493,-2.1825)(7.48493,-2.26488)(7.50493,-2.31247)(7.52493,-2.32517)(7.54493,-2.30311)(7.56493,-2.24667)(7.58493,-2.15644)(7.60493,-2.0332)(7.62493,-1.87799)(7.64493,-1.692)(7.66493,-1.47666)(7.68493,-1.23354)(7.70493,-0.964386)(7.72493,-0.671107)(7.74493,-0.355738)(7.76493,-0.0204397)(7.78493,0.332522)(7.80493,0.700791)(7.82493,1.08194)(7.84493,1.47349)(7.86493,1.87292)(7.88493,2.2777)(7.90493,2.68528)(7.92493,3.09315)(7.94493,3.49882)(7.96493,3.89985)(7.98493,4.29388)(8.00493,4.6786)(8.02493,5.05184)(8.04493,5.41151)(8.06493,5.75565)(8.08493,6.08245)(8.10493,6.39024)(8.12493,6.67749)(8.14493,6.94285)(8.16493,7.18515)(8.18493,7.40339)(8.20493,7.59673)(8.22493,7.76456)(8.24493,7.90643)(8.26493,8.02207)(8.28493,8.11141)(8.30493,8.17457)(8.32493,8.21184)(8.34493,8.22368)(8.36493,8.21073)(8.38493,8.17379)(8.40493,8.11382)(8.42493,8.0319)(8.44493,7.92927)(8.46493,7.80727)(8.48493,7.66738)(8.50493,7.51114)(8.52493,7.3402)(8.54493,7.15627)(8.56493,6.96113)(8.58493,6.75657)(8.60493,6.54445)(8.62493,6.32662)(8.64493,6.10492)(8.66493,5.88119)(8.68493,5.65723)(8.70493,5.43482)(8.72493,5.21565)(8.74493,5.00137)(8.76493,4.79354)(8.78493,4.59362)(8.80493,4.403)(8.82493,4.22293)(8.84493,4.05455)(8.86493,3.89888)(8.88493,3.75682)(8.90493,3.62911)(8.92493,3.51637)(8.94493,3.41908)(8.96493,3.33755)(8.98493,3.27199)(9.00493,3.22244)(9.02493,3.18879)(9.04493,3.17083)(9.06493,3.16817)(9.08493,3.18034)(9.10493,3.20672)(9.12493,3.24657)(9.14493,3.29908)(9.16493,3.3633)(9.18493,3.43822)(9.20493,3.52275)(9.22493,3.61572)(9.24493,3.71591)(9.26493,3.82207)(9.28493,3.9329)(9.30493,4.04708)(9.32493,4.16329)(9.34493,4.2802)(9.36493,4.39652)(9.38493,4.51097)(9.40493,4.62229)(9.42493,4.7293)(9.44493,4.83086)(9.46493,4.9259)(9.48493,5.01343)(9.50493,5.09255)(9.52493,5.16243)(9.54493,5.22238)(9.56493,5.27177)(9.58493,5.31011)(9.60493,5.33702)(9.62493,5.35221)(9.64493,5.35554)(9.66493,5.34697)(9.68493,5.32659)(9.70493,5.29459)(9.72493,5.25128)(9.74493,5.1971)(9.76493,5.13256)(9.78493,5.0583)(9.80493,4.97504)(9.82493,4.88359)(9.84493,4.78485)(9.86493,4.67977)(9.88493,4.56938)(9.90493,4.45475)(9.92493,4.337)(9.94493,4.21728)(9.96493,4.09675)(9.98493,3.9766)(10,3.88703) 
};

\addplot [
color=mycolor3,
solid
]
coordinates{
 (0,5)(0.00046341,5)(0.00301216,4.99997)(0.00764626,4.99983)(0.0258359,4.99804)(0.0449307,4.99451)(0.0649307,4.98911)(0.0849307,4.98186)(0.104931,4.97282)(0.124931,4.96211)(0.144931,4.94986)(0.164931,4.93621)(0.184931,4.92133)(0.204931,4.9054)(0.224931,4.88861)(0.244931,4.87115)(0.264931,4.85325)(0.284931,4.83511)(0.304931,4.81696)(0.324931,4.79901)(0.344931,4.78148)(0.364931,4.76457)(0.384931,4.74849)(0.404931,4.73343)(0.424931,4.71955)(0.444931,4.70702)(0.464931,4.69598)(0.484931,4.68655)(0.504931,4.6788)(0.524931,4.67282)(0.544931,4.66865)(0.564931,4.66628)(0.584931,4.6657)(0.604931,4.66686)(0.624931,4.66967)(0.644931,4.67401)(0.664931,4.67973)(0.684931,4.68665)(0.704931,4.69456)(0.724931,4.7032)(0.744931,4.7123)(0.764931,4.72155)(0.784931,4.73062)(0.804931,4.73916)(0.824931,4.74678)(0.844931,4.75308)(0.864931,4.75765)(0.884931,4.76006)(0.904931,4.75987)(0.924931,4.75663)(0.944931,4.7499)(0.964931,4.73924)(0.984931,4.7242)(1.00493,4.70434)(1.02493,4.67926)(1.04493,4.64855)(1.06493,4.61183)(1.08493,4.56876)(1.10493,4.519)(1.12493,4.46228)(1.14493,4.39834)(1.16493,4.32697)(1.18493,4.248)(1.20493,4.16132)(1.22493,4.06686)(1.24493,3.9646)(1.26493,3.85459)(1.28493,3.7369)(1.30493,3.6117)(1.32493,3.47919)(1.34493,3.33964)(1.36493,3.19337)(1.38493,3.04078)(1.40493,2.8823)(1.42493,2.71842)(1.44493,2.54971)(1.46493,2.37676)(1.48493,2.20023)(1.50493,2.0208)(1.52493,1.83924)(1.54493,1.65631)(1.56493,1.47283)(1.58493,1.28964)(1.60493,1.10762)(1.62493,0.92765)(1.64493,0.750637)(1.66493,0.577494)(1.68493,0.409134)(1.70493,0.246471)(1.72493,0.0904044)(1.74493,-0.0581773)(1.76493,-0.19841)(1.78493,-0.329456)(1.80493,-0.450512)(1.82493,-0.560817)(1.84493,-0.659653)(1.86493,-0.746354)(1.88493,-0.820312)(1.90493,-0.880981)(1.92493,-0.927879)(1.94493,-0.960598)(1.96493,-0.978804)(1.98493,-0.982241)(2.00493,-0.970735)(2.02493,-0.944197)(2.04493,-0.902623)(2.06493,-0.846096)(2.08493,-0.774788)(2.10493,-0.688962)(2.12493,-0.588967)(2.14493,-0.47524)(2.16493,-0.348304)(2.18493,-0.208769)(2.20493,-0.0573224)(2.22493,0.105268)(2.24493,0.278159)(2.26493,0.46044)(2.28493,0.651134)(2.30493,0.849202)(2.32493,1.05356)(2.34493,1.26306)(2.36493,1.47652)(2.38493,1.69274)(2.40493,1.91048)(2.42493,2.12847)(2.44493,2.34545)(2.46493,2.56014)(2.48493,2.77128)(2.50493,2.97761)(2.52493,3.17791)(2.54493,3.37096)(2.56493,3.55562)(2.58493,3.73075)(2.60493,3.89529)(2.62493,4.04825)(2.64493,4.18868)(2.66493,4.31572)(2.68493,4.4286)(2.70493,4.52661)(2.72493,4.60917)(2.74493,4.67577)(2.76493,4.726)(2.78493,4.75958)(2.80493,4.7763)(2.82493,4.77609)(2.84493,4.75899)(2.86493,4.72513)(2.88493,4.67479)(2.90493,4.60831)(2.92493,4.52619)(2.94493,4.42902)(2.96493,4.31748)(2.98493,4.19237)(3.00493,4.05459)(3.02493,3.90513)(3.04493,3.74504)(3.06493,3.5755)(3.08493,3.39772)(3.10493,3.21301)(3.12493,3.02271)(3.14493,2.82823)(3.16493,2.63102)(3.18493,2.43256)(3.20493,2.23435)(3.22493,2.0379)(3.24493,1.84475)(3.26493,1.6564)(3.28493,1.47435)(3.30493,1.30008)(3.32493,1.13504)(3.34493,0.980611)(3.36493,0.838143)(3.38493,0.708913)(3.40493,0.594128)(3.42493,0.494916)(3.44493,0.412314)(3.46493,0.347267)(3.48493,0.300615)(3.50493,0.273091)(3.52493,0.265312)(3.54493,0.27778)(3.56493,0.310869)(3.58493,0.364829)(3.60493,0.439783)(3.62493,0.535718)(3.64493,0.652495)(3.66493,0.789839)(3.68493,0.947346)(3.70493,1.12448)(3.72493,1.32059)(3.74493,1.53487)(3.76493,1.76644)(3.78493,2.01426)(3.80493,2.27722)(3.82493,2.55407)(3.84493,2.8435)(3.86493,3.14409)(3.88493,3.45433)(3.90493,3.77267)(3.92493,4.09748)(3.94493,4.42708)(3.96493,4.75975)(3.98493,5.09372)(4.00493,5.42724)(4.02493,5.75851)(4.04493,6.08576)(4.06493,6.40721)(4.08493,6.72112)(4.10493,7.02577)(4.12493,7.3195)(4.14493,7.60069)(4.16493,7.86779)(4.18493,8.11933)(4.20493,8.35393)(4.22493,8.57029)(4.24493,8.76721)(4.26493,8.94361)(4.28493,9.09852)(4.30493,9.23109)(4.32493,9.34061)(4.34493,9.42648)(4.36493,9.48826)(4.38493,9.52564)(4.40493,9.53845)(4.42493,9.52665)(4.44493,9.49036)(4.46493,9.42985)(4.48493,9.34551)(4.50493,9.23787)(4.52493,9.10762)(4.54493,8.95554)(4.56493,8.78259)(4.58493,8.58979)(4.60493,8.37832)(4.62493,8.14946)(4.64493,7.90456)(4.66493,7.64509)(4.68493,7.3726)(4.70493,7.08869)(4.72493,6.79505)(4.74493,6.49339)(4.76493,6.18548)(4.78493,5.87312)(4.80493,5.55813)(4.82493,5.24232)(4.84493,4.9275)(4.86493,4.61549)(4.88493,4.30805)(4.90493,4.00691)(4.92493,3.71376)(4.94493,3.43022)(4.96493,3.15786)(4.98493,2.89814)(5.00493,2.65246)(5.02493,2.4221)(5.04493,2.20826)(5.06493,2.012)(5.08493,1.83428)(5.10493,1.67594)(5.12493,1.53767)(5.14493,1.42005)(5.16493,1.3235)(5.18493,1.24832)(5.20493,1.19466)(5.22493,1.16254)(5.24493,1.15183)(5.26493,1.16226)(5.28493,1.19343)(5.30493,1.24479)(5.32493,1.31569)(5.34493,1.40532)(5.36493,1.51278)(5.38493,1.63704)(5.40493,1.77696)(5.42493,1.93131)(5.44493,2.09878)(5.46493,2.27795)(5.48493,2.46735)(5.50493,2.66544)(5.52493,2.87063)(5.54493,3.08128)(5.56493,3.29574)(5.58493,3.51232)(5.60493,3.72932)(5.62493,3.94505)(5.64493,4.15784)(5.66493,4.36603)(5.68493,4.56799)(5.70493,4.76216)(5.72493,4.94702)(5.74493,5.1211)(5.76493,5.28303)(5.78493,5.43151)(5.80493,5.56533)(5.82493,5.68339)(5.84493,5.78468)(5.86493,5.86831)(5.88493,5.93351)(5.90493,5.97962)(5.92493,6.00613)(5.94493,6.01264)(5.96493,5.99888)(5.98493,5.96473)(6.00493,5.91018)(6.02493,5.83537)(6.04493,5.74056)(6.06493,5.62617)(6.08493,5.4927)(6.10493,5.34081)(6.12493,5.17126)(6.14493,4.98493)(6.16493,4.78281)(6.18493,4.56598)(6.20493,4.33563)(6.22493,4.09303)(6.24493,3.83952)(6.26493,3.57651)(6.28493,3.30549)(6.30493,3.02797)(6.32493,2.74551)(6.34493,2.45972)(6.36493,2.17221)(6.38493,1.8846)(6.40493,1.59852)(6.42493,1.31558)(6.44493,1.03738)(6.46493,0.765473)(6.48493,0.50139)(6.50493,0.246601)(6.52493,0.00252153)(6.54493,-0.229502)(6.56493,-0.448197)(6.58493,-0.652374)(6.60493,-0.840933)(6.62493,-1.01287)(6.64493,-1.16729)(6.66493,-1.30338)(6.68493,-1.42047)(6.70493,-1.51799)(6.72493,-1.59549)(6.74493,-1.65263)(6.76493,-1.6892)(6.78493,-1.70511)(6.80493,-1.70039)(6.82493,-1.6752)(6.84493,-1.62979)(6.86493,-1.56455)(6.88493,-1.47998)(6.90493,-1.37668)(6.92493,-1.25536)(6.94493,-1.11681)(6.96493,-0.961951)(6.98493,-0.791749)(7.00493,-0.60727)(7.02493,-0.409648)(7.04493,-0.200078)(7.06493,0.0201854)(7.08493,0.249843)(7.10493,0.487555)(7.12493,0.731954)(7.14493,0.98165)(7.16493,1.23524)(7.18493,1.49132)(7.20493,1.74849)(7.22493,2.00536)(7.24493,2.26057)(7.26493,2.5128)(7.28493,2.76074)(7.30493,3.00316)(7.32493,3.23886)(7.34493,3.4667)(7.36493,3.68564)(7.38493,3.89468)(7.40493,4.0929)(7.42493,4.27948)(7.44493,4.45368)(7.46493,4.61485)(7.48493,4.76243)(7.50493,4.89597)(7.52493,5.0151)(7.54493,5.11956)(7.56493,5.20917)(7.58493,5.28389)(7.60493,5.34373)(7.62493,5.38883)(7.64493,5.4194)(7.66493,5.43575)(7.68493,5.43829)(7.70493,5.42749)(7.72493,5.40391)(7.74493,5.36819)(7.76493,5.32102)(7.78493,5.26318)(7.80493,5.19548)(7.82493,5.1188)(7.84493,5.03404)(7.86493,4.94216)(7.88493,4.84415)(7.90493,4.74101)(7.92493,4.63376)(7.94493,4.52342)(7.96493,4.41104)(7.98493,4.29763)(8.00493,4.18421)(8.02493,4.07178)(8.04493,3.96131)(8.06493,3.85373)(8.08493,3.74996)(8.10493,3.65084)(8.12493,3.5572)(8.14493,3.4698)(8.16493,3.38934)(8.18493,3.31646)(8.20493,3.25175)(8.22493,3.19573)(8.24493,3.14884)(8.26493,3.11145)(8.28493,3.08389)(8.30493,3.06637)(8.32493,3.05906)(8.34493,3.06205)(8.36493,3.07536)(8.38493,3.09893)(8.40493,3.13265)(8.42493,3.17633)(8.44493,3.22971)(8.46493,3.29248)(8.48493,3.36429)(8.50493,3.44469)(8.52493,3.53322)(8.54493,3.62936)(8.56493,3.73255)(8.58493,3.8422)(8.60493,3.95768)(8.62493,4.07832)(8.64493,4.20346)(8.66493,4.33241)(8.68493,4.46445)(8.70493,4.59887)(8.72493,4.73497)(8.74493,4.87203)(8.76493,5.00935)(8.78493,5.14624)(8.80493,5.28202)(8.82493,5.41605)(8.84493,5.5477)(8.86493,5.67636)(8.88493,5.80147)(8.90493,5.92249)(8.92493,6.03893)(8.94493,6.15033)(8.96493,6.25627)(8.98493,6.35638)(9.00493,6.45032)(9.02493,6.53781)(9.04493,6.6186)(9.06493,6.6925)(9.08493,6.75936)(9.10493,6.81906)(9.12493,6.87154)(9.14493,6.91678)(9.16493,6.95478)(9.18493,6.9856)(9.20493,7.00934)(9.22493,7.0261)(9.24493,7.03604)(9.26493,7.03934)(9.28493,7.03621)(9.30493,7.02688)(9.32493,7.01159)(9.34493,6.9906)(9.36493,6.96419)(9.38493,6.93265)(9.40493,6.89626)(9.42493,6.85533)(9.44493,6.81014)(9.46493,6.761)(9.48493,6.70818)(9.50493,6.65198)(9.52493,6.59265)(9.54493,6.53046)(9.56493,6.46566)(9.58493,6.39846)(9.60493,6.32909)(9.62493,6.25773)(9.64493,6.18456)(9.66493,6.10973)(9.68493,6.03337)(9.70493,5.9556)(9.72493,5.8765)(9.74493,5.79615)(9.76493,5.7146)(9.78493,5.63189)(9.80493,5.54803)(9.82493,5.46302)(9.84493,5.37685)(9.86493,5.28949)(9.88493,5.20092)(9.90493,5.11108)(9.92493,5.01993)(9.94493,4.92741)(9.96493,4.83348)(9.98493,4.73806)(10,4.66517) 
};

\end{axis}
\end{tikzpicture}

%% file: Simulations/consensus_with_way_I_too_large.tikz
%
%
\begin{tikzpicture}

\definecolor{mycolor1}{rgb}{0,0.5,0}
\definecolor{mycolor2}{rgb}{0,0.75,0.75}
\definecolor{mycolor3}{rgb}{0.75,0,0.75}

\begin{axis}[%
scale only axis,
width=\fwidth,
height=\fheight,
xmin=0, xmax=10,
ymin=-150, ymax=100,
xlabel={$t$},
ylabel={$x(t)$},
axis on top]
\addplot [
color=blue,
solid
]
coordinates{
 (0,5)(0.000189186,5)(0.00122971,4.99996)(0.00312157,4.99974)(0.0118109,4.99621)(0.0272976,4.9797)(0.045041,4.94452)(0.065041,4.88394)(0.085041,4.80118)(0.105041,4.69644)(0.125041,4.57011)(0.145041,4.42273)(0.165041,4.25502)(0.185041,4.06788)(0.205041,3.86236)(0.225041,3.63969)(0.245041,3.40123)(0.265041,3.14849)(0.285041,2.88311)(0.305041,2.60685)(0.325041,2.32158)(0.345041,2.02928)(0.365041,1.73198)(0.385041,1.43182)(0.405041,1.13097)(0.425041,0.831636)(0.445041,0.536064)(0.465041,0.246496)(0.485041,-0.0348302)(0.505041,-0.305703)(0.525041,-0.563953)(0.545041,-0.807466)(0.565041,-1.0342)(0.585041,-1.24221)(0.605041,-1.42964)(0.625041,-1.59477)(0.645041,-1.73599)(0.665041,-1.85185)(0.685041,-1.94105)(0.705041,-2.00245)(0.725041,-2.03511)(0.745041,-2.03824)(0.765041,-2.01128)(0.785041,-1.95383)(0.805041,-1.86573)(0.825041,-1.74699)(0.845041,-1.59786)(0.865041,-1.41877)(0.885041,-1.21037)(0.905041,-0.973521)(0.925041,-0.70927)(0.945041,-0.418864)(0.965041,-0.103733)(0.985041,0.234511)(1.00504,0.59409)(1.02504,0.973069)(1.04504,1.36937)(1.06504,1.78077)(1.08504,2.20495)(1.10504,2.63949)(1.12504,3.08186)(1.14504,3.52948)(1.16504,3.97971)(1.18504,4.4299)(1.20504,4.87736)(1.22504,5.31942)(1.24504,5.75342)(1.26504,6.17675)(1.28504,6.58686)(1.30504,6.98127)(1.32504,7.35762)(1.34504,7.71364)(1.36504,8.04718)(1.38504,8.35626)(1.40504,8.63905)(1.42504,8.89388)(1.44504,9.11929)(1.46504,9.31399)(1.48504,9.47691)(1.50504,9.60718)(1.52504,9.70417)(1.54504,9.76746)(1.56504,9.79686)(1.58504,9.79241)(1.60504,9.75439)(1.62504,9.6833)(1.64504,9.57986)(1.66504,9.44502)(1.68504,9.27995)(1.70504,9.08601)(1.72504,8.86477)(1.74504,8.61799)(1.76504,8.34758)(1.78504,8.05565)(1.80504,7.7444)(1.82504,7.41621)(1.84504,7.07355)(1.86504,6.71898)(1.88504,6.35514)(1.90504,5.98474)(1.92504,5.6105)(1.94504,5.23519)(1.96504,4.86154)(1.98504,4.49229)(2.00504,4.13011)(2.02504,3.77763)(2.04504,3.43739)(2.06504,3.11182)(2.08504,2.80324)(2.10504,2.51383)(2.12504,2.24563)(2.14504,2.00049)(2.16504,1.78011)(2.18504,1.58598)(2.20504,1.41937)(2.22504,1.28137)(2.24504,1.17283)(2.26504,1.09437)(2.28504,1.04638)(2.30504,1.02902)(2.32504,1.04221)(2.34504,1.08564)(2.36504,1.15876)(2.38504,1.26079)(2.40504,1.39074)(2.42504,1.54739)(2.44504,1.72931)(2.46504,1.9349)(2.48504,2.16236)(2.50504,2.40971)(2.52504,2.67484)(2.54504,2.95549)(2.56504,3.24927)(2.58504,3.5537)(2.60504,3.86621)(2.62504,4.18416)(2.64504,4.50488)(2.66504,4.82567)(2.68504,5.14381)(2.70504,5.45662)(2.72504,5.76143)(2.74504,6.05566)(2.76504,6.33677)(2.78504,6.60235)(2.80504,6.85007)(2.82504,7.07776)(2.84504,7.28338)(2.86504,7.46505)(2.88504,7.62109)(2.90504,7.74997)(2.92504,7.85039)(2.94504,7.92126)(2.96504,7.96168)(2.98504,7.97101)(3.00504,7.9488)(3.02504,7.89487)(3.04504,7.80924)(3.06504,7.69218)(3.08504,7.54418)(3.10504,7.36595)(3.12504,7.15844)(3.14504,6.92278)(3.16504,6.66032)(3.18504,6.3726)(3.20504,6.06134)(3.22504,5.72841)(3.24504,5.37584)(3.26504,5.00579)(3.28504,4.62055)(3.30504,4.2225)(3.32504,3.81409)(3.34504,3.39786)(3.36504,2.97636)(3.38504,2.55218)(3.40504,2.12793)(3.42504,1.70617)(3.44504,1.28945)(3.46504,0.880268)(3.48504,0.481028)(3.50504,0.0940642)(3.52504,-0.278404)(3.54504,-0.634276)(3.56504,-0.971591)(3.58504,-1.28854)(3.60504,-1.58346)(3.62504,-1.8549)(3.64504,-2.10154)(3.66504,-2.3223)(3.68504,-2.51625)(3.70504,-2.68271)(3.72504,-2.82116)(3.74504,-2.93131)(3.76504,-3.01306)(3.78504,-3.06654)(3.80504,-3.09203)(3.82504,-3.09005)(3.84504,-3.06128)(3.86504,-3.00658)(3.88504,-2.92698)(3.90504,-2.82366)(3.92504,-2.69796)(3.94504,-2.55132)(3.96504,-2.38533)(3.98504,-2.20166)(4.00504,-2.00206)(4.02504,-1.78835)(4.04504,-1.56242)(4.06504,-1.32616)(4.08504,-1.08151)(4.10504,-0.830387)(4.12504,-0.574698)(4.14504,-0.316315)(4.16504,-0.0570606)(4.18504,0.201309)(4.20504,0.457115)(4.22504,0.708775)(4.24504,0.954809)(4.26504,1.19386)(4.28504,1.42468)(4.30504,1.64618)(4.32504,1.8574)(4.34504,2.05753)(4.36504,2.24592)(4.38504,2.42206)(4.40504,2.5856)(4.42504,2.73635)(4.44504,2.87428)(4.46504,2.9995)(4.48504,3.11226)(4.50504,3.21295)(4.52504,3.3021)(4.54504,3.38034)(4.56504,3.44842)(4.58504,3.5072)(4.60504,3.5576)(4.62504,3.60063)(4.64504,3.63735)(4.66504,3.66886)(4.68504,3.69629)(4.70504,3.72078)(4.72504,3.74348)(4.74504,3.7655)(4.76504,3.78792)(4.78504,3.8118)(4.80504,3.83809)(4.82504,3.86769)(4.84504,3.90143)(4.86504,3.93999)(4.88504,3.98399)(4.90504,4.0339)(4.92504,4.09007)(4.94504,4.15272)(4.96504,4.22192)(4.98504,4.29761)(5.00504,4.37958)(5.02504,4.46748)(5.04504,4.56082)(5.06504,4.65898)(5.08504,4.76118)(5.10504,4.86654)(5.12504,4.97407)(5.14504,5.08265)(5.16504,5.19108)(5.18504,5.29808)(5.20504,5.40231)(5.22504,5.50235)(5.24504,5.59679)(5.26504,5.68417)(5.28504,5.76305)(5.30504,5.83199)(5.32504,5.88962)(5.34504,5.93459)(5.36504,5.96568)(5.38504,5.9817)(5.40504,5.98164)(5.42504,5.96457)(5.44504,5.92974)(5.46504,5.87655)(5.48504,5.80459)(5.50504,5.71363)(5.52504,5.60365)(5.54504,5.47484)(5.56504,5.3276)(5.58504,5.16257)(5.60504,4.98062)(5.62504,4.78284)(5.64504,4.57055)(5.66504,4.3453)(5.68504,4.10884)(5.70504,3.86317)(5.72504,3.61044)(5.74504,3.35303)(5.76504,3.09346)(5.78504,2.83443)(5.80504,2.57874)(5.82504,2.32932)(5.84504,2.0892)(5.86504,1.86143)(5.88504,1.64915)(5.90504,1.45545)(5.92504,1.28344)(5.94504,1.13614)(5.96504,1.01652)(5.98504,0.927411)(6.00504,0.871502)(6.02504,0.851306)(6.04504,0.869127)(6.06504,0.92703)(6.08504,1.02681)(6.10504,1.16999)(6.12504,1.35773)(6.14504,1.5909)(6.16504,1.86997)(6.18504,2.19505)(6.20504,2.56582)(6.22504,2.9816)(6.24504,3.44125)(6.26504,3.94323)(6.28504,4.48553)(6.30504,5.06576)(6.32504,5.68109)(6.34504,6.32825)(6.36504,7.0036)(6.38504,7.70308)(6.40504,8.42231)(6.42504,9.1565)(6.44504,9.9006)(6.46504,10.6492)(6.48504,11.3968)(6.50504,12.1375)(6.52504,12.8653)(6.54504,13.574)(6.56504,14.2574)(6.58504,14.9094)(6.60504,15.5235)(6.62504,16.0936)(6.64504,16.6137)(6.66504,17.0777)(6.68504,17.48)(6.70504,17.8152)(6.72504,18.078)(6.74504,18.2638)(6.76504,18.3682)(6.78504,18.3873)(6.80504,18.3177)(6.82504,18.1567)(6.84504,17.902)(6.86504,17.5521)(6.88504,17.1058)(6.90504,16.563)(6.92504,15.9242)(6.94504,15.1904)(6.96504,14.3636)(6.98504,13.4463)(7.00504,12.4419)(7.02504,11.3546)(7.04504,10.1891)(7.06504,8.95099)(7.08504,7.64641)(7.10504,6.28225)(7.12504,4.86599)(7.14504,3.40566)(7.16504,1.90987)(7.18504,0.387684)(7.20504,-1.15141)(7.22504,-2.69754)(7.24504,-4.24055)(7.26504,-5.77003)(7.28504,-7.27545)(7.30504,-8.74614)(7.32504,-10.1715)(7.34504,-11.5409)(7.36504,-12.8439)(7.38504,-14.0703)(7.40504,-15.2103)(7.42504,-16.2544)(7.44504,-17.1935)(7.46504,-18.0193)(7.48504,-18.7238)(7.50504,-19.3)(7.52504,-19.7414)(7.54504,-20.0426)(7.56504,-20.199)(7.58504,-20.2067)(7.60504,-20.0632)(7.62504,-19.7668)(7.64504,-19.3168)(7.66504,-18.7137)(7.68504,-17.9593)(7.70504,-17.0561)(7.72504,-16.0081)(7.74504,-14.8202)(7.76504,-13.4985)(7.78504,-12.0502)(7.80504,-10.4836)(7.82504,-8.80781)(7.84504,-7.03307)(7.86504,-5.1705)(7.88504,-3.23205)(7.90504,-1.23043)(7.92504,0.82094)(7.94504,2.90805)(7.96504,5.01638)(7.98504,7.13097)(8.00504,9.23658)(8.02504,11.3178)(8.04504,13.3589)(8.06504,15.3446)(8.08504,17.2594)(8.10504,19.0881)(8.12504,20.816)(8.14504,22.4287)(8.16504,23.9127)(8.18504,25.2549)(8.20504,26.4432)(8.22504,27.4663)(8.24504,28.3141)(8.26504,28.9774)(8.28504,29.4485)(8.30504,29.7206)(8.32504,29.7885)(8.34504,29.6485)(8.36504,29.2982)(8.38504,28.7368)(8.40504,27.9649)(8.42504,26.9849)(8.44504,25.8006)(8.46504,24.4174)(8.48504,22.8424)(8.50504,21.0841)(8.52504,19.1525)(8.54504,17.0591)(8.56504,14.8167)(8.58504,12.4396)(8.60504,9.94316)(8.62504,7.34395)(8.64504,4.65955)(8.66504,1.90848)(8.68504,-0.889937)(8.70504,-3.71572)(8.72504,-6.54836)(8.74504,-9.36693)(8.76504,-12.1503)(8.78504,-14.8772)(8.80504,-17.5264)(8.82504,-20.0771)(8.84504,-22.5086)(8.86504,-24.8009)(8.88504,-26.9346)(8.90504,-28.8913)(8.92504,-30.6534)(8.94504,-32.2046)(8.96504,-33.5299)(8.98504,-34.6155)(9.00504,-35.4492)(9.02504,-36.0207)(9.04504,-36.3212)(9.06504,-36.3437)(9.08504,-36.0832)(9.10504,-35.5367)(9.12504,-34.7033)(9.14504,-33.5839)(9.16504,-32.1819)(9.18504,-30.5023)(9.20504,-28.5527)(9.22504,-26.3423)(9.24504,-23.8827)(9.26504,-21.1872)(9.28504,-18.2712)(9.30504,-15.1518)(9.32504,-11.8478)(9.34504,-8.37971)(9.36504,-4.76939)(9.38504,-1.04011)(9.40504,2.78371)(9.42504,6.67662)(9.44504,10.6123)(9.46504,14.5639)(9.48504,18.5038)(9.50504,22.4044)(9.52504,26.2379)(9.54504,29.9766)(9.56504,33.593)(9.58504,37.0604)(9.60504,40.3527)(9.62504,43.4447)(9.64504,46.3122)(9.66504,48.9325)(9.68504,51.2844)(9.70504,53.3483)(9.72504,55.1063)(9.74504,56.5427)(9.76504,57.6439)(9.78504,58.3984)(9.80504,58.7972)(9.82504,58.8336)(9.84504,58.5036)(9.86504,57.8056)(9.88504,56.7408)(9.90504,55.313)(9.92504,53.5286)(9.94504,51.3966)(9.96504,48.9287)(9.98504,46.1391)(10,43.8523) 
};

\addplot [
color=mycolor1,
solid
]
coordinates{
 (0,-6)(0.000189186,-6)(0.00122971,-5.9999)(0.00312157,-5.99939)(0.0118109,-5.99124)(0.0272976,-5.95304)(0.045041,-5.8717)(0.065041,-5.73166)(0.085041,-5.54042)(0.105041,-5.29844)(0.125041,-5.00654)(0.145041,-4.66597)(0.165041,-4.27831)(0.185041,-3.84555)(0.205041,-3.37001)(0.225041,-2.85437)(0.245041,-2.30164)(0.265041,-1.71513)(0.285041,-1.09843)(0.305041,-0.455397)(0.325041,0.209877)(0.345041,0.893099)(0.365041,1.5898)(0.385041,2.29536)(0.405041,3.00506)(0.425041,3.7141)(0.445041,4.41764)(0.465041,5.11084)(0.485041,5.78886)(0.505041,6.44697)(0.525041,7.08052)(0.545041,7.68499)(0.565041,8.25606)(0.585041,8.78957)(0.605041,9.28165)(0.625041,9.72867)(0.645041,10.1273)(0.665041,10.4745)(0.685041,10.7678)(0.705041,11.0046)(0.725041,11.1834)(0.745041,11.3024)(0.765041,11.3607)(0.785041,11.3577)(0.805041,11.2932)(0.825041,11.1675)(0.845041,10.9814)(0.865041,10.7359)(0.885041,10.4328)(0.905041,10.0741)(0.925041,9.66224)(0.945041,9.20016)(0.965041,8.69114)(0.985041,8.13881)(1.00504,7.54718)(1.02504,6.92054)(1.04504,6.2635)(1.06504,5.5809)(1.08504,4.87781)(1.10504,4.15949)(1.12504,3.43135)(1.14504,2.69888)(1.16504,1.96766)(1.18504,1.24331)(1.20504,0.531389)(1.22504,-0.162558)(1.24504,-0.833104)(1.26504,-1.47496)(1.28504,-2.08302)(1.30504,-2.65239)(1.32504,-3.17844)(1.34504,-3.65684)(1.36504,-4.08357)(1.38504,-4.45499)(1.40504,-4.76784)(1.42504,-5.01928)(1.44504,-5.20691)(1.46504,-5.32879)(1.48504,-5.38346)(1.50504,-5.36996)(1.52504,-5.28781)(1.54504,-5.13708)(1.56504,-4.91831)(1.58504,-4.63257)(1.60504,-4.28144)(1.62504,-3.86699)(1.64504,-3.39178)(1.66504,-2.85883)(1.68504,-2.27161)(1.70504,-1.63402)(1.72504,-0.950357)(1.74504,-0.225278)(1.76504,0.53622)(1.78504,1.32884)(1.80504,2.14703)(1.82504,2.98502)(1.84504,3.83685)(1.86504,4.69644)(1.88504,5.55761)(1.90504,6.41414)(1.92504,7.25982)(1.94504,8.08849)(1.96504,8.89407)(1.98504,9.67065)(2.00504,10.4125)(2.02504,11.1141)(2.04504,11.7702)(2.06504,12.3759)(2.08504,12.9266)(2.10504,13.4183)(2.12504,13.847)(2.14504,14.2096)(2.16504,14.5031)(2.18504,14.7254)(2.20504,14.8745)(2.22504,14.9492)(2.24504,14.9488)(2.26504,14.8732)(2.28504,14.7226)(2.30504,14.4981)(2.32504,14.2012)(2.34504,13.8338)(2.36504,13.3986)(2.38504,12.8986)(2.40504,12.3373)(2.42504,11.7188)(2.44504,11.0475)(2.46504,10.3283)(2.48504,9.5662)(2.50504,8.76687)(2.52504,7.93605)(2.54504,7.07977)(2.56504,6.20425)(2.58504,5.31583)(2.60504,4.42095)(2.62504,3.52612)(2.64504,2.63781)(2.66504,1.76245)(2.68504,0.906343)(2.70504,0.0756635)(2.72504,-0.723633)(2.74504,-1.48584)(2.76504,-2.20554)(2.78504,-2.87766)(2.80504,-3.49749)(2.82504,-4.06075)(2.84504,-4.56359)(2.86504,-5.00265)(2.88504,-5.37506)(2.90504,-5.67848)(2.92504,-5.91112)(2.94504,-6.07174)(2.96504,-6.15965)(2.98504,-6.17477)(3.00504,-6.11754)(3.02504,-5.98902)(3.04504,-5.7908)(3.06504,-5.525)(3.08504,-5.19431)(3.10504,-4.80188)(3.12504,-4.35137)(3.14504,-3.84688)(3.16504,-3.29292)(3.18504,-2.69438)(3.20504,-2.0565)(3.22504,-1.38481)(3.24504,-0.685083)(3.26504,0.0367014)(3.28504,0.774408)(3.30504,1.5218)(3.32504,2.27258)(3.34504,3.02047)(3.36504,3.75922)(3.38504,4.4827)(3.40504,5.18492)(3.42504,5.86011)(3.44504,6.50273)(3.46504,7.10755)(3.48504,7.66967)(3.50504,8.18457)(3.52504,8.64816)(3.54504,9.05677)(3.56504,9.40725)(3.58504,9.69691)(3.60504,9.92363)(3.62504,10.0858)(3.64504,10.1824)(3.66504,10.2129)(3.68504,10.1775)(3.70504,10.0768)(3.72504,9.912)(3.74504,9.68488)(3.76504,9.39775)(3.78504,9.05342)(3.80504,8.65519)(3.82504,8.2068)(3.84504,7.71244)(3.86504,7.17666)(3.88504,6.60439)(3.90504,6.00083)(3.92504,5.37146)(3.94504,4.72194)(3.96504,4.05812)(3.98504,3.38594)(4.00504,2.7114)(4.02504,2.04049)(4.04504,1.37915)(4.06504,0.733225)(4.08504,0.108389)(4.10504,-0.489888)(4.12504,-1.05639)(4.14504,-1.58622)(4.16504,-2.07482)(4.18504,-2.51801)(4.20504,-2.91207)(4.22504,-3.25371)(4.24504,-3.54016)(4.26504,-3.76914)(4.28504,-3.93894)(4.30504,-4.04838)(4.32504,-4.09685)(4.34504,-4.08433)(4.36504,-4.01135)(4.38504,-3.87902)(4.40504,-3.68902)(4.42504,-3.44355)(4.44504,-3.14536)(4.46504,-2.79767)(4.48504,-2.40421)(4.50504,-1.9691)(4.52504,-1.49689)(4.54504,-0.99248)(4.56504,-0.461069)(4.58504,0.0918763)(4.60504,0.660681)(4.62504,1.23952)(4.64504,1.82245)(4.66504,2.40352)(4.68504,2.97677)(4.70504,3.53632)(4.72504,4.07641)(4.74504,4.5915)(4.76504,5.07627)(4.78504,5.52571)(4.80504,5.93515)(4.82504,6.30033)(4.84504,6.61743)(4.86504,6.88313)(4.88504,7.09461)(4.90504,7.24964)(4.92504,7.34654)(4.94504,7.38427)(4.96504,7.36239)(4.98504,7.28112)(5.00504,7.14129)(5.02504,6.94441)(5.04504,6.6926)(5.06504,6.38859)(5.08504,6.03575)(5.10504,5.638)(5.12504,5.19981)(5.14504,4.72616)(5.16504,4.2225)(5.18504,3.69471)(5.20504,3.14902)(5.22504,2.59199)(5.24504,2.03041)(5.26504,1.47129)(5.28504,0.921728)(5.30504,0.38891)(5.32504,-0.120015)(5.34504,-0.597978)(5.36504,-1.03807)(5.38504,-1.43362)(5.40504,-1.77823)(5.42504,-2.06592)(5.44504,-2.29109)(5.46504,-2.44869)(5.48504,-2.5342)(5.50504,-2.54372)(5.52504,-2.47404)(5.54504,-2.32265)(5.56504,-2.08779)(5.58504,-1.76849)(5.60504,-1.36463)(5.62504,-0.876879)(5.64504,-0.306802)(5.66504,0.343201)(5.68504,1.06986)(5.70504,1.86908)(5.72504,2.73591)(5.74504,3.66461)(5.76504,4.64867)(5.78504,5.68085)(5.80504,6.75323)(5.82504,7.85726)(5.84504,8.98386)(5.86504,10.1234)(5.88504,11.2659)(5.90504,12.4011)(5.92504,13.5182)(5.94504,14.6066)(5.96504,15.6555)(5.98504,16.654)(6.00504,17.5916)(6.02504,18.4577)(6.04504,19.2423)(6.06504,19.9356)(6.08504,20.5286)(6.10504,21.0127)(6.12504,21.3799)(6.14504,21.6234)(6.16504,21.7368)(6.18504,21.715)(6.20504,21.5538)(6.22504,21.25)(6.24504,20.8017)(6.26504,20.2081)(6.28504,19.4696)(6.30504,18.5878)(6.32504,17.5658)(6.34504,16.4077)(6.36504,15.1189)(6.38504,13.7063)(6.40504,12.1778)(6.42504,10.5425)(6.44504,8.81061)(6.46504,6.99362)(6.48504,5.10386)(6.50504,3.15462)(6.52504,1.16006)(6.54504,-0.864923)(6.56504,-2.90481)(6.58504,-4.94355)(6.60504,-6.96467)(6.62504,-8.95146)(6.64504,-10.887)(6.66504,-12.7545)(6.68504,-14.5372)(6.70504,-16.2185)(6.72504,-17.7825)(6.74504,-19.2137)(6.76504,-20.4973)(6.78504,-21.6192)(6.80504,-22.5667)(6.82504,-23.3277)(6.84504,-23.8916)(6.86504,-24.2491)(6.88504,-24.3924)(6.90504,-24.3151)(6.92504,-24.0126)(6.94504,-23.4818)(6.96504,-22.7214)(6.98504,-21.7322)(7.00504,-20.5166)(7.02504,-19.0787)(7.04504,-17.425)(7.06504,-15.5633)(7.08504,-13.5036)(7.10504,-11.2577)(7.12504,-8.83905)(7.14504,-6.26273)(7.16504,-3.54553)(7.18504,-0.705713)(7.20504,2.23709)(7.22504,5.26198)(7.24504,8.34694)(7.26504,11.469)(7.28504,14.6043)(7.30504,17.7286)(7.32504,20.8168)(7.34504,23.844)(7.36504,26.7848)(7.38504,29.6142)(7.40504,32.3075)(7.42504,34.8405)(7.44504,37.1896)(7.46504,39.3323)(7.48504,41.2472)(7.50504,42.9142)(7.52504,44.3148)(7.54504,45.432)(7.56504,46.2507)(7.58504,46.758)(7.60504,46.9428)(7.62504,46.7967)(7.64504,46.3134)(7.66504,45.489)(7.68504,44.3225)(7.70504,42.8152)(7.72504,40.9713)(7.74504,38.7975)(7.76504,36.3032)(7.78504,33.5006)(7.80504,30.4045)(7.82504,27.032)(7.84504,23.403)(7.86504,19.5395)(7.88504,15.4659)(7.90504,11.2086)(7.92504,6.79594)(7.94504,2.25801)(7.96504,-2.37358)(7.98504,-7.06587)(8.00504,-11.7848)(8.02504,-16.4956)(8.04504,-21.1626)(8.06504,-25.7502)(8.08504,-30.2223)(8.10504,-34.5434)(8.12504,-38.6782)(8.14504,-42.5924)(8.16504,-46.2524)(8.18504,-49.6264)(8.20504,-52.6837)(8.22504,-55.3957)(8.24504,-57.736)(8.26504,-59.6802)(8.28504,-61.2068)(8.30504,-62.2969)(8.32504,-62.9346)(8.34504,-63.1071)(8.36504,-62.8049)(8.38504,-62.0219)(8.40504,-60.7555)(8.42504,-59.0068)(8.44504,-56.7804)(8.46504,-54.0847)(8.48504,-50.9317)(8.50504,-47.3369)(8.52504,-43.3195)(8.54504,-38.9023)(8.56504,-34.1114)(8.58504,-28.9758)(8.60504,-23.5281)(8.62504,-17.8034)(8.64504,-11.8394)(8.66504,-5.67641)(8.68504,0.643281)(8.70504,7.07553)(8.72504,13.5747)(8.74504,20.0938)(8.76504,26.5852)(8.78504,33.0006)(8.80504,39.2916)(8.82504,45.4099)(8.84504,51.308)(8.86504,56.9391)(8.88504,62.2578)(8.90504,67.2205)(8.92504,71.7854)(8.94504,75.9131)(8.96504,79.5671)(8.98504,82.7138)(9.00504,85.3229)(9.02504,87.3677)(9.04504,88.8253)(9.06504,89.6769)(9.08504,89.9081)(9.10504,89.5089)(9.12504,88.4736)(9.14504,86.8015)(9.16504,84.4966)(9.18504,81.5677)(9.20504,78.0281)(9.22504,73.8964)(9.24504,69.1953)(9.26504,63.9526)(9.28504,58.2001)(9.30504,51.974)(9.32504,45.3146)(9.34504,38.266)(9.36504,30.8755)(9.38504,23.194)(9.40504,15.2749)(9.42504,7.17434)(9.44504,-1.04971)(9.46504,-9.33743)(9.48504,-17.6278)(9.50504,-25.859)(9.52504,-33.969)(9.54504,-41.8955)(9.56504,-49.5773)(9.58504,-56.9538)(9.60504,-63.9663)(9.62504,-70.5579)(9.64504,-76.6741)(9.66504,-82.2633)(9.68504,-87.2774)(9.70504,-91.6718)(9.72504,-95.4061)(9.74504,-98.4443)(9.76504,-100.755)(9.78504,-102.312)(9.80504,-103.095)(9.82504,-103.089)(9.84504,-102.283)(9.86504,-100.674)(9.88504,-98.2652)(9.90504,-95.0636)(9.92504,-91.0838)(9.94504,-86.3463)(9.96504,-80.8771)(9.98504,-74.7083)(10,-69.6591) 
};

\addplot [
color=red,
solid
]
coordinates{
 (0,8)(0.000189186,8)(0.00122971,7.99993)(0.00312157,7.99956)(0.0118109,7.99369)(0.0272976,7.9662)(0.045041,7.90767)(0.065041,7.80698)(0.085041,7.66963)(0.105041,7.49609)(0.125041,7.28716)(0.145041,7.04392)(0.165041,6.76779)(0.185041,6.46048)(0.205041,6.124)(0.225041,5.76061)(0.245041,5.37285)(0.265041,4.9635)(0.285041,4.53555)(0.305041,4.09219)(0.325041,3.6368)(0.345041,3.17289)(0.365041,2.7041)(0.385041,2.23415)(0.405041,1.76684)(0.425041,1.30599)(0.445041,0.855411)(0.465041,0.418903)(0.485041,0.000191895)(0.505041,-0.397089)(0.525041,-0.769433)(0.545041,-1.11349)(0.565041,-1.42609)(0.585041,-1.70428)(0.605041,-1.94535)(0.625041,-2.14686)(0.645041,-2.30665)(0.665041,-2.42286)(0.685041,-2.494)(0.705041,-2.51888)(0.725041,-2.49671)(0.745041,-2.42703)(0.765041,-2.3098)(0.785041,-2.14533)(0.805041,-1.93435)(0.825041,-1.67796)(0.845041,-1.37764)(0.865041,-1.03525)(0.885041,-0.653004)(0.905041,-0.233488)(0.925041,0.220392)(0.945041,0.70541)(0.965041,1.21805)(0.985041,1.75452)(1.00504,2.3108)(1.02504,2.88264)(1.04504,3.46564)(1.06504,4.05521)(1.08504,4.64669)(1.10504,5.23533)(1.12504,5.81633)(1.14504,6.38492)(1.16504,6.93634)(1.18504,7.46592)(1.20504,7.96913)(1.22504,8.44156)(1.24504,8.87901)(1.26504,9.27752)(1.28504,9.63338)(1.30504,9.94318)(1.32504,10.2038)(1.34504,10.4126)(1.36504,10.5672)(1.38504,10.6657)(1.40504,10.7066)(1.42504,10.6888)(1.44504,10.6117)(1.46504,10.4752)(1.48504,10.2796)(1.50504,10.0259)(1.52504,9.71526)(1.54504,9.34949)(1.56504,8.93086)(1.58504,8.46207)(1.60504,7.94625)(1.62504,7.38697)(1.64504,6.78816)(1.66504,6.15415)(1.68504,5.48959)(1.70504,4.79942)(1.72504,4.08885)(1.74504,3.36334)(1.76504,2.62851)(1.78504,1.89012)(1.80504,1.15405)(1.82504,0.426203)(1.84504,-0.287492)(1.86504,-0.981161)(1.88504,-1.64902)(1.90504,-2.28542)(1.92504,-2.88491)(1.94504,-3.44226)(1.96504,-3.95252)(1.98504,-4.41106)(2.00504,-4.81363)(2.02504,-5.15635)(2.04504,-5.43579)(2.06504,-5.64901)(2.08504,-5.79353)(2.10504,-5.86742)(2.12504,-5.86929)(2.14504,-5.79832)(2.16504,-5.65424)(2.18504,-5.43739)(2.20504,-5.1487)(2.22504,-4.78967)(2.24504,-4.3624)(2.26504,-3.86956)(2.28504,-3.31436)(2.30504,-2.70057)(2.32504,-2.03247)(2.34504,-1.31481)(2.36504,-0.5528)(2.38504,0.247924)(2.40504,1.08136)(2.42504,1.94119)(2.44504,2.82078)(2.46504,3.71329)(2.48504,4.61171)(2.50504,5.50886)(2.52504,6.39755)(2.54504,7.27052)(2.56504,8.12061)(2.58504,8.94073)(2.60504,9.72398)(2.62504,10.4637)(2.64504,11.1533)(2.66504,11.787)(2.68504,12.3589)(2.70504,12.8638)(2.72504,13.2969)(2.74504,13.654)(2.76504,13.9315)(2.78504,14.1262)(2.80504,14.2357)(2.82504,14.2583)(2.84504,14.1929)(2.86504,14.0391)(2.88504,13.7973)(2.90504,13.4685)(2.92504,13.0545)(2.94504,12.5577)(2.96504,11.9814)(2.98504,11.3295)(3.00504,10.6063)(3.02504,9.81718)(3.04504,8.96777)(3.06504,8.06439)(3.08504,7.11385)(3.10504,6.12341)(3.12504,5.10074)(3.14504,4.05386)(3.16504,2.99107)(3.18504,1.92089)(3.20504,0.851995)(3.22504,-0.206846)(3.24504,-1.24686)(3.26504,-2.25932)(3.28504,-3.23563)(3.30504,-4.1674)(3.32504,-5.0465)(3.34504,-5.86514)(3.36504,-6.61595)(3.38504,-7.29201)(3.40504,-7.88695)(3.42504,-8.39499)(3.44504,-8.81097)(3.46504,-9.13046)(3.48504,-9.34975)(3.50504,-9.46591)(3.52504,-9.47681)(3.54504,-9.38117)(3.56504,-9.17856)(3.58504,-8.86943)(3.60504,-8.4551)(3.62504,-7.93776)(3.64504,-7.32047)(3.66504,-6.60715)(3.68504,-5.80255)(3.70504,-4.91222)(3.72504,-3.94248)(3.74504,-2.90037)(3.76504,-1.79361)(3.78504,-0.630534)(3.80504,0.579966)(3.82504,1.8285)(3.84504,3.10526)(3.86504,4.40008)(3.88504,5.70253)(3.90504,7.00202)(3.92504,8.28781)(3.94504,9.54918)(3.96504,10.7755)(3.98504,11.9562)(4.00504,13.0811)(4.02504,14.1402)(4.04504,15.124)(4.06504,16.0236)(4.08504,16.8306)(4.10504,17.5371)(4.12504,18.1361)(4.14504,18.6214)(4.16504,18.9875)(4.18504,19.2301)(4.20504,19.3456)(4.22504,19.3314)(4.24504,19.186)(4.26504,18.9091)(4.28504,18.5012)(4.30504,17.964)(4.32504,17.3002)(4.34504,16.5139)(4.36504,15.6097)(4.38504,14.5936)(4.40504,13.4726)(4.42504,12.2545)(4.44504,10.948)(4.46504,9.56259)(4.48504,8.10869)(4.50504,6.59726)(4.52504,5.03993)(4.54504,3.44882)(4.56504,1.83649)(4.58504,0.215826)(4.60504,-1.40007)(4.62504,-2.99794)(4.64504,-4.56455)(4.66504,-6.08671)(4.68504,-7.55147)(4.70504,-8.94615)(4.72504,-10.2585)(4.74504,-11.4768)(4.76504,-12.5899)(4.78504,-13.5874)(4.80504,-14.4597)(4.82504,-15.1982)(4.84504,-15.795)(4.86504,-16.2435)(4.88504,-16.5383)(4.90504,-16.6748)(4.92504,-16.65)(4.94504,-16.462)(4.96504,-16.1102)(4.98504,-15.5953)(5.00504,-14.9195)(5.02504,-14.0861)(5.04504,-13.0999)(5.06504,-11.9669)(5.08504,-10.6943)(5.10504,-9.29068)(5.12504,-7.76566)(5.14504,-6.13003)(5.16504,-4.39557)(5.18504,-2.57501)(5.20504,-0.681937)(5.22504,1.26931)(5.24504,3.26374)(5.26504,5.28584)(5.28504,7.3197)(5.30504,9.34907)(5.32504,11.3576)(5.34504,13.3289)(5.36504,15.2466)(5.38504,17.0946)(5.40504,18.8574)(5.42504,20.5195)(5.44504,22.0666)(5.46504,23.4846)(5.48504,24.7608)(5.50504,25.8831)(5.52504,26.8407)(5.54504,27.6239)(5.56504,28.2243)(5.58504,28.635)(5.60504,28.8503)(5.62504,28.8662)(5.64504,28.6803)(5.66504,28.2916)(5.68504,27.701)(5.70504,26.9106)(5.72504,25.9247)(5.74504,24.7487)(5.76504,23.39)(5.78504,21.8573)(5.80504,20.1608)(5.82504,18.3125)(5.84504,16.3252)(5.86504,14.2135)(5.88504,11.9928)(5.90504,9.67973)(5.92504,7.29182)(5.94504,4.84741)(5.96504,2.3655)(5.98504,-0.134376)(6.00504,-2.63233)(6.02504,-5.10826)(6.04504,-7.54199)(6.06504,-9.91348)(6.08504,-12.203)(6.10504,-14.3911)(6.12504,-16.4591)(6.14504,-18.389)(6.16504,-20.1637)(6.18504,-21.7672)(6.20504,-23.1848)(6.22504,-24.4028)(6.24504,-25.4093)(6.26504,-26.1938)(6.28504,-26.7476)(6.30504,-27.0636)(6.32504,-27.1366)(6.34504,-26.9634)(6.36504,-26.5427)(6.38504,-25.8752)(6.40504,-24.9634)(6.42504,-23.8122)(6.44504,-22.4282)(6.46504,-20.82)(6.48504,-18.9981)(6.50504,-16.975)(6.52504,-14.7648)(6.54504,-12.3833)(6.56504,-9.848)(6.58504,-7.17763)(6.60504,-4.39238)(6.62504,-1.51353)(6.64504,1.43661)(6.66504,4.43493)(6.68504,7.45764)(6.70504,10.4805)(6.72504,13.4791)(6.74504,16.4287)(6.76504,19.3049)(6.78504,22.0836)(6.80504,24.7411)(6.82504,27.2546)(6.84504,29.6021)(6.86504,31.7626)(6.88504,33.7166)(6.90504,35.446)(6.92504,36.9342)(6.94504,38.1664)(6.96504,39.1299)(6.98504,39.8139)(7.00504,40.2097)(7.02504,40.311)(7.04504,40.1136)(7.06504,39.616)(7.08504,38.8188)(7.10504,37.7252)(7.12504,36.341)(7.14504,34.6741)(7.16504,32.7351)(7.18504,30.5368)(7.20504,28.0942)(7.22504,25.4245)(7.24504,22.5471)(7.26504,19.4832)(7.28504,16.2556)(7.30504,12.8889)(7.32504,9.4092)(7.34504,5.84356)(7.36504,2.22023)(7.38504,-1.4318)(7.40504,-5.08297)(7.42504,-8.70336)(7.44504,-12.263)(7.46504,-15.732)(7.48504,-19.0809)(7.50504,-22.2809)(7.52504,-25.3041)(7.54504,-28.1237)(7.56504,-30.7141)(7.58504,-33.0515)(7.60504,-35.1137)(7.62504,-36.8807)(7.64504,-38.3345)(7.66504,-39.4594)(7.68504,-40.2424)(7.70504,-40.673)(7.72504,-40.7434)(7.74504,-40.4488)(7.76504,-39.7873)(7.78504,-38.7599)(7.80504,-37.3706)(7.82504,-35.6264)(7.84504,-33.5373)(7.86504,-31.1164)(7.88504,-28.3793)(7.90504,-25.3446)(7.92504,-22.0336)(7.94504,-18.47)(7.96504,-14.6799)(7.98504,-10.6915)(8.00504,-6.53504)(8.02504,-2.24241)(8.04504,2.15295)(8.06504,6.61636)(8.08504,11.1122)(8.10504,15.6041)(8.12504,20.0553)(8.14504,24.4291)(8.16504,28.6886)(8.18504,32.7977)(8.20504,36.721)(8.22504,40.4243)(8.24504,43.8744)(8.26504,47.0403)(8.28504,49.8924)(8.30504,52.4037)(8.32504,54.5494)(8.34504,56.3075)(8.36504,57.6589)(8.38504,58.5874)(8.40504,59.0803)(8.42504,59.1282)(8.44504,58.7249)(8.46504,57.8683)(8.48504,56.5596)(8.50504,54.8038)(8.52504,52.6097)(8.54504,49.9896)(8.56504,46.9594)(8.58504,43.5387)(8.60504,39.7505)(8.62504,35.6208)(8.64504,31.1789)(8.66504,26.4571)(8.68504,21.4901)(8.70504,16.3153)(8.72504,10.9719)(8.74504,5.50124)(8.76504,-0.0538973)(8.78504,-5.64955)(8.80504,-11.2409)(8.82504,-16.7825)(8.84504,-22.2289)(8.86504,-27.5348)(8.88504,-32.6555)(8.90504,-37.5473)(8.92504,-42.1677)(8.94504,-46.4759)(8.96504,-50.4333)(8.98504,-54.0036)(9.00504,-57.153)(9.02504,-59.8509)(9.04504,-62.07)(9.06504,-63.7864)(9.08504,-64.9801)(9.10504,-65.6351)(9.12504,-65.7394)(9.14504,-65.2855)(9.16504,-64.2703)(9.18504,-62.695)(9.20504,-60.5655)(9.22504,-57.8924)(9.24504,-54.6905)(9.26504,-50.9793)(9.28504,-46.7826)(9.30504,-42.1284)(9.32504,-37.0486)(9.34504,-31.579)(9.36504,-25.7592)(9.38504,-19.6317)(9.40504,-13.2424)(9.42504,-6.63952)(9.44504,0.126234)(9.46504,7.00229)(9.48504,13.9346)(9.50504,20.8678)(9.52504,27.7463)(9.54504,34.5139)(9.56504,41.1148)(9.58504,47.4939)(9.60504,53.5972)(9.62504,59.3722)(9.64504,64.7687)(9.66504,69.7387)(9.68504,74.2372)(9.70504,78.2226)(9.72504,81.6567)(9.74504,84.5055)(9.76504,86.7393)(9.78504,88.333)(9.80504,89.2664)(9.82504,89.5243)(9.84504,89.0969)(9.86504,87.9799)(9.88504,86.1743)(9.90504,83.6868)(9.92504,80.5298)(9.94504,76.721)(9.96504,72.2837)(9.98504,67.2467)(10,63.1072) 
};

\addplot [
color=mycolor2,
solid
]
coordinates{
 (0,4)(0.000189186,4)(0.00122971,4.00002)(0.00312157,4.00012)(0.0118109,4.00175)(0.0272976,4.00938)(0.045041,4.0256)(0.065041,4.05341)(0.085041,4.09116)(0.105041,4.13857)(0.125041,4.1952)(0.145041,4.26048)(0.165041,4.33374)(0.185041,4.41417)(0.205041,4.50084)(0.225041,4.59272)(0.245041,4.68868)(0.265041,4.7875)(0.285041,4.88788)(0.305041,4.98845)(0.325041,5.08781)(0.345041,5.18448)(0.365041,5.277)(0.385041,5.36387)(0.405041,5.44359)(0.425041,5.51471)(0.445041,5.57578)(0.465041,5.62544)(0.485041,5.66238)(0.505041,5.68535)(0.525041,5.69323)(0.545041,5.68501)(0.565041,5.65978)(0.585041,5.61678)(0.605041,5.55542)(0.625041,5.47524)(0.645041,5.37596)(0.665041,5.25747)(0.685041,5.11985)(0.705041,4.96337)(0.725041,4.78848)(0.745041,4.59582)(0.765041,4.38623)(0.785041,4.16074)(0.805041,3.92055)(0.825041,3.66704)(0.845041,3.40179)(0.865041,3.1265)(0.885041,2.84304)(0.905041,2.55343)(0.925041,2.25978)(0.945041,1.96434)(0.965041,1.66943)(0.985041,1.37745)(1.00504,1.09084)(1.02504,0.812095)(1.04504,0.543697)(1.06504,0.288126)(1.08504,0.0478254)(1.10504,-0.174818)(1.12504,-0.377497)(1.14504,-0.558008)(1.16504,-0.714272)(1.18504,-0.844354)(1.20504,-0.94649)(1.22504,-1.0191)(1.24504,-1.06082)(1.26504,-1.07049)(1.28504,-1.04721)(1.30504,-0.99033)(1.32504,-0.899457)(1.34504,-0.774487)(1.36504,-0.615595)(1.38504,-0.423253)(1.40504,-0.198225)(1.42504,0.0584251)(1.44504,0.345338)(1.46504,0.660863)(1.48504,1.00306)(1.50504,1.36972)(1.52504,1.75835)(1.54504,2.16623)(1.56504,2.59038)(1.58504,3.02762)(1.60504,3.47457)(1.62504,3.92768)(1.64504,4.38324)(1.66504,4.83744)(1.68504,5.28637)(1.70504,5.72604)(1.72504,6.15247)(1.74504,6.56165)(1.76504,6.94961)(1.78504,7.31247)(1.80504,7.64644)(1.82504,7.94786)(1.84504,8.21325)(1.86504,8.43934)(1.88504,8.62308)(1.90504,8.7617)(1.92504,8.85271)(1.94504,8.89397)(1.96504,8.88366)(1.98504,8.82035)(2.00504,8.703)(2.02504,8.53097)(2.04504,8.30407)(2.06504,8.02252)(2.08504,7.68701)(2.10504,7.29867)(2.12504,6.85909)(2.14504,6.3703)(2.16504,5.8348)(2.18504,5.2555)(2.20504,4.63576)(2.22504,3.97934)(2.24504,3.29037)(2.26504,2.57338)(2.28504,1.83321)(2.30504,1.075)(2.32504,0.304198)(2.34504,-0.473547)(2.36504,-1.25237)(2.38504,-2.02626)(2.40504,-2.78908)(2.42504,-3.53464)(2.44504,-4.25673)(2.46504,-4.94917)(2.48504,-5.60583)(2.50504,-6.22075)(2.52504,-6.78811)(2.54504,-7.30235)(2.56504,-7.75815)(2.58504,-8.15052)(2.60504,-8.47486)(2.62504,-8.72695)(2.64504,-8.90303)(2.66504,-8.99984)(2.68504,-9.01463)(2.70504,-8.94522)(2.72504,-8.79001)(2.74504,-8.54802)(2.76504,-8.21888)(2.78504,-7.8029)(2.80504,-7.30101)(2.82504,-6.71483)(2.84504,-6.04664)(2.86504,-5.29937)(2.88504,-4.47663)(2.90504,-3.58262)(2.92504,-2.62221)(2.94504,-1.60083)(2.96504,-0.524487)(2.98504,0.600287)(3.00504,1.76647)(3.02504,2.96659)(3.04504,4.19278)(3.06504,5.43684)(3.08504,6.69025)(3.10504,7.94429)(3.12504,9.19005)(3.14504,10.4185)(3.16504,11.6207)(3.18504,12.7874)(3.20504,13.9099)(3.22504,14.9794)(3.24504,15.9873)(3.26504,16.9254)(3.28504,17.7858)(3.30504,18.5612)(3.32504,19.2447)(3.34504,19.8298)(3.36504,20.3108)(3.38504,20.6828)(3.40504,20.9412)(3.42504,21.0825)(3.44504,21.104)(3.46504,21.0037)(3.48504,20.7806)(3.50504,20.4344)(3.52504,19.9659)(3.54504,19.3767)(3.56504,18.6693)(3.58504,17.8473)(3.60504,16.915)(3.62504,15.8777)(3.64504,14.7415)(3.66504,13.5133)(3.68504,12.2008)(3.70504,10.8126)(3.72504,9.35771)(3.74504,7.84602)(3.76504,6.28786)(3.78504,4.69408)(3.80504,3.07598)(3.82504,1.44519)(3.84504,-0.186409)(3.86504,-1.80673)(3.88504,-3.4036)(3.90504,-4.96486)(3.92504,-6.47841)(3.94504,-7.93236)(3.96504,-9.31506)(3.98504,-10.6153)(4.00504,-11.8221)(4.02504,-12.9254)(4.04504,-13.9154)(4.06504,-14.7832)(4.08504,-15.5207)(4.10504,-16.1204)(4.12504,-16.5762)(4.14504,-16.8825)(4.16504,-17.0349)(4.18504,-17.0303)(4.20504,-16.8664)(4.22504,-16.5423)(4.24504,-16.0581)(4.26504,-15.4152)(4.28504,-14.6162)(4.30504,-13.6648)(4.32504,-12.5661)(4.34504,-11.3262)(4.36504,-9.95222)(4.38504,-8.45264)(4.40504,-6.83685)(4.42504,-5.11523)(4.44504,-3.2991)(4.46504,-1.40067)(4.48504,0.567132)(4.50504,2.59065)(4.52504,4.65563)(4.54504,6.74735)(4.56504,8.85067)(4.58504,10.9502)(4.60504,13.0303)(4.62504,15.0753)(4.64504,17.0696)(4.66504,18.9979)(4.68504,20.8449)(4.70504,22.5959)(4.72504,24.2367)(4.74504,25.7538)(4.76504,27.1343)(4.78504,28.3662)(4.80504,29.4384)(4.82504,30.341)(4.84504,31.065)(4.86504,31.6029)(4.88504,31.9482)(4.90504,32.0958)(4.92504,32.0421)(4.94504,31.7849)(4.96504,31.3234)(4.98504,30.6583)(5.00504,29.7919)(5.02504,28.728)(5.04504,27.4718)(5.06504,26.03)(5.08504,24.4108)(5.10504,22.6237)(5.12504,20.6796)(5.14504,18.5907)(5.16504,16.3703)(5.18504,14.0329)(5.20504,11.594)(5.22504,9.06992)(5.24504,6.47791)(5.26504,3.83582)(5.28504,1.16207)(5.30504,-1.52451)(5.32504,-4.20478)(5.34504,-6.85948)(5.36504,-9.46933)(5.38504,-12.0152)(5.40504,-14.4782)(5.42504,-16.8399)(5.44504,-19.0825)(5.46504,-21.1888)(5.48504,-23.1424)(5.50504,-24.928)(5.52504,-26.5314)(5.54504,-27.9395)(5.56504,-29.1405)(5.58504,-30.1242)(5.60504,-30.8816)(5.62504,-31.4056)(5.64504,-31.6906)(5.66504,-31.7327)(5.68504,-31.5299)(5.70504,-31.0816)(5.72504,-30.3895)(5.74504,-29.4569)(5.76504,-28.2886)(5.78504,-26.8917)(5.80504,-25.2748)(5.82504,-23.448)(5.84504,-21.4233)(5.86504,-19.2141)(5.88504,-16.8354)(5.90504,-14.3035)(5.92504,-11.6359)(5.94504,-8.85118)(5.96504,-5.96914)(5.98504,-3.01026)(6.00504,0.00421046)(6.02504,3.05246)(6.04504,6.11228)(6.06504,9.16122)(6.08504,12.1768)(6.10504,15.1365)(6.12504,18.0182)(6.14504,20.8002)(6.16504,23.4615)(6.18504,25.9815)(6.20504,28.3409)(6.22504,30.5212)(6.24504,32.5055)(6.26504,34.2778)(6.28504,35.8239)(6.30504,37.1311)(6.32504,38.1883)(6.34504,38.9865)(6.36504,39.5183)(6.38504,39.7783)(6.40504,39.7633)(6.42504,39.4719)(6.44504,38.9049)(6.46504,38.0652)(6.48504,36.9576)(6.50504,35.5892)(6.52504,33.9689)(6.54504,32.1076)(6.56504,30.0181)(6.58504,27.715)(6.60504,25.2145)(6.62504,22.5344)(6.64504,19.6942)(6.66504,16.7143)(6.68504,13.6166)(6.70504,10.4238)(6.72504,7.15955)(6.74504,3.84811)(6.76504,0.514259)(6.78504,-2.81696)(6.80504,-6.12037)(6.82504,-9.37091)(6.84504,-12.5438)(6.86504,-15.6146)(6.88504,-18.5597)(6.90504,-21.3562)(6.92504,-23.9824)(6.94504,-26.4175)(6.96504,-28.6423)(6.98504,-30.639)(7.00504,-32.3915)(7.02504,-33.8854)(7.04504,-35.1082)(7.06504,-36.0495)(7.08504,-36.7009)(7.10504,-37.0562)(7.12504,-37.1116)(7.14504,-36.8651)(7.16504,-36.3175)(7.18504,-35.4715)(7.20504,-34.3324)(7.22504,-32.9074)(7.24504,-31.2063)(7.26504,-29.2407)(7.28504,-27.0243)(7.30504,-24.5729)(7.32504,-21.9041)(7.34504,-19.0372)(7.36504,-15.9929)(7.38504,-12.7935)(7.40504,-9.46261)(7.42504,-6.02466)(7.44504,-2.50513)(7.46504,1.06985)(7.48504,4.67368)(7.50504,8.27948)(7.52504,11.8603)(7.54504,15.3893)(7.56504,18.84)(7.58504,22.1866)(7.60504,25.4038)(7.62504,28.4675)(7.64504,31.3546)(7.66504,34.0435)(7.68504,36.5139)(7.70504,38.7473)(7.72504,40.7269)(7.74504,42.4381)(7.76504,43.868)(7.78504,45.0063)(7.80504,45.8446)(7.82504,46.3769)(7.84504,46.5999)(7.86504,46.5123)(7.88504,46.1155)(7.90504,45.4131)(7.92504,44.4112)(7.94504,43.1185)(7.96504,41.5455)(7.98504,39.7052)(8.00504,37.6127)(8.02504,35.285)(8.04504,32.741)(8.06504,30.0012)(8.08504,27.0878)(8.10504,24.0241)(8.12504,20.8349)(8.14504,17.5457)(8.16504,14.1829)(8.18504,10.7733)(8.20504,7.34428)(8.22504,3.92307)(8.24504,0.536922)(8.26504,-2.78728)(8.28504,-6.0232)(8.30504,-9.14529)(8.32504,-12.129)(8.34504,-14.9509)(8.36504,-17.5891)(8.38504,-20.023)(8.40504,-22.234)(8.42504,-24.2053)(8.44504,-25.9221)(8.46504,-27.3718)(8.48504,-28.544)(8.50504,-29.4307)(8.52504,-30.0262)(8.54504,-30.3275)(8.56504,-30.3339)(8.58504,-30.0472)(8.60504,-29.4717)(8.62504,-28.6144)(8.64504,-27.4842)(8.66504,-26.0926)(8.68504,-24.4535)(8.70504,-22.5825)(8.72504,-20.4974)(8.74504,-18.2178)(8.76504,-15.7648)(8.78504,-13.1613)(8.80504,-10.431)(8.82504,-7.59911)(8.84504,-4.69136)(8.86504,-1.73426)(8.88504,1.24533)(8.90504,4.22039)(8.92504,7.16397)(8.94504,10.0495)(8.96504,12.8508)(8.98504,15.5427)(9.00504,18.1009)(9.02504,20.5024)(9.04504,22.7255)(9.06504,24.7503)(9.08504,26.5585)(9.10504,28.1341)(9.12504,29.463)(9.14504,30.5333)(9.16504,31.3357)(9.18504,31.8632)(9.20504,32.1113)(9.22504,32.0783)(9.24504,31.7649)(9.26504,31.1745)(9.28504,30.3131)(9.30504,29.1892)(9.32504,27.8139)(9.34504,26.2004)(9.36504,24.3646)(9.38504,22.3242)(9.40504,20.0991)(9.42504,17.7108)(9.44504,15.1826)(9.46504,12.5392)(9.48504,9.80649)(9.50504,7.01131)(9.52504,4.18122)(9.54504,1.34427)(9.56504,-1.4713)(9.58504,-4.23728)(9.60504,-6.92578)(9.62504,-9.5095)(9.64504,-11.962)(9.66504,-14.2578)(9.68504,-16.373)(9.70504,-18.2851)(9.72504,-19.9734)(9.74504,-21.4194)(9.76504,-22.6066)(9.78504,-23.5211)(9.80504,-24.1512)(9.82504,-24.4883)(9.84504,-24.5264)(9.86504,-24.2624)(9.88504,-23.696)(9.90504,-22.8301)(9.92504,-21.6705)(9.94504,-20.2259)(9.96504,-18.5081)(9.98504,-16.5316)(10,-14.8944) 
};

\addplot [
color=mycolor3,
solid
]
coordinates{
 (0,5)(0.000189186,5)(0.00122971,5)(0.00312157,4.99998)(0.0118109,4.99965)(0.0272976,4.99812)(0.045041,4.99488)(0.065041,4.98932)(0.085041,4.98178)(0.105041,4.97231)(0.125041,4.96101)(0.145041,4.94799)(0.165041,4.93338)(0.185041,4.91736)(0.205041,4.9001)(0.225041,4.88182)(0.245041,4.86272)(0.265041,4.84305)(0.285041,4.82305)(0.305041,4.80297)(0.325041,4.78307)(0.345041,4.76361)(0.365041,4.74485)(0.385041,4.72704)(0.405041,4.7104)(0.425041,4.69518)(0.445041,4.68157)(0.465041,4.66976)(0.485041,4.6599)(0.505041,4.65213)(0.525041,4.64653)(0.545041,4.64318)(0.565041,4.64209)(0.585041,4.64325)(0.605041,4.64658)(0.625041,4.65199)(0.645041,4.65933)(0.665041,4.66838)(0.685041,4.67891)(0.705041,4.69062)(0.725041,4.70318)(0.745041,4.71619)(0.765041,4.72923)(0.785041,4.74182)(0.805041,4.75347)(0.825041,4.76362)(0.845041,4.7717)(0.865041,4.77712)(0.885041,4.77924)(0.905041,4.77742)(0.925041,4.77102)(0.945041,4.75938)(0.965041,4.74183)(0.985041,4.71773)(1.00504,4.68644)(1.02504,4.64735)(1.04504,4.59987)(1.06504,4.54345)(1.08504,4.47757)(1.10504,4.40178)(1.12504,4.31567)(1.14504,4.21889)(1.16504,4.11117)(1.18504,3.9923)(1.20504,3.86217)(1.22504,3.72073)(1.24504,3.56805)(1.26504,3.40425)(1.28504,3.22958)(1.30504,3.04438)(1.32504,2.84909)(1.34504,2.64425)(1.36504,2.43049)(1.38504,2.20857)(1.40504,1.97931)(1.42504,1.74365)(1.44504,1.50263)(1.46504,1.25736)(1.48504,1.00903)(1.50504,0.758937)(1.52504,0.508416)(1.54504,0.258879)(1.56504,0.01179)(1.58504,-0.231342)(1.60504,-0.468973)(1.62504,-0.699536)(1.64504,-0.92145)(1.66504,-1.13313)(1.68504,-1.33301)(1.70504,-1.51954)(1.72504,-1.69121)(1.74504,-1.84654)(1.76504,-1.98413)(1.78504,-2.10266)(1.80504,-2.20086)(1.82504,-2.2776)(1.84504,-2.33182)(1.86504,-2.36262)(1.88504,-2.36919)(1.90504,-2.35089)(1.92504,-2.30721)(1.94504,-2.23783)(1.96504,-2.14254)(1.98504,-2.02136)(2.00504,-1.87445)(2.02504,-1.70216)(2.04504,-1.50502)(2.06504,-1.28375)(2.08504,-1.03924)(2.10504,-0.772569)(2.12504,-0.485004)(2.14504,-0.177973)(2.16504,0.146923)(2.18504,0.487924)(2.20504,0.843116)(2.22504,1.21045)(2.24504,1.58773)(2.26504,1.97266)(2.28504,2.36283)(2.30504,2.75576)(2.32504,3.14887)(2.34504,3.53954)(2.36504,3.92512)(2.38504,4.30293)(2.40504,4.67028)(2.42504,5.02454)(2.44504,5.36308)(2.46504,5.68334)(2.48504,5.98286)(2.50504,6.25927)(2.52504,6.5103)(2.54504,6.73385)(2.56504,6.92798)(2.58504,7.09089)(2.60504,7.22101)(2.62504,7.31696)(2.64504,7.3776)(2.66504,7.40201)(2.68504,7.38954)(2.70504,7.33978)(2.72504,7.25259)(2.74504,7.12812)(2.76504,6.96679)(2.78504,6.76931)(2.80504,6.53665)(2.82504,6.27009)(2.84504,5.97119)(2.86504,5.64174)(2.88504,5.28386)(2.90504,4.89987)(2.92504,4.49236)(2.94504,4.06415)(2.96504,3.61827)(2.98504,3.15793)(3.00504,2.68653)(3.02504,2.20763)(3.04504,1.72492)(3.06504,1.24218)(3.08504,0.763273)(3.10504,0.292137)(3.12504,-0.167283)(3.14504,-0.61104)(3.16504,-1.03522)(3.18504,-1.43596)(3.20504,-1.80951)(3.22504,-2.15221)(3.24504,-2.46058)(3.26504,-2.7313)(3.28504,-2.96125)(3.30504,-3.14756)(3.32504,-3.28762)(3.34504,-3.37907)(3.36504,-3.4199)(3.38504,-3.40838)(3.40504,-3.34317)(3.42504,-3.22325)(3.44504,-3.04801)(3.46504,-2.81719)(3.48504,-2.53097)(3.50504,-2.18988)(3.52504,-1.79491)(3.54504,-1.34742)(3.56504,-0.849177)(3.58504,-0.302361)(3.60504,0.290466)(3.62504,0.926364)(3.64504,1.60203)(3.66504,2.3138)(3.68504,3.05771)(3.70504,3.82947)(3.72504,4.62451)(3.74504,5.43802)(3.76504,6.26497)(3.78504,7.10013)(3.80504,7.93813)(3.82504,8.77347)(3.84504,9.60057)(3.86504,10.4138)(3.88504,11.2076)(3.90504,11.9762)(3.92504,12.7143)(3.94504,13.4165)(3.96504,14.0774)(3.98504,14.692)(4.00504,15.2556)(4.02504,15.7637)(4.04504,16.2119)(4.06504,16.5964)(4.08504,16.9138)(4.10504,17.1609)(4.12504,17.3351)(4.14504,17.4342)(4.16504,17.4565)(4.18504,17.4008)(4.20504,17.2663)(4.22504,17.053)(4.24504,16.7613)(4.26504,16.3919)(4.28504,15.9465)(4.30504,15.427)(4.32504,14.8359)(4.34504,14.1763)(4.36504,13.4519)(4.38504,12.6665)(4.40504,11.8249)(4.42504,10.9318)(4.44504,9.99277)(4.46504,9.01348)(4.48504,8.00003)(4.50504,6.95881)(4.52504,5.89647)(4.54504,4.81987)(4.56504,3.73605)(4.58504,2.65217)(4.60504,1.57543)(4.62504,0.513097)(4.64504,-0.527628)(4.66504,-1.53962)(4.68504,-2.51587)(4.70504,-3.44959)(4.72504,-4.33421)(4.74504,-5.16344)(4.76504,-5.93134)(4.78504,-6.63235)(4.80504,-7.26133)(4.82504,-7.8136)(4.84504,-8.285)(4.86504,-8.67191)(4.88504,-8.97126)(4.90504,-9.18062)(4.92504,-9.29814)(4.94504,-9.32264)(4.96504,-9.25359)(4.98504,-9.09113)(5.00504,-8.83607)(5.02504,-8.48989)(5.04504,-8.05476)(5.06504,-7.53349)(5.08504,-6.92955)(5.10504,-6.24702)(5.12504,-5.49061)(5.14504,-4.66558)(5.16504,-3.77776)(5.18504,-2.83346)(5.20504,-1.83947)(5.22504,-0.802998)(5.24504,0.268385)(5.26504,1.36678)(5.28504,2.48403)(5.30504,3.61177)(5.32504,4.74148)(5.34504,5.86456)(5.36504,6.97239)(5.38504,8.05636)(5.40504,9.108)(5.42504,10.119)(5.44504,11.0812)(5.46504,11.9868)(5.48504,12.8284)(5.50504,13.5989)(5.52504,14.2916)(5.54504,14.9006)(5.56504,15.4203)(5.58504,15.8457)(5.60504,16.1725)(5.62504,16.3973)(5.64504,16.5171)(5.66504,16.5298)(5.68504,16.4341)(5.70504,16.2293)(5.72504,15.9158)(5.74504,15.4944)(5.76504,14.9671)(5.78504,14.3364)(5.80504,13.6058)(5.82504,12.7795)(5.84504,11.8622)(5.86504,10.8597)(5.88504,9.77811)(5.90504,8.62449)(5.92504,7.4063)(5.94504,6.13157)(5.96504,4.80886)(5.98504,3.4471)(6.00504,2.05562)(6.02504,0.644042)(6.04504,-0.777777)(6.06504,-2.19982)(6.08504,-3.61199)(6.10504,-5.00415)(6.12504,-6.36623)(6.14504,-7.6883)(6.16504,-8.96061)(6.18504,-10.1737)(6.20504,-11.3185)(6.22504,-12.3862)(6.24504,-13.3686)(6.26504,-14.2581)(6.28504,-15.0476)(6.30504,-15.7306)(6.32504,-16.3014)(6.34504,-16.7551)(6.36504,-17.0875)(6.38504,-17.2953)(6.40504,-17.376)(6.42504,-17.328)(6.44504,-17.1507)(6.46504,-16.8441)(6.48504,-16.4096)(6.50504,-15.8491)(6.52504,-15.1655)(6.54504,-14.3628)(6.56504,-13.4455)(6.58504,-12.4193)(6.60504,-11.2903)(6.62504,-10.0658)(6.64504,-8.75354)(6.66504,-7.36187)(6.68504,-5.89985)(6.70504,-4.37705)(6.72504,-2.8035)(6.74504,-1.18963)(6.76504,0.453835)(6.78504,2.11589)(6.80504,3.7854)(6.82504,5.45111)(6.84504,7.10179)(6.86504,8.72628)(6.88504,10.3136)(6.90504,11.8529)(6.92504,13.3339)(6.94504,14.7464)(6.96504,16.0809)(6.98504,17.3283)(7.00504,18.4803)(7.02504,19.5291)(7.04504,20.4676)(7.06504,21.2897)(7.08504,21.9899)(7.10504,22.5637)(7.12504,23.0075)(7.14504,23.3186)(7.16504,23.4952)(7.18504,23.5367)(7.20504,23.4431)(7.22504,23.2157)(7.24504,22.8567)(7.26504,22.3691)(7.28504,21.757)(7.30504,21.0254)(7.32504,20.1801)(7.34504,19.2277)(7.36504,18.1756)(7.38504,17.032)(7.40504,15.8056)(7.42504,14.5059)(7.44504,13.1427)(7.46504,11.7264)(7.48504,10.2677)(7.50504,8.77777)(7.52504,7.26773)(7.54504,5.74897)(7.56504,4.23292)(7.58504,2.73092)(7.60504,1.25422)(7.62504,-0.186172)(7.64504,-1.57952)(7.66504,-2.91548)(7.68504,-4.18417)(7.70504,-5.3762)(7.72504,-6.48282)(7.74504,-7.49593)(7.76504,-8.40813)(7.78504,-9.21284)(7.80504,-9.90427)(7.82504,-10.4775)(7.84504,-10.9286)(7.86504,-11.2544)(7.88504,-11.4528)(7.90504,-11.5227)(7.92504,-11.464)(7.94504,-11.2773)(7.96504,-10.9645)(7.98504,-10.5282)(8.00504,-9.97217)(8.02504,-9.30086)(8.04504,-8.51965)(8.06504,-7.63474)(8.08504,-6.65308)(8.10504,-5.58231)(8.12504,-4.43075)(8.14504,-3.20725)(8.16504,-1.92119)(8.18504,-0.582382)(8.20504,0.799025)(8.22504,2.2126)(8.24504,3.64774)(8.26504,5.0937)(8.28504,6.53977)(8.30504,7.97524)(8.32504,9.38961)(8.34504,10.7725)(8.36504,12.1141)(8.38504,13.4046)(8.40504,14.6349)(8.42504,15.7964)(8.44504,16.881)(8.46504,17.8814)(8.48504,18.7909)(8.50504,19.6035)(8.52504,20.3141)(8.54504,20.9184)(8.56504,21.413)(8.58504,21.7952)(8.60504,22.0635)(8.62504,22.2169)(8.64504,22.2557)(8.66504,22.1807)(8.68504,21.9939)(8.70504,21.698)(8.72504,21.2964)(8.74504,20.7936)(8.76504,20.1944)(8.78504,19.5046)(8.80504,18.7306)(8.82504,17.8793)(8.84504,16.9581)(8.86504,15.9747)(8.88504,14.9375)(8.90504,13.8549)(8.92504,12.7357)(8.94504,11.5885)(8.96504,10.4225)(8.98504,9.24638)(9.00504,8.06896)(9.02504,6.89883)(9.04504,5.74432)(9.06504,4.61347)(9.08504,3.51388)(9.10504,2.45274)(9.12504,1.43671)(9.14504,0.471871)(9.16504,-0.436297)(9.18504,-1.28298)(9.20504,-2.06404)(9.22504,-2.77609)(9.24504,-3.41648)(9.26504,-3.98329)(9.28504,-4.47542)(9.30504,-4.89248)(9.32504,-5.23489)(9.34504,-5.50378)(9.36504,-5.70101)(9.38504,-5.82915)(9.40504,-5.89141)(9.42504,-5.89163)(9.44504,-5.83419)(9.46504,-5.72401)(9.48504,-5.56646)(9.50504,-5.36729)(9.52504,-5.13258)(9.54504,-4.86866)(9.56504,-4.58204)(9.58504,-4.27935)(9.60504,-3.96722)(9.62504,-3.65225)(9.64504,-3.34092)(9.66504,-3.03949)(9.68504,-2.75397)(9.70504,-2.49002)(9.72504,-2.25288)(9.74504,-2.04732)(9.76504,-1.8776)(9.78504,-1.74736)(9.80504,-1.65963)(9.82504,-1.61676)(9.84504,-1.62039)(9.86504,-1.67143)(9.88504,-1.77004)(9.90504,-1.91562)(9.92504,-2.10678)(9.94504,-2.34139)(9.96504,-2.61656)(9.98504,-2.92869)(10,-3.18375) 
};

\end{axis}
\end{tikzpicture}